
\input amstex
\input xy
\xyoption{all}
\documentstyle{amsppt}
\document
\magnification=1200
\NoBlackBoxes
\nologo
\vsize16cm

\def\Id{{Id}}
\def\Z{{\Bbb Z}}
\def\L{{\Bbb L}}
\def\colim{{colim}}

\def\t{{\tau}}
\def\twomor{{\text{ }\rightarrow\text{ }}}
\def\monpro{{\mu}}

\def\ass{{\alpha}}

\def\cohom{{\underline{cohom}}} 
\def\hom{{Hom}}
\def\tensor{{\otimes}} 
\def\white{{\bigcirc}}
\def\triple{{T}}
\def\U{{\frak U}}
\def\operad{{P}}
\def\Operad{{\Cal{OP}}} 

\def\C{{\Cal C}} 
\def\D{{\Cal D}} 
\def\A{{\Cal A}}
\def\As{{\frak A}} 
\def\as{{A}} 

\def\Cat{{\Cal{CAT}}} 
\def\E{{E}} 
\def\Endo{{\frak E}} 
\def\Fun{{Fun}} 
\def\OpSt{{\Cal{OP}_{st}}} 

\def\Tot{{Tot}}
\def\colfun{{\chi}}
\def\D{{\Cal D}}

\def\PsOpSt{{\Psi\Cal{OP}_{st}}} 
\def\PsOpLC{{\Psi\Cal{OP}}} 
\def\opcomp{{\Upsilon}} 
\def\operad{{P}} 
\def\opmorf{{F}} 
\def\opmort{{\zeta}} 
\def\ntcomp{{*}} 
\def\opmoro{{\alpha}} 
\def\PsOpP{{\Psi\Cal{OP}(\Cat)/ /\operad}} 
\def\PsOpPS{{\Psi\Cal{OP}_{st}(\Cat)/ /\operad}} 
\def\opunit{{U}} 
\def\opuntr{{\eta}} 
\def\OpP{{\Operad(\Cat)/ /\operad}} 
\def\OpPS{{\Operad_{st}(\Cat)/ /\operad}} 
\def\opelem{{G}} 
\def\opunel{{e}} 
\def\opmorff{{G}}
\def\opmortt{{\psi}}
\def\opuntrr{{\phi}}
\def\GenTri{{\frak T}}
\def\OpTrip{{\frak P}}
\def\Triple{{\Cal T}}
\def\c{{C}}
\def\Colors{{\Omega}}
\def\color{{\omega}}
\def\Graphs{{\frak G}}
\def\colorm{{f}}
\def\Sym{{\Sigma}}
\def\sym{{\sigma}}
\def\oper{{p}}
\def\SYM{{\underline{\Sym}}}
\def\symact{{S}}
\def\OpTriS{{\OpTrip^s}}
\def\corolla{{\sigma}}
\def\cormor{{\phi}}
\def\graph{{\tau}}
\def\Vertices{{V}}
\def\vertex{{v}}
\def\grarep{{\rho}}
\def\coldia{{F}}
\def\hopf{{\tau}}

\bigskip

\centerline{\bf GENERALIZED OPERADS}

\smallskip

\centerline{\bf AND THEIR INNER COHOMOMORPHISMS\footnotemark}
\footnotetext{{\it 2000 Mathematics Subject Classification:} 18D50, 18D10, 20C30.
{\it Keywords and phrases:} Operads, algebras,  inner cohomomorphisms,
symmetry and deformations in noncommutative geometry.}
\medskip

\centerline{\bf Dennis V. Borisov, Yuri I. Manin}

\medskip

\centerline{\it Northwestern University, Evanston, USA}

\centerline{\it Max--Planck--Institut f\"ur Mathematik, Bonn, Germany,}

\bigskip

{\bf Abstract.} In this paper we introduce  a notion
of {\it generalized operad} containing as special cases
various kinds of operad--like objects: ordinary, cyclic, modular,
pro\-perads etc.   We then construct
inner cohomomorphism objects in their categories (and 
categories of algebras over them).  We argue that 
they provide an approach to
symmetry and moduli objects in non-commutative geometries
based upon these  ``ring--like'' structures.
We give a unified axiomatic treatment of generalized operads
as functors on categories of abstract labeled graphs. 
Finally, we extend 
inner cohomomorphism constructions to more
general categorical contexts.
This version differs from the previous ones by several local
changes (including the title) and two extra references.

\bigskip

\centerline{\bf \S 0. Introduction}

\medskip

{\bf 0.1.  Inner cohomomorphisms of associative algebras.}
Let $k$ be a field. Consider pairs $\Cal{A}=(A,A_1)$
consisting of an associative $k$--algebra $A$ and a
finite dimensional subspace $A_1$ generating $A$. For two
such pairs $\Cal{A}=(A,A_1)$ and $\Cal{B}=(B,B_1)$,
define the category $\Cal{A}\Rightarrow \Cal{B}$
by the following data.

\smallskip

An object of $\Cal{A}\Rightarrow \Cal{B}$ is a pair
$(F,u)$ where $F$ is a $k$--algebra and $u:\,A\to F\otimes B$
is a homomorphism of algebras such that $u(A_1)\subset F\otimes B_1$
(all tensor products being taken over $k$).

\smallskip

A morphism $(F,u)\to (F^{\prime},u^{\prime})$ in $\Cal{A}\Rightarrow \Cal{B}$
is a homomorphism of algebras $v:\,F\to F^{\prime}$
such that $u^{\prime} = (v\otimes \roman{id}_B)\circ u$.

\smallskip

The following result was proved in [Ma3] (see. Prop. 2.3 in Chapter 4):

\medskip

{\bf 0.1.1. Theorem.} {\it The category $\Cal{A}\Rightarrow \Cal{B}$
has an initial object
$$
(E,\,\delta :\, A\to E\otimes B)
$$
defined uniquely up to unique isomorphism,
together with a finite dimensional subspace $E_1\subset E$
generating $E$ and satisfying $\delta (A_1)\subset E_1\otimes B_1$.}

\medskip

This result can be reinterpreted as follows. Consider
another category $PAlg$ whose objects are finitely generated
$k$--algebras {\it together with a presentation $P$}, i.~e. a surjection
$\varphi_A:\,T(A_1)\to A$ where $A_1$ is a finite dimensional linear space,
$T(A_1)$ is its tensor algebra, and such that $\roman{Ker}\,\varphi_A
\cap A_1 =\{0\}$ so that $A_1$ can be considered as a subspace
of $A$ (this condition is not really necessary and
can be omitted as is done in Sec. 2.)
Morphisms $(A,\varphi_A)\to (B,\varphi_B)$
are algebra homomorphisms $u:\,A\to B$ such that $u(A_1)\subset B_1.$

\smallskip

This category has a monoidal symmetric structure
given by $(A,A_1)\bigcirc (B,B_1) :=(C,C_1)$ where $C_1=A_1\otimes B_1$
and $C$ is the subalgebra of $A\otimes B$ generated by $C_1$.

\smallskip

Theorem 0.1.1 establishes a functorial bijection between $\roman{Hom}$'s
in this category
$$
\roman{Hom}\,((A,A_1), (F,F_1)\bigcirc (B,B_1)) \cong
\roman{Hom}\,((E,E_1), (F,F_1)).
$$
When there is no risk of confusion, we will omit $A_1, B_1$ etc
in notation, and denote $(E,E_1)$ by $\underline{cohom}\,(A,B)$
so that we have the standard functorial isomorphism in $(PAlg,\bigcirc )$:
$$
\roman{Hom}\,(A, F\bigcirc B) \cong
\roman{Hom}\,(\underline{cohom}\,(A,B), F)
$$
defining  inner cohomomorphism objects.

\smallskip

The usual reasoning produces functorial
comultiplication maps between these objects
$$
\Delta_{A,B,C}:\ \underline{cohom}\,(A,C)\to
\underline{cohom}\,(A,B)\bigcirc \underline{cohom}\,(B,C),
\eqno(0.1)
$$
which are coassociative (compatible with the ordinary
associativity constraints for $\bigcirc$).

\smallskip

In particular $\underline{coend}\,(A):=\underline{cohom}\,(A,A)$
has the canonical structure of a bialgebra.

\medskip

{\bf 0.2. Interpretation and motivation.}
Theorem 0.1.1 was the base of the approach to
quantum groups as symmetry objects in noncommutative
geometry discussed in [Ma1]--[Ma4]. Namely, consider
$PAlg$ as a category of function algebras on
``quantum linear spaces'' so that the category of quantum linear
spaces themselves will be $PAlg^{op}.$ Then
cohomomorphism algebras correspond to ``matrix quantum spaces'',
and coendomorphism algeras, after passing to Hopf envelopes,
become Hopf algebras of symmetries. (In fact, to obtain the conventional
quantum groups, one has to add some ``missing relations'', cf. [Ma2],
which also can be done functorially).

\smallskip

In this paper we present several layers of generalizations of
Theorem 0.1.1. The first step consists in extending it to operads
with presentation and algebras over them, with an appropriate
monoidal structure. We are motivated by the same desire to
understand symmetry objects (``quantum semi--groups'') in
non--commutative geometry based upon operads, or algebras over an
operad different from $ASS$. In fact, as a bonus we also get an
unconventional approach to the deformation theory of operadic
algebras.

\smallskip

Namely, let $\Cal{P}$ be an operad, $V$ a linear space,
and $OpEnd\,V$ the operad of endomorphisms of $V$.
The set of structures of a $\Cal{P}$--algebra upon $V$
is then
$$
\roman{Hom}_{Oper}\,(\Cal{P}, OpEnd\,V).
\eqno(0.2)
$$
We suggest to consider the object (defined after a choice
of spaces of linear generators of both operads)
$$
\underline{cohom}\,(\Cal{P}, OpEnd\,V)
\eqno(0.3)
$$
as (an operad of functions upon) a noncommutative space.
Morphisms of (0.3) to the unit object of the
monoidal category of operads will then constitute
its set of ``classical points'' (0.2).

\medskip

{\bf 0.3. The phantom of the operad.} The next extension of
Theorem 0.1.1 involves replacing operads by any of the related
structures a representative list of which the reader can find, for
example, in [Mar]: May and Markl operads, cyclic operads, modular
operads, PROPS, properads, dioperads etc. In this paper, we use
for all of them the generic name ``generalized operad'',
or simply ``operad'', and call operads like
May's and Markl's ones ``ordinary operads'', or ``classical
operads''.

\smallskip

It was long recognized that
one variable part of the definition of all these structures
is the combinatorics and decoration of underlying
graphs (``pasting schemes''
of [Mar]), whereas another is the category in which
components of the respective operad are supposed to lie.
Operad itself for us is a functor from a category
of labeled graphs to another symmetric monoidal category,
as was stressed already in [KoMa], [GeKa2] and many other
works. We decided to spell out the underlying formalism
in the Sec. 1 of this paper. If we appear to be too fussy
e.g. in the Definition 1.3, this is because we found out
that uncritical reliance on illustrative pictures
can be really misleading.

\smallskip

One can and must approach operadic constructions
from various directions and with various stocks of analogies.
In this paper, we look at operads, especially those with
values in abelian categories,
as analogs of associative rings; collections are analogs
of their generating spaces. We imagine various noncommutative
geometries based upon operads, and are interested in
naturally emerging symmetry and moduli objects in these
noncommutative geometries.

\smallskip

But of course there are many more different intuitive ideas related
to operads.

\smallskip

a) Operads provide  tools for studying general
algebraic structures determined by a basic set, a family
of composition laws, and a family of constraints
imposed upon these laws.

\smallskip

b) Operads embody a categorification of graph theory
which can be used to study knot invariants, Feynman
perturbation series etc.

\smallskip

c) Operads and their algebras are a formalization
of computational processes and devices, in particular,
tensor networks and quantum curcuits, cf. [MarkSh],
[Zo] and references therein. With this in mind, we describe
general endomorphism operads in 2.5 below.

\smallskip

It is interesting to notice
that the classical theory of recursive functions
must refer to a very special and in a sense
universal algebra over a non--linear ``computational
operad'', but nobody so far was able to formalize the latter.
Main obstacle is this: a standard description of any
partially recursive function produces a circuit
that may contain cycles of an a priori unknown multiplicity and
eventually infinite subprocesses producing no output at all.

\medskip

{\bf 0.4. Plan of the paper.} In Sec. 1, we discuss the background
topics. The centerpieces of the first
part related to graphs are definitions 1.2.1 and 1.3,
and the rest is devoted to collections and operads.

\smallskip

In Sec. 2, we state and prove our main theorems in two contexts:
for operads in abelian categories and for algebras. However,
the latter requires serious additional restrictions.
We also discuss in 2.7, 2.8 some
explicit descriptions of cohomomorphism objects,
whenever they are known, in particular, for quadratic
and more general $N$--homogeneous algebras and operads.

\smallskip

Notice that Theorem 0.1.1
was extended in [GrM] to include the case of {\it twisted}
tensor products of algebras: see also [Ma4] where the latter
appeared in the construction of the De Rham complex
of quantum groups and spaces. For further developments see
[GrM] and references therein. Similar generalizations
might exist for operads as well.

\smallskip

Sec. 3 is dedicated to the next layer of generalizations.
Namely, operads considered in Sec. 2 and Sec. 3
can be viewed as algebras over a triple whose
main component is the endofunctor $\Cal{F}$ on a category
of collections, described
in 1.5.5 and 1.5.6. Categories of collections
in this context are abelian and endowed with
a symmetric monoidal structure. In Sec. 3,
we consider more general triples and, in particular, do not
assume that the category on which the relevant endofuctor acts,
is abelian.

\smallskip

This line of thought is continued in Sec. 4, where in particular
operads with values in various categories of algebras are considered.
This is partly motivated by the quantum cohomology operad.

\smallskip

Finally, in Appendix we briefly treat Markl's list ([Mar], p. 45)
in terms of labeled graphs and functors on them.

\medskip

{\it Acknowledgements.} Yu.~M. gratefully acknowledges Bruno
Vallette's comments on the preliminary drafts of this paper and
discussions with him, in particular related to [Va2]. His numerous
suggestions are incorporated in the text. D.~B. would like to
thank Ezra Getzler for his remarks relating to model categories
and the theory of operads in $2$-categorical setting.

\bigskip

\centerline{\bf \S 1. Background}

\medskip

{\bf 1.1. Graphs.} We define objects of the category of (finite) graphs
as in [KoMa], [BeMa], [GeKa2]. Geometric realizations of our
graphs are not necessarily connected. This allows us
to introduce a
monoidal structure ``disjoint union''
on graphs (cf. 1.2.4), and to consider certain morphisms
such as graftings and mergers which were not needed in [GeKa2]
but arise naturally in more general types of operads.
Moreover, our notion of a graph morphism
is strictly finer than that considered in the literature:
as a part of a morphism, we consider the involution $j_h$
in the Definition 1.2.1 below. Our basic category of sets
is assumed to be small.

\medskip

{\bf 1.1.1. Definition.} {\it A graph $\tau$ is a family
of finite sets and maps
$(F_{\tau}, V_{\tau}, \partial_{\tau}, j_{\tau})$
Elements of $F_{\tau}$ are called flags of $\tau$,
elements of $V_{\tau}$ are called vertices of $\tau$.
The map $\partial_{\tau}:\, F_{\tau} \to V_{\tau}$
associates to each flag a vertex, its boundary.
The map $j_{\tau}:\, F_{\tau}\to F_{\tau}$
is an involution: $j_{\tau}^2= id.$}

\smallskip

{\it Marginal cases.} If $V_{\tau}$ is empty, $F_{\tau}$ must be
empty as well. This defines an {\it empty graph.}
To the contrary, $F_{\tau}$ might be empty whereas
$V_{\tau}$ is not. In [GeKa2] and other places,
in order to treat various units,
a ``non--graph'' with one flag and no vertices is considered.
Its role in our constructions sometimes
can be played by the empty graph.

\smallskip

{\it Edges, tails, corollas.} One vertex graphs with identical $j_{\tau}$
are called {\it corollas}.
Let $v$ be a vertex of $\tau$,
$F_{\tau}(v) :=\partial_{\tau}^{-1}(v)$. Then
$\tau_v:= (F_{\tau}(v), \{v\}, \roman{evident}\ \partial ,
\roman{identical}\ j)$ is a corolla, which
is called the corolla of $v$ in $\tau$.

\smallskip

Flags fixed by $j_{\tau}$ form the set of
{\it tails} of $\tau$ denoted $T_{\tau}.$

\smallskip

Two--element orbits of $j_{\tau}$ form the set
$E_{\tau}$ of {\it edges} of $\tau$.   Elements of
such an orbit are called {\it halves} of the respective
edge.

\medskip

{\bf 1.1.2. Geometric realization of a graph.}
First, let $\tau$ be a corolla. If its set of flags is empty,
its geometric realization $|\tau |$ is, by definition, a point.
Otherwise construct a disjoint union of segments $[0,\dfrac{1}{2}]$
and identify in it all points $0$. This is $|\tau |$. The image of
$0$ thus becomes the geometric realization of the unique vertex of
$\tau$.

\smallskip

Generally, to construct $|\tau |$ take a disjoint union
of geometric realizations of corollas of all vertices
and identify points $\dfrac{1}{2}$ of any two flags
forming an orbit of $j_{\tau}$.

\smallskip

A graph $\tau$ is called connected (resp. simply connected,
resp. tree etc) iff its geometric realization is such.
In the same vein, we can speak about connected components of a graph
etc. Vertices $v$ with empty $F_{\tau}(v)$ are
considered as connected components.

\medskip

{\bf 1.2. Morphisms of graphs and monoidal category $Gr$.}
Let $\tau$, $\sigma$ be two graphs.

\medskip

{\bf 1.2.1. Definition.} {\it A morphism $h:\,\tau\to \sigma$
is a triple $(h^F,h_V, j_h)$, where
$h^F:\,F_{\sigma}\to F_{\tau}$ is a contravariant map,
$h_V:\,V_{\tau}\to V_{\sigma}$ is a covariant map,
and $j_h$ is an involution
on the set of tails of $\tau$ contained
in  $F_{\tau}\setminus h^F(F_{\sigma})$. This data
must satisfy the following conditions.

\smallskip

(i) $h^F$ is injective, $h_V$ is surjective.

\smallskip

(ii) The image $h^F(F_{\sigma})$ and its complement $F_{\tau}\setminus h^F(F_{\sigma})$
are $j_{\tau}$--invariant subsets of flags.
The involution $j_h$ fixes no tail in
$F_{\tau}\setminus h^F(F_{\sigma})$.

\smallskip

It will be convenient to extend $j_h$ to other flags in
$F_{\tau}$ by identity.

\smallskip

We will say that $h$ contracts all flags in
$F_{\tau}\setminus h^F(F_{\sigma})$.
If two flags in $F_{\tau}\setminus h^F(F_{\sigma})$
form an edge, we say that this edge is contracted by $h$.
If two tails in $F_{\tau}\setminus h^F(F_{\sigma})$
form an orbit of $j_h$, we say that
it is a virtual edge contracted by $h$.

\smallskip

(iii) If a flag $f_{\tau}$ is not contracted by $h$,
that is, has the form $h^F(f_{\sigma})$, then
$h_V$ sends $\partial_{\tau}f_{\tau}$ to
$\partial_{\sigma}f_{\sigma}$. Two vertices
of a contracted edge (actual or virtual) must have the same $h_V$--image.

\smallskip

(iv) The bijection $h_F^{-1}:\, h^F(F_{\sigma})\to F_{\sigma}$
maps edges of $\tau$ to edges of $\sigma$.

If it maps a pair
of tails of $\tau$ to an edge of $\sigma$,
we will say that $h$ grafts these tails.

\smallskip

The composition of two morphisms corresponds to the set--theoretic
composition of the respective maps $h^F$ and $h_V$, and taking
the union of two sets of virtual edges.

\smallskip

The resulting category is denoted $Gr$.}

\medskip

{\bf 1.2.2. Geometric realization of a morphism.}
On geometric realizations, the action of $h$
can be visualized as follows: we construct
a subgraph of $\tau$ consisting only of flags in
$h^F(F_{\sigma})$, then produce its quotient, and then
identify this subquotient with $\sigma$
using $(h^F)^{-1}$.

\smallskip

The shortest description of this subquotient is this:
{\it merge in $|\tau |$ all vertices belonging to each one
fiber of $h_V$, then delete all flags which are contracted by $h$.}

\smallskip

It is easy to see that $(h^F)^{-1}$ identifies
the geometric graph thus obtained with $|\sigma |$.

\smallskip

This short and intuitive description may be misleading for
important concrete categories $\Gamma$
of labeled graphs (see Definition 1.3). It might happen that
such a category does not allow morphisms which simply delete flags,
and/or does not allow morphisms that merge two vertices
without contracting a path of edges
(actual or virtual) connecting such a pair.

\smallskip

The following sequence of steps has more chances
to represent a sequence of morphisms in $\Gamma$.

\smallskip

a) In each $j_h$--orbit, graft tails of $|\tau |$ belonging to
$F_{\tau}\setminus h^F(F_{\sigma})$ making an actual edge
from the virtual one.

\smallskip

Contract to a point each connected component
of the union of all actual or virtual edges whose halves belong to
$F_{\tau}\setminus h^F(F_{\sigma})$. This point becomes a new vertex.

\smallskip

Besides
contracted halves of edges, all other flags
adjacent to various vertices of the contracted
component are retained and become adjacent to
the vertex that is the image of this component.
Thus the set of remaining flags consists exactly of the
(geometric realizations of) $h^F(F_{\sigma})$.

\smallskip

From Def. 1.2.1 (iii) it follows that all flags adjacent to the new
vertex are sent by $(h^F)^{-1}$ to a subset of flags
adjacent to one and the same vertex of $\sigma$.

\smallskip

b) Graft loose ends of each pair of remaining tails
that are grafted by $h$.

\smallskip

c) Merge together those vertices of the obtained (geometric) graph
whose preimages are sent by $h_V$ to one and the same vertex
of $\tau$.

\smallskip

Finally, identify the resulting subquotient of $|\tau |$
with $|\sigma |$ in such a way that on flags
this map becomes $(h^F)^{-1}:\, h^F(F_{\sigma})\to F_{\tau}.$

\smallskip

Notice that all steps a) -- c) could have been done in arbitrary order,
for example c), b), a), with the same result.
Notice also that in the geometric realization,
no trace of $j_h$ is remained: the virtual edges
vanished, but we do not know how their
halves were paired. This information is encoded
only in the combinatorial description involving $j_h$.

\smallskip

Each step above is in fact a geometric realization of
a morphism in $Gr$. We will now describe formally the respective
classes of morphisms.

\medskip

{\bf 1.2.3. Contractions, graftings, mergers.}
a) {\it Virtual contractions.} A morphism $h:\,\tau \to \sigma$
is called {\it a
virtual contraction}, if
$F_{\tau}\setminus  h^F(F_{\sigma})$ consists only of tails,
restriction of $j_{\tau}$ on $ h^F(F_{\sigma})$ coincides
with the image of $j_{\sigma}$, and
$h_V$ is a bijection.

\smallskip

b) {\it Contractions and full contractions.}
A morphism is called {\it a contraction}, if $h^F$ is bijective on tails,
and if for any $v\in V_{\sigma}$, any two different
vertices in $h_V^{-1}(v)$ are connected by a path consisting
of edges contracted by $h$.
\smallskip

Let $\sigma$ be a graph. Define $\tau$ as follows:
$F_{\tau}:=T_{\sigma}$, $V_{\tau}$:=\{connected components of $\sigma$\}.
Let $h^F$ be the identical injection $T_{\sigma} \to F_{\sigma}$.
Let $h_V$ send any vertex to the connected component in which it is
contained. The resulting morphism is called {\it the full contraction.}
Its image is a union of corollas, tails of $\sigma$ are
distributed among them as they are among connected components of $\sigma$.
Morphisms isomorphic to such ones are also called full contractions
(of their source).

\smallskip

c) {\it Grafting and total grafting.} A morphism $h$ is called
{\it a grafting}, if $h^F$
and $h_V$ are bijections.

\smallskip

Let $\tau$ be a graph.
Denote by $\sigma := \coprod_{v\in V_{\tau}} \tau_v$
the disjoint union of corollas of all its vertices.
Formally, $F_{\sigma}=F_{\tau}$, $V_{\sigma}=V_{\tau}$,
$\partial_{\sigma}=\partial_{\tau}$, and $j_{\sigma} =id.$
Let $h:\,\sigma \to \tau$ consist of identical maps.
Such a morphism is called {\it total grafting}, and we will reserve for it
a special notation:
$$
\circ_{\tau}:\,\coprod_{v\in V_{\tau}} \tau_v \to \tau
\eqno(1.1)
$$
It is defined uniquely by its target $\tau$,
up to unique isomorphism identical on $\tau$.
Any isomorphism of targets induces an isomorphism of
such morphisms.
Its formal inversion cuts all edges of $\tau$ in half.
Morphisms isomorphic to such one are also called total grafting
(of their target).

\smallskip

d) {\it Mergers and full mergers.}  A morphism $h$ is called
{\it a merger}, if $h^F$ is bijective and identifies
$j_{\sigma}$ with $j_{\tau}.$

\smallskip

{\it A full merger} projects all vertices into one.
All edges become loops; tails remain tails.

\smallskip

Mergers play no role in the theory of ordinary operads, but are essential
for treating PROPs, cf. Appendix.

\smallskip

Each graph $\tau$ admits at least one morphism to a corolla 
$$
con:\,\tau \to \sigma
\eqno(1.2)
$$
It can be described as a grafting, followed by the full contraction,
followed by the full merger. Notice that in general $con$ is not unique, moreover, there is no meaningful canonical choice of $con$. This is due to the fact that in our approach, it is not the objects of $Gr$, that represent concatenations of operations, but morphisms into corollas.

\smallskip

Isomorphisms constitute the intersection of all four classes.

\smallskip

A combinatorial argument
imitating 1.2.2 shows that any morphism can be decomposed into a product
of a virtual contraction, a contraction, a grafting and a merger.
As soon as the order of the types of morphisms
is chosen, one can define such a decomposition in a
canonical way, up to a unique isomorphism.

\medskip

{\bf 1.2.4. Disjoint union as a monoidal structure on $Gr$.}
Disjoint union of two abstract sets having no common elements 
is an obvious notion.
If we want to extend it to ``all'' sets, an appropriate formalization
is that of a symmetric monoidal structure ``direct sum'' with empty set as
the unit object. It exists, but is neither unique, nor completely obvious:
what is the ``disjoint union of a set with itself''?
One way to introduce such a structure is described in [Bo2], Example 6.1.9.

\smallskip

We will focus on a small category of finite sets of all cardinalities 
and sketch the following method which neatly accounts
for proliferation of combinatorics of symmetric groups
in the standard treatments of operads.

\smallskip

The small category of finite sets of all cardinalities consisting of
$$
\emptyset , \{1\}, \{1,2\}, \dots ,\{1,2, \dots ,n\}, \dots
$$
admits a monoidal structure ``disjoint union'' $\coprod$ given by
$$
\{1,2, \dots ,m\}\coprod
\{1,2, \dots ,n\} :=
 \{1,2, \dots ,m+n\}
$$
and evident commutativity and associativity constraints.
For example, identification of $X\coprod Y$ with
$Y\coprod X$ proceeds by putting all elements of $Y$
before those of $X$ and retaining the order inside groups.
Empty set is the unit of this monoidal structure.

\smallskip

In these terms, it is clear how to extend
this construction to the category of ``all''
totally ordered finite sets, and then to drop
the orderings by passing to appropriate colimits.
Thus, we can endow the category of finite sets by a monoidal
structure which we keep denoting  $\coprod$
and calling ``disjoint union''.

\smallskip

This monoidal structure can then be extended to $Gr$:
$\sigma\coprod\tau$ is determined by disjoint unions of
their respective flag and vertex sets, and $\partial$,
$j$ act on both parts as they used to.

\smallskip

Finally, for any finite family $\{\tau_s\,|\,s\in S\}$,
we can define $\coprod_s \tau_s$ functorially in $\{\tau_s\}$
and $S$ as is spelled out in [DeMi].

\medskip

{\bf 1.2.5. Atomization of a morphism.} Let $h:\,\tau \to \sigma$
be a morphism of graphs. We define
its {\it atomization} as a commutative diagram of the following
form:
$$
\xymatrix{\coprod_{v\in V_{\sigma}}\tau_v \ar[d]_{k}
\ar[r]^{\coprod h_v}
& \coprod_{v\in V_{\sigma}} \sigma_v
\ar[d]^{\circ_{\sigma}} \\
\tau \ar[r]^h & \sigma}
\eqno(1.3)
$$
Here $\sigma_v$ is the corolla of a vertex $v\in V_{\sigma}$,
$\circ_{\sigma}$ is the  total grafting morphism,
and the remaining data are constructed as follows.

\smallskip

{\it Graph $\tau_v$.} We put for $v\in V_{\sigma}$:
$$
F_{\tau_v}:=\{f\in F_{\tau}\,|\,h_V(\partial_{\tau}f)=v\},\
V_{\tau_v}:=\{w\in V_{\tau}\,|\,h_V(w)=v\},
$$
$$
\partial_{\tau_v} =\partial_{\tau}\,|_{\tau_v},\
j_{\tau_v} =j_{\tau}\,|_{\tau_v}.
$$

\smallskip

{\it Morphism $h_v:\,\tau_v\to\sigma_v$.} We put
$$
h_v^F:=h^F\,|_{F_{\sigma_v}}:\, F_{\sigma_v}\to F_{\tau_v},\
h_{v,V}:=h_V\,|_{V_{\tau_v}}:\, V_{\tau_v}\to V_{\sigma_v},\
j_{h_v}:=j_h\,|_{F_{\tau_v}}.
$$
\smallskip

{\it Morphism $k$.} By definition, $k^F$ and $k_V$ are identical
maps, hence $k$ is a grafting.

\medskip

{\bf 1.2.6. Heredity.} Let now $\circ_{\sigma}:\,\coprod_{v\in V_{\sigma}} \sigma_v \to \sigma$ be a total grafting morphism. Assume that we
are given  a family of morphisms
$h_v:\,\tau_v\to\sigma_v$, $v\in V_{\sigma}$. Then this data
can be uniquely extended to the atomization diagram (1.3)
of a morphism $h:\,\sigma \to \tau$.

\medskip

{\bf 1.3. Definition.} {\it An abstract category of labelled
graphs is a category $\Gamma$ endowed with a functor $\psi
:\,\Gamma \to Gr$ satisfying the following conditions:

\smallskip

\smallskip

(i) $\Gamma$ is endowed with a monoidal structure
which $\psi$ maps to the disjoint union in $Gr$. It will be denoted
by the same sign $\coprod$.

\smallskip

(ii) $\psi$ is faithful: if two morphisms with common source and target
become equal after applying $\psi$, they are equal.

\smallskip

(iii) Call  a $\Gamma$--corolla any object $\tau$ of $\Gamma$
such that $\psi (\tau )$ is a corolla.
Any object  $\tau \in \Gamma$
admits at least one morphism to
a $\Gamma$--corolla
$$
con:\,\tau \to \sigma
$$
which is a lift to $\Gamma$ of the diagram of the form (1.2)
with the source $\psi (\tau )$.

\smallskip

(iv) Any object of $\Gamma$ is the target of a morphism from
a disjoint union of $\Gamma$--corollas
$$
\circ_{\tau}:\,\coprod_{v\in V_{\psi (\tau )}} \tau_v \to \tau
$$
which is a lift to $\Gamma$ of the diagram of the form (1.1)
with the target $\psi (\tau )$. It is defined uniquely
up to unique isomorphism identical on $\tau$.

\smallskip

(v) Any morphism $h:\,\tau\to\sigma$ can be embedded
into a commutative diagram of the form (1.3)
which is a lift of the atomization diagram
for the morphism $\psi (h):\,\psi(\tau )\to\psi(\sigma )$.
Such a diagram is defined uniquely up to unique isomorphism.

\smallskip

(vi) Moreover, assume that we are given $\sigma$ in $\Gamma$,
and
for each $\Gamma$--corolla $\sigma_v$, $v\in V_{\psi (\sigma )}$,
a morphism $h_v:\,\tau_v\to\sigma_v$ where $\{\tau_v\}$
is a family of objects of $\Gamma$. Then there exists a morphism
$\tau\to\sigma$ in $\Gamma$ such that all this data
fit into (a lift of) the atomization diagram (1.3).
Moreover, $\tau$ and $\sigma$ are defined uniquely
up to unique isomorphism. }

\medskip

The last requirement formalizes what Markl calls
``hereditary'' property in [Mar], p. 45.

\medskip

{\bf 1.3.1. Comments.} a) It is helpful (and usually
realistic) to imagine
any object $\sigma \in \Gamma$ as a pair consisting of
the ``underlying graph'' $\psi (\sigma )$ and an additional
structure on the components of $\psi (\sigma )$ such as decorating
vertices by integers, a cyclic
order on flags adjoining to a vertex, etc (see examples below).
We will generally refer to such a structure as ``labeling''.

\smallskip

As a rule, existence of labeling of a given type
and/or additional algebraic properties required
for a treatment a certain type of operadic structure
put some restrictions upon underlying graphs
so that $\psi$ need not be surjective on objects.
On the other hand, if these restrictions are satisfied,
there might be many different compatible labeling
on the same underlying graph so that $\psi$ need not be
injective on objects either.

\smallskip

The functor $\psi$ on objects simply forgets labeling.

\smallskip

b) In the same vein, any morphism $\sigma \to \tau$ in $\Gamma$
should be imagined as a morphism $\psi (\sigma )\to \psi (\tau )$
of underlying graphs satisfying some constraints of two
types: purely geometric ones which can be
stated in $Gr$, and certain compatibility conditions
with labelings (see below). This is the content
of condition (ii) above.

\smallskip

c) Let $\sigma \in \Gamma$. Slightly abusing the language,
we will call flags, vertices, edges etc. of $\psi (\sigma )$
the respective components of $\sigma$, and write, say, $V_{\sigma}$
in place of $V_{\psi (\sigma )}$. Similar conventions will apply
to morphisms in $\Gamma$. By extension,
$\sigma$ is called a $\Gamma$--corolla (resp. tree etc), if $\psi (\sigma )$
is a corolla (resp. tree etc).

\medskip

{\bf 1.3.2. Examples of labeling.} a) {\it Oriented graphs.}
Any map $F_{\sigma} \to \{in, out\}$ such that
halves of any edge are oriented by different labels, is called {\it an
orientation of $\sigma$.} On the geometric realization,
a flag marked by $in$ (resp. $out$) is oriented towards
(resp. outwards) its vertex.

\smallskip

Tails of $\sigma$ oriented $in$ (resp. $out$) are called
{\it inputs} (resp. {\it outputs}) of $\sigma$.
Similarly, $F_{\sigma}(v)$ is partitioned into inputs and outputs
of the vertex $v$.

\smallskip

Consider an orientation of $\sigma$. Its edge is
called {\it an oriented loop}, if both its halves belong
to the same vertex. Otherwise an oriented edge
starts at a source vertex and ends at a different
target vertex.

\smallskip

More generally, a sequence of pairwise distinct
edges $e_1, \dots ,e_n$, is called {\it a simple path}
of length $n$, if $e_i$ and $e_{i+1}$ have a common vertex, and
the $n-1$ vertices obtained in this way are
pairwise distinct. If moreover $e_1$ and $e_n$
also have a common vertex distinct from the mentioned ones,
this path is {\it a wheel} of length $n$.
A loop is a wheel of length one.
Edges in a wheel are endowed only with a cyclic order up
to inversion.

\smallskip

Clearly, all edges in a path (resp. a wheel)
can be oriented so that the source of $e_{i+1}$
is the target of $e_i$.

\smallskip

If the graph is already oriented,
the induced orientation on any path (resp. wheel) either
has this property or not. Respectively, the wheel is called oriented
or not.

\smallskip

{\it A morphism of oriented graphs} $h$
is a morphism of graphs such that $h^F$
is compatible with orientations.

\smallskip

b) {\it Directed graphs.} An oriented graph
$\sigma$ is called {\it directed} if it
satisfies the following condition:

{\it On each connected component of the geometric realization,
one can define a continuous real valued function (``height'')
in such a way that
moving in the direction of orientation along each flag
decreases the value of this function.}

\smallskip

In particular, a directed graph has no oriented wheels.

\smallskip

Notice that, somewhat counterintuitively,
a directed graph is not necessarily
oriented ``from its inputs to its outputs'' as is usually
shown on illustrating pictures. In effect,
take a corolla with only {\it in} flags and
another corolla with only {\it out}
flags, and graft one input to one output.
The resulting graph is directed (check this) although
its only edge is oriented from global
outputs to global inputs.

\smallskip

This is one reason why it is sometimes sensible
to include in a category $\Gamma$ of directed
graphs only those, which have at least
one input and least one output at each
vertex (cf. the definition of reduced bimodules in
Sec. 1.1 ov [Va1]).

\smallskip

Another reason for excluding certain marginal
{\it (``unstable'')} types of labeled corollas
might be our desire to ensure essential finiteness of the categories
denoted $\Rightarrow \sigma$ in Sec. 1.5.5
(cf. a description of unstable modular corollas below).
In a category of (disjoint unions of) directed trees, for example,
this leads to the additional requirement: corolla of any
vertex has at least three flags.

\smallskip

This requirement might lead to some technical
problems if we want to consider unital versions of our operads.

\smallskip

c) {\it Genus labeling.} A genus labeling of $\sigma$ is a
map $g: V_{\sigma}\to \bold{Z}_{\ge 0}, v\mapsto g_v.$
The genus of a {\it connected} labeled graph $\sigma$
is defined as
$$
g(\sigma ):= \sum_{v\in V_{\sigma}} g_v +
\roman{rk}\,H_1(|\sigma|) =
\sum_{v\in V_{\sigma}} (g_v -1) + \roman{card}\,E_{\sigma} +1.
$$
Genus labeled graphs (or only connected
ones) are called {\it modular graphs}.
Corolla of any vertex of a modular graph is a modular graph.

\smallskip

A morphism of modular graphs is a morphism of graphs
compatible with labeling in the following sense.
Contraction of a looping edge raises the genus of its vertex by one.
Contraction of a non--looping edge prescribes
to the emerging vertex the sum of genera
of two ends of the edge. Finally, grafting
flags does not change genera of vertices.

\smallskip

Thus, a morphism between two connected modular
graphs can exist only if their genera coincide.

\smallskip

A modular corolla with vertex of genus $g$
and $n$ flags is called {\it stable} iff
$2g-2+n>0.$ A modular graph is called {\it stable},
iff corollas of all its vertices are stable.

\smallskip

d) {\it Colored graphs.} Let $I$ be an abstract set
(elements of which are called colors).
{\it An $I$--colored graph} is a graph $\tau$
together with a map $F_{\tau}\to I$ such that
two halves of each edge get the same color.
Morphisms are restricted by the condition
that $h^F$ preserves color.

\smallskip

In [LoMa], a topological operad was studied
governed by a category of colored graphs
with two--element $I=\{black,\,white\}$.
Halves of an edge in this category are always white.

\smallskip

e) {\it Cyclic labeling.} A cyclic labeling of
a graph $\tau$ is a choice of cyclic order upon each
set $F_{\tau}(v)$. Alternatively, it is a family of
bijections $F_{\tau}(v)\to \mu_{|v|}$ where $\mu_{|v|}$
is the group of roots of unity (in $\bold{C}$)
of degree $|v|:=F_{\tau}(v)$; two maps which
differ by a multiplication by a root of unity
define the same labeling.
 Yet another description
identifies a cyclic labeling of $\tau$ with
a choice of planar structure for each corolla $\tau_v$,
that is, an isotopy class of embeddings of $|\tau_v|$
into an {\it oriented} plane.

\smallskip

In a category of cyclic labeled graphs,
mergers are not allowed, whereas contractions,
say, of one edge, should
be compatible with cyclic labeling in an evident way:
say $(0, 1,\dots  ,m)$ and $(0,\bar{1},\dots ,\bar{n})$
turn into $(1,\dots  ,m,\bar{1},\dots ,\bar{n})$, where by
$0$ we denoted the contracted halves of the same edge in two
corollas.

\smallskip

An interesting variant of cyclic labeling is
{\it unoriented cyclic labeling} (it has nothing
to do with orientation of the graph itself):
the cyclic orders $(0, 1,\dots  ,m)$ and
$(m, m-1,\dots  ,0)$ are considered equivalent.
In this version, contraction of an edge leads generally to
two different morphisms in $\Gamma$, with two
different targets.

\smallskip

Combinatorics of unoriented cyclic labeled trees is very
essential in the description of the topological operad
of real points $\overline{M}_{0,*}(\bold{R})$,
cf. [GoMa].

\medskip

{\bf 1.4. Ground categories $\Cal{G}$.} Operads
of various types in this paper will be defined
as certain functors from a category of labeled graphs $\Gamma$ to
a symmetric monoidal category $(\Cal{G}, \otimes )$ which will be called
{\it ground category.} The simplest example is that
of finite--dimensional vector spaces over a field,
or that of finite complexes of such spaces.

\smallskip

In order to ensure validity of various constructions
we will have to postulate (locally) some additional properties of
$(\Cal{G}, \otimes )$, the most important of which are
contained in the following list. At this stage, we do not assume
that all of them, or some subset of them, hold
simultaneously.

\smallskip

a) Existence of a unit object.

\smallskip

b) Existence of internal cohom objects $\underline{cohom}\,(X,Z)$
for any objects $X,Z$ in $\Cal{G}$. By definition,
they fit into functorial isomorphisms
$$
\roman{Hom}_{\Cal{G}}(X, Y\otimes Z)=
\roman{Hom}_{\Cal{G}}(\underline{cohom}\,(X,Z),Y).
$$
These isomorphisms are established by composition with coevaluation
morphisms
$$
c=c_{X,Z}:\, X\to \underline{cohom}\,(X,Z)\otimes Z
$$
(cf. the diagram (2.3) in Section 2 below).

\smallskip

c) Existence of countable coproducts such that $\otimes$ is distributive
with respect to these coproducts.

\smallskip

d) Existence of finite (and sometimes infinite) colimits.

\smallskip

e) $\Cal{G}$ is an abelian category, $\otimes$ is an
additive bifunctor exact in each argument.

\smallskip

f) $\Cal{G}$ is a closed model category.

\medskip

{\bf 1.5. $\Gamma\Cal{G}$--collections.}
We will denote by $\Gamma COR$ the subcategory
(groupoid) of $\Gamma$ consisting of $\Gamma$--corollas
and {\it isomorphisms} between them.

\medskip

{\bf 1.5.1. Definition.} {\it
A $\Gamma\Cal{G}$--collection $A_1$ is a functor
$A_1:\, \Gamma COR \to \Cal{G}$.  A morphism of
$\Gamma\Cal{G}$--collections $A_1\to B_1$ is
a functor morphism (natural transformation.)}

\smallskip

The category of $\Gamma\Cal{G}$--collections
will be denoted $\Gamma\Cal{G}COLL$.

\smallskip

{\bf 1.5.2. Examples.} a) If $\Gamma$ is the category of
stable modular graphs, the category
$\Gamma\Cal{G}COLL$ is equivalent to the category
of double sequences $A_1((g,n))$ of
objects of $\Cal{G}$ endowed with an action of $\bold{S}_n$
upon $A_1((g,n))$. A morphism $A_1\to B_1$ is a sequence
of morphisms in $\Cal{G}$, $A_1((g,n))\to B_1((g,n))$,
compatible with $\bold{S}_n$--actions.

\smallskip

In effect, any $\Gamma$--collection $A_1$ is determined
up to an isomorphism
by its restriction to the category of modular corollas
with flags $\{1,\dots ,n\}$ and a vertex labeled by $g$.
Their isomorphisms correspond to permutations of
$\{1,\dots ,n\}$.

\smallskip

Such collections are called stable $\bold{S}$--modules
in [GeKa2], (2.1).

\smallskip

b) Let $\Gamma$ be a category of oriented graphs,
containing all oriented corollas. Then the category
$\Gamma\Cal{G}$--collections is equivalent to the category
of double sequences  $A_1(m,n)$ of objects in $\Cal{G}$
endowed with actions of $\bold{S}_m\times \bold{S}_n$,
and equivariant componentwise isomorphisms.

\smallskip

This is clear: look at oriented corollas with
inputs $\{1,\dots ,n\}$ and outputs $\{1,\dots ,m\}$:
see e.g. [Va1], 1.1.

\medskip

{\bf 1.5.3. White product of collections.} For two collections
$A_1,B_1$, define their white product by
$$
(A_1\bigcirc B_1) (\sigma ):= A_1(\sigma )\otimes B_1(\sigma)
$$
for any $\Gamma$--corolla $\sigma$,
with obvious extension to morphisms.

\smallskip

This determines a symmetric monoidal structure $\bigcirc$ on
$\Gamma\Cal{G}COLL$.

\smallskip

If $(\Cal{G},\otimes )$ is endowed with a unit object $u$,
then the collection $U$, $U(\sigma )=u$, sending each
morphism of corollas to $id_u$ is a unit object of
$(\Gamma\Cal{G}COLL, \bigcirc )$.

\medskip

{\bf 1.5.4.  Inner cohomomorphisms for collections.}
If $(\Cal{G}, \otimes )$ admits internal cohom
objects, the same holds for $(\Gamma\Cal{G}COLL, \bigcirc )$:
just work componentwise.

\medskip

{\bf 1.5.5. Endofunctor $\Cal{F}$ on $\Gamma\Cal{G}COLL$.}
Let $A_1$ be a $\Gamma\Cal{G}$--collection. Consider
a $\Gamma$--corolla $\sigma$ and denote
by $\Rightarrow\sigma$  the category
whose objects are $\Gamma$--morphisms
$\Gamma$--graphs $\tau\to\sigma$, and whose
morphisms are $\Gamma$--isomorphisms of morphisms identical
on $\sigma$. Put
$$
\Cal{F}(A_1)(\sigma ) := \roman{colim}_{\Rightarrow\sigma}
\otimes_{v\in V_{\tau}}A_1 (\tau_v)
\eqno(1.4)
$$
The existence of appropriate colimits in $\Cal{G}$
such that $\otimes$ is distributive with respect to them
should be postulated at this stage.

\smallskip

In some important cases (e.g. stable modular graphs)
any category $\Rightarrow\sigma$ is essentially finite
(equivalent to a category with finitely many objects
and morphisms). Therefore existence of finite colimits
in $\Cal{G}$ suffices.

\smallskip

Clearly, this construction is functorial with respect to isomorphisms
of $\Gamma$--corollas so that we actually get
a new collection $\Cal{F}(A_1)$.
Moreover, the map $A_1\mapsto \Cal{F}(A_1)$
extends to an endofunctor of $\Gamma\Cal{G}COLL$.

\smallskip

As is well known, functor composition
endows the category of endofunctors by the structure
of strict monoidal category with identity.

\medskip

{\bf 1.5.6. Proposition.} {\it The endofunctor $\Cal{F}$
 has a natural structure
of a triple, that is, a monoid with identity
in the category of endofunctors.}

\smallskip

{\bf Proof {\it (sketch).}} The argument is essentially the same
as in (2.17) of [GeKa2].
We have to construct a multiplication morphism
$\mu :\, \Cal{F}\cdot\Cal{F}\to\Cal{F}$,
an identity morphism $\eta :\,\roman{Id} \to \Cal{F}$,
and to check the commutativity of several diagrams.

\smallskip

In other words, for a variable collection $A_1$,
we need functorial morphisms of collections
$\mu_{A_1}:\,\Cal{F}^2(A_1)\to \Cal{F}(A_1)$,
$\eta_{A_1}:\,A_1 \to \Cal{F}(A_1)$ fitting the relevant commutative
diagrams. In turn, to define them, we have
to give their values (in $\Cal{G}$) on any $\Gamma$--corolla
$\sigma$, functorially with respect to $\sigma$.

\smallskip

The construction of $\mu$ essentially uses
the hereditary property 1.3(vi) of the category $\Gamma$.
In fact, we have
$$
\Cal{F}^2(A_1)(\sigma )= \roman{colim}_{\tau \to\sigma}
\otimes_{v\in V_{\tau}} \Cal{F}(\tau_v) =
$$
$$
\roman{colim}_{\tau \to\sigma}
[\roman{colim}_{\rho_v\to\tau_v} \otimes_{w\in V_{\rho_v}}
A_1(\rho_{v,w})]
\eqno(1.5)
$$
where $\rho_v$ are objects of $\Gamma$, and $\rho_{v,w}$
is the $\Gamma$--corolla of a vertex $w$ of $\rho_v$.
Using heredity, we can produce from each family $\rho_v\to \tau_v$
a morphism $\rho \to \tau$, and replace the right hand side
of (1.5) by
$$
\roman{colim}_{\rho\to \tau \to\sigma}
\otimes_{w\in V_{\rho}} A_1(\rho_w)
$$
The latter colimit maps to $\Cal{F}(A_1)(\sigma )$
via composition of two arrows in $\rho\to\sigma$.

\smallskip

As Getzler and Kapranov suggest, this construction
and various similar ones needed to produce
$\eta_{A_1}$ and to check axioms, become
more transparent if one uses the simplicial formalism.

\smallskip

Given $\sigma$ and $k\ge 0$, define
the category $\Rightarrow_k\sigma$:
its objects are sequences of  morphisms
$(f_1, \dots ,f_k)$ in $\Gamma$,
$f_1:\,\tau_0\to \tau_1, \dots ,f_k:\,\tau_{k-1}\to\tau_{k}$
together with an augmentation morphism $\tau_k\to \sigma$.
Morphisms in
$\Rightarrow_k\sigma$ are isomorphisms of such sequences
compatible with augmentation.

\smallskip

Categories $\Rightarrow_k\sigma$ are interconnected
by the standard face and degeneracy functors
turning them into components of a simplicial
category.

\smallskip

Namely, $d_i :\, \Rightarrow_k\sigma
\to\, \Rightarrow_{k-1}\sigma$  skips $\tau_0,f_1$ (resp.
$f_k,\tau_k$)
for $i=0$ (resp. $i=k$); skips $\tau_i$ and composes
$f_i, f_{i+1}$ for $1\le i\le k-1$.

\smallskip

Similarly,  $s_i
:\, \Rightarrow_k\sigma
\to\, \Rightarrow_{k+1}\sigma$ inserts $id:\,\tau_i \to \tau_i$.

\smallskip

An argument similar to one which we sketched above
for $k=1$ will convince the reader that
one can identify the value of the
$(k+1)$--th iteration  $\Cal{F}^{k+1}(A_1)$ at $\sigma$
with the functor sending $\sigma$ to
$$
\roman{colim}_{\Rightarrow_k\sigma}
\otimes_{v\in V_{\tau_0}}A_1 (\tau_{0,v}) .
$$

\smallskip

Thus the multiplication
$\mu$ is the functor induced on colimits by $d_1 :\,
\Rightarrow_1\sigma \to\, \Rightarrow_0\sigma$.
Monoidal identity maps $A_1$ to $\Cal{F}(A_1)$
by sending $A_1(\sigma )$ to the diagram
$\roman{id}:\,A_1(\sigma )\to A_1(\sigma )$.

\smallskip

Two morphisms $\Cal{F}^3\to\Cal{F}^2$
corresponding to two configurations of brackets
are induced by $d_1$ and $d_2$ respectively, and
the associativity is expressed by the simplicial
identity $d_1d_1=d_1d_2.$ Identity $\eta$ is treated similarly.

\medskip

{\bf 1.6.  $\Gamma\Cal{G}$--operads as functors.} Traditional approaches
to operads in our context lead to three
different but equivalent definitions of it.
Very briefly, they can be summarized as follows.

\smallskip

(I) An operad is a tensor functor $(\Gamma ,\coprod )
\to (\Cal{G},\otimes )$ satisfying certain additional constraints.

\smallskip

(II) An operad is a collection together with a structure
of an algebra over the triple $(\Cal{F}, \mu ,\eta )$.

\smallskip

(III) An operad is a collection together
with a composition law which makes it a monoid with respect to
an appropriate symmetric monoidal structure
upon $\Gamma\Cal{G}COLL$ (to be described).

\smallskip

\smallskip

We will start with the first description,
and proceed to the second one in Sec. 1.7.
As for the third description which requires first
construction of a special
monoidal structure on collections,
it seems to be less universal. In the Appendix, we will sketch
it for orientation labeling and directed graphs, following
[Va1].

\smallskip

 Consider the category whose objects are
functors $A:\, \Gamma \to \Cal{G}$
{\it compatible
with the monoidal structures $\coprod$ and $\otimes$}
in the following sense: we are given functorial isomorphisms
$$
a_{\sigma ,\tau}:\, A(\sigma
\coprod \tau )\to A(\sigma )\otimes A(\tau )
\eqno(1.6)
$$
for all $\sigma , \tau \in \Gamma$ such that
inverse isomorphisms $a^{-1}_{\sigma ,\tau}$
satisfy
conditions spelled out in Def. 1.8 of [DeMi]. Such functors form
a category, morphisms in which are 
functor morphisms compatible with $a_{\sigma ,\tau}$.

\smallskip

In [DeMi] such functors are called {\it tensor functors,}
we will also use this terminology.

\medskip

{\bf 1.6.1. Definition.} {\it The category
$\Gamma\Cal{G}OPER$ of $\Gamma\Cal{G}$--operads is the category of
those tensor functors $(A,a):\,\Gamma \to \Cal{G}$ that
send any grafting morphism, in particular $\circ_{\tau}$,
in $\Gamma$ to an isomorphism.

\smallskip

Morphisms are functor morphisms.}

\smallskip

Informally, making grafting morphisms invertible
means that $A(\tau )$  for any
$\tau \in \Gamma$ can be canonically identified with
the tensor product $\otimes_v A(\tau_v)$ where $\tau_v$
runs over $\Gamma$--corollas of all vertices of
$\tau$ (or rather, of $\psi (\tau )$). To see this, one should
apply (1.6) (functorially extended to disjoint unions
of arbitrary finite families) to the l.h.s. of (1.1).

\smallskip

Moreover, a morphism into $\Gamma$-corolla $\tau\to\sigma$, postulated in Definition 1.3 (iii),
combined
with just described tensor decomposition
produces a morphism in $\Cal{G}$
$$
\otimes_{v\in V_{\tau}} A(\tau_v) \to A(\sigma).
\eqno(1.7)
$$
This is  our embodiment of operadic
compositions.

\smallskip

We will omit $a$ in the notation $(A,a)$ for brevity,
and simply treat (1.6) as identical map, as well as its
extensions to arbitrary families and their disjoint
unions.

\medskip

{\bf 1.7. From functors to algebras over the triple
$(\Cal{F},\mu ,\eta )$.} Let now $A:\,\Gamma \to \Cal{G}$
be an operad. Denote by $A_1$ its restriction
to the subcategory of $\Gamma$--corollas and their
isomorphisms. Let $\sigma$ be a $\Gamma$--corolla.
We can treat any object $\tau \to \sigma$ in $\Rightarrow \sigma$
(cf. 1.5.5) as a morphism in $\Gamma$ and apply
to it the functor $A$. We will get a morphism
$\otimes_{v\in V_{\tau}}A_1 (\tau_v) \to A_1(\sigma )$
functorial with respect to isomorphisms of $\tau$
identical on $\sigma$. Due to the universal property
of colimits, these morphisms induce a morphism
in $\Cal{G}$
$$
\Cal{F}(A_1)(\sigma )\to A_1(\sigma ).
$$
The system of these morphisms is functorial in $\sigma$,
so we get finally a morphism of collections
$$
\alpha_{A_1}:\, \Cal{F}(A_1)\to A_1 .
\eqno(1.8)
$$

Similarly, among the objects of  $\Rightarrow \sigma$
there is the identical morphism $id:\,\sigma \to \sigma$.
As above, it produces a morphism of collections
$$
\eta_{A_1}:\, A_1\to \Cal{F}(A_1).
\eqno(1.9)
$$

\medskip

{\bf 1.7.1. Proposition.} {\it (i) The data
$(A_1,\alpha_{A_1} ,\eta_{A_1})$ constitute an algebra over
the triple $(\Cal{F},\mu,\eta )$
(see e.g. [MarShSt], pp. 88--89).

\smallskip

(ii) The map $A\mapsto (A_1,\alpha_{A_1} ,\eta_{A_1})$
extends to a functor $coll$ establishing equivalence of
the category of $\Gamma\Cal{G}$--operads and the category
of algebras over the triple $(\Cal{F},\mu,\eta )$.}

\medskip

{\bf Proof {\it (sketch)}}. The proof is essentially the same
as that of Proposition (2.23) in [GeKa2].
The first statement reduces to the check of commutativity
of several diagrams.

\smallskip

The main problem in the second
statement consists in extending each morphism
of operads as algebras to a morphism of operads as functors.
In other words, knowing the operadic
compositions induced by full contractions and mergers
(if the latter occur in $\Gamma$),
we want to reconstruct operadic
compositions induced by partial contractions
such as $\tau\to\sigma$ in (a $\Gamma$--version of)
a diagram (1.3). To this end,
complete (1.3) by a morphism
$$
f:\,\coprod_{w\in V_{\tau}}\tau_w\to \coprod_{v\in V_{\sigma}}\tau_v
$$
which is the disjoint union of morphisms of total
graftings with targets $\tau_v$ (so that
$\tau_w$ are corollas whereas $\tau_v$ generally
are not).  Since $k\cdot f$ is
a total grafting with target $\tau$, in order to
calculate the value of our functor on $\tau\to\sigma$
it suffices to know its value on the composition
$\circ_{\sigma}(\coprod h_v) f$ which involves
only graftings, full contractions
and mergers of corollas.

\medskip

{\bf 1.7.2. Algebra $\Cal{F} (A_1).$} As a general formalism shows, for any $A_1$, the collection
$\Cal{F}(A_1)$ has the canonical structure
of an operad, with $\alpha_{\Cal{F}(A_1)}:=\mu_{A_1}:\, \Cal{F}^2(A_1)\to \Cal{F}(A_1)$.

\medskip

{\bf 1.7.3. Proposition.} {\it For any collection $A_1$ and an operad $B$,
there exists a canonical identification
$$
\roman{Hom}_{{}_{\Gamma\Cal{G}COLL}}(A_1, coll (B)) =
\roman{Hom}_{{}_{\Gamma\Cal{G}OPER}}(\Cal{F}(A_1), B)
$$
functorial in both arguments.}

\medskip

This means that the functor $\Cal{F}:\, \Gamma\Cal{G}COLL \to
\Gamma\Cal{G}OPER$ is a construction of the free operad freely generated
by a collection, and is thus an analog of tensor algebra of
a linear space.

\medskip

{\bf 1.8. White product of operads.} For two $\Gamma\Cal{G}$--operads
$A,B$, define their white product by
$$
(A\bigcirc B) (\sigma ):= A(\sigma )\otimes B(\sigma)
$$
for any object  $\sigma\in \Gamma$.
The extension to morphisms is evident.

\smallskip

This determines a symmetric monoidal structure $\bigcirc$ on
$\Gamma\Cal{G}OPER$.

\smallskip

As before, if $(\Cal{G},\otimes )$ is endowed with a unit object $u$,
then the functor $U:\, \sigma\mapsto u$ sending each
morphism of labeled graphs to $id_u$ is a unit object of
$(\Gamma\Cal{G}OPER, \bigcirc )$.

\medskip

{\bf 1.9.  Morphism $j$.}
Morphism of operads
$j:\,\Cal{F}(E_1\bigcirc B_1)\to \Cal{F}(E_1)\bigcirc \Cal{F}(B_1)$
that is, a family of morphisms
$$
\roman{colim}_{\Rightarrow \sigma} \otimes_{v\in V_{\tau}} E_1(\tau_v)\otimes
B_1(\tau_v) \to (\roman{colim}_{\Rightarrow \sigma}
\otimes_{v\in V_{\tau{'}}}E_1(\tau_v{'}))\otimes (\roman{colim}_{\Rightarrow \sigma} \otimes_{v\in V_{\tau{''}}} B_1(\tau_v{''}))
$$
comes from the ``diagonal'' part of the right hand side:
$\tau = \tau{'}=\tau{''}.$

\bigskip

\centerline{\bf \S 2.  Inner cohomomorphism operads}

\medskip

{\bf 2.1. Preparation.} Fix a graph category $\Gamma$ and
a ground category $\Cal{G}$ as above.

\smallskip

In this Section, we will assume that $(\Cal{G},\otimes )$
has  inner cohomomorphism objects. Moreover, for validity
of the main Theorem 2.2, $\Cal{G}$ must be
abelian, with tensor product exact in each argument.

\smallskip

Let $A$ be a  $\Gamma\Cal{G}$--operad, $A_1$ a
$\Gamma\Cal{G}$--collection, and $i_A:\,A_1\to A$
a morphism of collections such that the respective
morphism of operads $f_A:\,\Cal{F}(A_1)\to A$
is surjective. Denote by $\Cal{A}$ the diagram $i_A:\,
A_1\to A$. Such a diagram can be thought of as
a presentation of $A$.

\smallskip

Let $\Cal{B}$ be a similar data $i_B:\,B_1\to B.$

\smallskip

As in 0.1, denote by $\Cal{A}\Rightarrow \Cal{B}$
the category whose objects are commutative diagrams
$$
\xymatrix{A_1 \ar[d]_{i_A} \ar[r]^{u_1} & F\bigcirc B_1
\ar[d]^{id_F\bigcirc i_B} \\
A \ar[r]^u & F\bigcirc B}
\eqno(2.1)
$$
where $F$ is an operad, $u$ is a morphism of operads,
and $u_1$ is a morphism of collections. In the upper row,
and below in similar situations, we write $F$ in place
$coll (F)$  for brevity, identifying an operad with its
underlying collection.

\smallskip

Notice that, unlike in 0.1, we do not assume
that $i_A,i_B$ are injective. Hence the upper row
of (2.1) has to be given explicitly.

\smallskip

Such a diagram can be denoted $(F,u,u_1)$ since the remaining
data are determined by $\Cal{A}$, $\Cal{B}$.

\smallskip

A morphism $(F,u,u_1)\to (F^{\prime},u^{\prime},u_1^{\prime})$
in $\Cal{A}\Rightarrow \Cal{B}$ is a morphism of
operads $F\to F^{\prime}$ inducing a morphism of
commutative diagrams (2.1) constructed for $F$
and $F^{\prime}$ respectively.

\medskip

{\bf 2.2. Theorem.} {\it The category $\Cal{A}\Rightarrow \Cal{B}$
has an initial object
$$
\xymatrix{A_1 \ar[d]_{i_A} \ar[r]^{\delta_1} & E\bigcirc B_1
\ar[d]^{id_E\bigcirc i_B} \\
A \ar[r]^{\delta} & E\bigcirc B}
\eqno(2.2)
$$
defined uniquely up to unique isomorphism.

\smallskip

Moreover, $E$ comes together with a presentation $\Cal{E}$,\
$i_E:\,E_1 \to E$ in which $E_1 = \underline{cohom} (A_1,B_1)$,
 inner cohomomorphism being taken in the category of $\Gamma\Cal{G}$--collections.

\smallskip

If $F$ is given together with its presentation $\Cal{F}$, that is
$i_F:\,F_1\to F$,
and $u$ is induced by $u_1:\,A_1\to F_1\bigcirc B_1$,
then the canonical
homomorphism $E\to F$ is induced by a unique
morphism in the category $\Cal{E}\Rightarrow \Cal{F}$.}

\medskip

{\bf Proof.} {\it (i) Preparation.} The morphism $u_1:\,A_1\to F\bigcirc B_1$
corresponds to a morphism $\tilde{u}_1:\,E_1\to F$
where   $E_1 = \underline{cohom} (A_1,B_1)$ as above.
Recall that  inner cohomomorphism collections
here can be constructed componentwise.
Let $c:\,A_1\to E_1\bigcirc B_1$ be the coevaluation
morphism. Then the diagram
$$
\xymatrix{A_1 \ar[d]_{id_{A_1}} \ar[r]^{c} & E_1\bigcirc B_1
\ar[d]^{\tilde{u}_1\bigcirc id_{B_1}} \\
A_1 \ar[r]^{u_1} & F\bigcirc B_1}
\eqno(2.3)
$$
is commutative. Composing (2.3) with (2.1), we get
a commutative square of collections
$$
\xymatrix{A_1 \ar[d]_{i_A} \ar[r]^{c} & E_1\bigcirc B_1
\ar[d] \\
A \ar[r]^u & F\bigcirc B}
\eqno(2.4)
$$
which produces a morphism of operads with presentations
$$
\xymatrix{\Cal{F}(A_1) \ar[d]_{f_A} \ar[r]^{\Cal{F}(c)} &
\Cal{F}(E_1\bigcirc B_1)
\ar[d]^g \\
A \ar[r]^u & F\bigcirc B}
\eqno(2.5)
$$
This diagram can be completed by the commutative triangle
$$
\xymatrix{\Cal{F}(A_1) \ar[d]_{f_A} \ar[r]^{\Cal{F}(c)} &
\Cal{F}(E_1\bigcirc B_1) \ar[d]^g
\ar[r]^{h}&
\Cal{F}(E_1)\bigcirc B \ar[ld]^{\tilde{v}\bigcirc id_B}\\
A \ar[r]^u & F\bigcirc B}
\eqno(2.6)
$$
where $h$ is the composition
$$
\xymatrix{\Cal{F}(E_1\bigcirc B_1) \ar[r]^{j} & \Cal{F}(E_1)\bigcirc \Cal{F}(B_1)
\ar[r]^{id\bigcirc f_B} & \Cal{F}(E_1)\bigcirc B}
$$
and $j$ is described in Section 1.9.

\smallskip

{\it (ii) Construction of main objects.} Now we will construct in turn an ideal $(\tilde{R})\subset \Cal{F}(E_1)$,
the operad $E:=\Cal{F}(E_1)/(\tilde{R})$ together with a morphism of operads
$\delta:\, A\to E\bigcirc B$ and a morphism $v:\,E\to F.$

\smallskip

Starting with this point, we will have to use our assumption
that assume that $\Cal{G}$ is abelian,
and $\otimes$ is exact.

\smallskip

Choose a subcollection $R\subset \roman{Ker}\,f_A$ generating
$\roman{Ker}\,f_A$ as an ideal in the free operad
$\Cal{F}(A_1)$. Replace
the morphism $h\circ \Cal{F}(c):\,\Cal{F}(A_1)\to \Cal{F}(E_1)\bigcirc B$
by the morphism canonically corresponding to it
$$
\underline{cohom}(\Cal{F}(A_1),B)\to \Cal{F}(E_1)
$$
where the  inner cohomomorphisms here and  below are taken
in the category of collections.

\smallskip

Since  inner cohomomorphisms are covariant functorial with
respect to the first argument, we have a commutative diagram
$$
\xymatrix{\underline{cohom}\,(R,B) \ar[d] \ar[rd] \\
\underline{cohom}\,(\Cal{F}(A_1),B) \ar[d] \ar[r] & \Cal{F}(E_1)
\ar[d]^{\tilde{v}}\\
\underline{cohom}\,(A,B) \ar[r] & F}
\eqno(2.7)
$$
where $\tilde{v}$ is taken from (2.6).

\smallskip

Denote by $\tilde{R}$ the image of the skew arrow in (2.7).
Since the composition $R\to \Cal{F}(A_1)\to A$
is zero, the same holds for the composition of the two
left vertical arrows in (2.7). This implies
that the composition $\tilde{R}\to \Cal{F}(E_1)\to F$
is zero. Since $\tilde{v}$ is a morphism of operads,
its kernel contains the ideal $(\tilde{R})$ generated
by $\tilde{R}.$

\smallskip

Rewrite the upper triangle in (2.7) as a commutative square
$$
\xymatrix{\underline{cohom}\,(R,B) \ar[d] \ar[r]&\tilde{R} \ar[d] \\
\underline{cohom}\,(\Cal{F}(A_1),B)  \ar[r] & \Cal{F}(E_1)}
$$

Replacing the horizontal arrows with the help of
coevaluation morphisms, we get the commutative square
$$
\xymatrix{R \ar[r] \ar[d] & \tilde{R}\bigcirc B \ar[d]\\
\Cal{F}(A_1) \ar[r] & \Cal{F}(E_1)\bigcirc B}
\eqno(2.8)
$$
which induces a morphism of operads
$$
\delta:\, A=\Cal{F}(A_1)/(R) \to \Cal{F}(E_1)\bigcirc B/(\tilde{R}\bigcirc
B)\cong \Cal{F}(E_1)/(\tilde{R})\bigcirc B = E\bigcirc B
\eqno (2.9)
$$
(use the exactness of $\otimes$ in $\Cal{G}$.)

\smallskip

{\it (iii) Completion of the proof.} It remains to show that
$\delta:\, A\to E\bigcirc B$ has all the properties stated in
Theorem 2.2.

\smallskip

From our construction, it is clear that it fits into
the diagram of the form (2.2), and that it comes
with the presentation $i_E:\,E_1\to E.$

\smallskip

The morphism $\tilde{v}:\, \Cal{F}(E_1)\to F$ descends to
a morphism of operads $v:\, E\to F$, and from the commutativity
of the diagram (2.6) one can infer that it induces
a morphism of objects $(E,\delta ,\delta_1)\to
(F,u,u_1)$.

\smallskip

We leave the remaining checks to the reader.

\medskip

{\bf 2.3. Remark.} Associative algebras (without unit) can be treated
as operads: functors on linear oriented trees
and contractions with values in linear spaces.
The Theorem 2.2 in this case reduces to the Theorem
0.1.1, or rather its extension where presentations
are not supposed to be injective on components of degree 1.

\medskip

{\bf 2.4.  Inner cohomomorphisms for operads with
presentation.} We can now reformulate
Theorem 2.2 in the same way as it was done
in the Introduction for associative algebras.

\smallskip

Consider the following category $OP$ of $\Gamma\Cal{G}$
operads given together with their presentation.

\smallskip

Objects of $OP$ are pairs $\Cal{A}$ as in 2.1.
A morphism $\Cal{A}\to \Cal{B}$ is a pair
consisting of a morphism of collections
$A_1\to B_1$ and a morphism of operads $A\to B$
compatible with $i_A,i_B$.

\smallskip

This category has a symmetric monoidal product $\bigodot$
induced by $\bigcirc$ in the following sense:
$\Cal{A}\bigodot\Cal{B}=\Cal{C}$ where $\Cal{C}$
is represented by $C_1 :=A_1\bigcirc B_1$,
$C:=$ the minimal suboperad containing the
image $(i_A\bigcirc i_B)(C_1)$ in $A\bigcirc B$,
and $i_C:=$ restriction of $i_A\bigcirc i_B$.

\smallskip

Theorem 2.2 can now be read as a statement that, functorially
in all arguments, we have
$$
\roman{Hom}_{OP} (\Cal{A}, \Cal{F}\bigodot\Cal{B})=
\roman{Hom}_{OP} (\Cal{E}, \Cal{F})
\eqno(2.10)
$$
that is, $\Cal{E}$ is an  inner cohomomorphism object in $OP$:
$$
\Cal{E} = \underline{cohom}_{OP} (\Cal{A},\Cal{B}).
\eqno(2.11)
$$

General categorical formalism produces canonical
comultiplication morphisms in $OP$
$$
\Delta_{\Cal{A},\Cal{B},\Cal{C}}:\,
\underline{cohom}_{OP} (\Cal{A},\Cal{C})\to
\underline{cohom}_{OP} (\Cal{A},\Cal{B})\bigodot
\underline{cohom}_{OP} (\Cal{B},\Cal{C})
\eqno(2.12)
$$
coassociative in an evident sense.

\smallskip

This is an operadic version of quantum matrices and
their comultiplication.

\smallskip

In particular, the operad
$$
\underline{coend}_{OP} \Cal{A}:=
\underline{cohom}_{OP} (\Cal{A},\Cal{A})
$$
is endowed with a canonical coassociative
comultiplication, morphism of operads
$$
\Delta_\Cal{A}:=\Delta_{\Cal{A},\Cal{A},\Cal{A}}:\,\underline{coend}_{OP} \Cal{A}\to
\underline{coend}_{OP} \Cal{A}\bigodot\underline{coend}_{OP} \Cal{A}.
\eqno(2.13)
$$
It is generally not cocommutative, as the case of associative rings
amply demonstrates.

\smallskip

We get thus a supply of ``quantum semigroups'',
or Hopf algebras in the category of operads
(ignoring antipode).
\smallskip

Notice finally that if $\Cal{G}$ has a unit object $u$, we can
sometimes define a unit object $U$ in $(OP,\bigodot )$. Taking
$\Cal{F}=U$ in the adjunction formula (2.10), we see that the
space of ``classical homomorphisms'' $\Cal{A}\to \Cal{B}$ in $OP$
coincides with space of points of $\Cal{E}$ with values in $U$
(here we implicitly imagine operadic affine quantum spaces as
objects of the dual category).

\medskip

{\bf 2.5. Algebras over an operad and their deformations.}
For certain labelings and respective graph categories
$\Gamma$ one can define a class of ``natural''
$\Gamma\Cal{G}$--collections.

\smallskip

The basic example is this.
Let $J$ be an abstract set of ``flavors''
such that a labeling of $\tau$ consists
of a map $F_{\tau}\to J$ (and possibly other data).
(Imagine orientations, colors, or pairs (orientation, color)).

\smallskip

Assume furthermore that automorphisms of any $\Gamma$--corolla
$\sigma$ form a subgroup of permutations of its flags preserving
flavors of flags. In this case any family of objects
$\bold{V}:=\{V_j\,|\,j\in J\}$ determines
the following $\Gamma$--collection:
$$
Coll\, (\bold{V})(\sigma ):=
\otimes_{j\in J}V_j^{\otimes F_{\tau}^{(j)}}
\eqno(2.14)
$$
where  $F_{\tau}^{(j)}$ is the subset of flags of flavour $j$.
Automorphisms of $\sigma$ act in an evident way.

\smallskip

Such a collection naturally extends to a functor
on the groupoid of $\Gamma$--graphs $\tau$
and their isomorphisms: simply replace $\sigma$ by $\tau$
in (2.14).

\smallskip

In order to extend it to arbitrary $\Gamma$--morphisms,
consider separately three classes of morphisms.

\smallskip

a) {\it Graftings}. Graftings correspond to bijections
of sets of flags preserving flavors. Hence they extend
to (2.14).

\smallskip

b) {\it Mergers}. They have the same property.

\smallskip

c) {\it Contractions.} In principle, in order to accommodate
contractions, we have to impose on $\bold{V}$ an additional
structure, namely, a set of polylinear forms on $(V_j)$ which we
denote $\bold{v}$ and axiomatize as follows:

\smallskip

{\it Let $S$ be a finite set and $\kappa :=\{(j_s,k_s)\,|\,s\in S\}$
be a family of pairs of flavors such that in some $\Gamma$--graph
there exist two vertices (perhaps coinciding) and
connecting them edges which are simultaneously
contracted by a $\Gamma$--morphism. The respective
component $v_{\kappa}$ of $\bold{v}$ is a morphism in $\Cal{G}$
$$
v_{\kappa}:\, \otimes_{s\in S} (V_{j_s}\otimes V_{k_s}) \to u
\eqno(2.15)
$$
where $u$ is the unit object of $\Cal{G}$.}

\smallskip

Given $\bold{v}$, the prescription for extending
our functor to contractions looks as follows.
To map a product (2.14) for the source of
a contraction to the similar product for its target,
we must map identically factors corresponding
to uncontracted flags, and to ``kill''
factors of the type $\otimes_{s\in S} V_{j_s}\otimes V_{k_s}$
with the help of (2.15).

\smallskip

This prescription will not necessarily describe
a functor $\Gamma \to \Cal{G}$: one should
to impose upon $\bold{v}$ coherence conditions
which we do not bother to spell out here.

\smallskip

In real life, this problem is avoided by
specifying only bilinear forms corresponding
to one--edge contractions, and then tensoring them
to obtain full scale (2.15).

\smallskip

We can now give the main definition of this Section.

\medskip

{\bf 2.5.1. Definition.} {\it a) Any functor as above
$\Gamma \to \Cal{G}$ with underlying collection
$Coll\, (\bold{V})$ and structure forms $\bold{v}$
is called the endomorphism operad of $(\bold{V},\bold{v})$
and denoted $OpEnd\,(\bold{V},\bold{v})$.

\smallskip

b) Let $\Cal{P}$ be a $\Gamma\Cal{G}$--operad.
Any morphism $P\to OpEnd\,(\bold{V},\bold{v})$
is called a structure of $P$--algebra on $(\bold{V},\bold{v})$.}

\medskip

{\it Examples.} We will illustrate this on three types of
labelling discussed in 1.3.2 above.

\smallskip

a) Working with a subcategory
of $Gr$ itself (as in the case of cyclic operads)
we should choose one object $V$ of $\Cal{G}$
and a symmetric pairing $g:\,V\otimes V\to u$.

\smallskip

b) Let now $\Gamma$ be a subcategory of
oriented graphs. In that case one usually chooses
$V_{out}=V$, $V_{in}=V^t$ (the dual object),
and takes for $v$ the canonical pairing $V^t\otimes V\to u$.

\smallskip

c) Finally, let $\Gamma$ be a category of colored graphs, $I$ the
set of colors. In this case, one applies the full machinery of the
definition above, simplifying it by caring only about one--edge
contractions. So we need $V_i\otimes V_i\to u$ in the unoriented
case, or else choose $V_{out,i}=V_{in,i}^t$.

\smallskip

Let us now return to our definition 2.5.1.
The whole space of structures of  $P$--algebra on $(\bold{V},\bold{v})$
is thus
$$
\roman{Hom}_{\Gamma\Cal{G}OPER} (P, OpEnd\,(\bold{V},\bold{v})).
\eqno(2.16)
$$

\smallskip

In the standard approach to the
deformation theory of operadic algebras one chooses
a space of basic operations, that is, a presentation
of $P$ (and then replaces the respective structure
by a differential in an appropriate Hochschild--type
complex, the step that we will not discuss here).
Sometimes, one can choose a compatible presentation of
$OpEnd\,(\bold{V},\bold{v})$. For example, if connected graphs in
$\Gamma$ consist of oriented trees and  all basic operations are binary,
one can choose for generators the collection $P_1$ which coincides
with $P$ on corollas with three flags and is zero otherwise.
Respectively, $OpEnd\,V$ is (hopefully) generated by
$\underline{hom}(V^{\otimes 2},V)$ (at least, this is the case
for $\Cal{G}=Vec_k$.)

\smallskip

Accepting this, we suggest to replace $P$ by
$\Cal{P}$ which is $i:\,P_1\to P$, to
augment $ OpEnd\,V$ accordingly, and to replace (2.16)
by an appropriate set of morphisms in the category $OP$.
After that, it is only natural to consider the operad
$$
\underline{cohom}_{OP}( \Cal{P}, \Cal{O}p\Cal{E}nd\,V)
\eqno(2.17)
$$
as a non--commutative space parameterizing deformations of (a
chosen collection of generators of) $P$--algebra structures.

\medskip

{\bf 2.6.  Inner cohomomorphisms for algebras over an operad.}
 $P$--algebras form a category $PALG$.
One can try to play with $P$--algebras the same game as we did with
associative algebras and operads of various types,
and to study existence of  inner cohomomorphisms
in a category of $P$--algebras with a presentation.
However, even in order to state
the problem we need at least two preliminary
constructions:

\smallskip

a) Free $P$--algebra $F_P(\bold{V},\bold{v})$ generated by a family
$(\bold{V},\bold{v})$ as above. Generally, it will take values
in a monoidal category larger than $\Cal{G}$,
e.g. that of inductive systems.

\smallskip

This will allow us to define the notion of
an algebra with presentation.

\smallskip

b) Symmetric white product on the category of $P$--algebras
extending $\otimes$.

\smallskip

To this end, we need a symmetric comultiplication
$\Delta :\, P\to P\bigcirc P$ and an analog
of the morphism $j$ from Sec. 1.9 for free
$P$--algebras.

\smallskip

Unfortunately, as our description of $P$--algebras
above shows, we cannot do even this preliminary work
in the same generality as we treated operads themselves.
Therefore we step back, and for the remainder of Section 2.6
work with a version of ordinary operads.

\medskip

{\bf 2.6.1. Ordinary operads and their algebras.}
Let $\Gamma$ be the category of graphs whose connected
components are directed trees with exactly one output
and at least one input at each vertex.
Morphisms are contractions and graftings; mergers are not allowed.
Let $(\Cal{G},\otimes ,u)$ be an abelian symmetric monoidal
category, such that $\otimes$ is exact in both arguments, and
endowed with finite colimits and cohomomorphism objects.

\smallskip

Consider a $\Gamma\Cal{G}$ operad $P$.
Denote by $P(n)$  the value of $P$ on the $\Gamma$--corolla
with inputs $\{1,\dots ,n\}$, and let $\bold{S}_n$
be the automorphism group of this directed corolla.
Assume that $P(1)=u$ and that contracting
an edge one end of which carries a corolla with one input
produces canonical identifications $Q\otimes u\to Q$.

\smallskip

It is well known that a free $P$--algebra freely
generated by $V_1$ exists in an appropriate category of inductive limits,
and its underlying object is
$$
F_P(V_1)=\oplus_{n=1}^{\infty} P(n)\otimes_{\bold{S}_n}V_1^{\otimes n}.
\eqno(2.18)
$$
This construction is functorial in $V_1$.
{\it A presentation} of an algebra $V$ is a surjective
morphism $F_P(V_1)\to V$; it can be reconstructed
from its restriction $i_V:\,V_1\to V$.

\smallskip

Assume that $P$ is endowed with a symmetric comultiplication
$\Delta :\, P\to P\bigcirc P$.

\smallskip

Then, given two objects $E_1, W_1$ of $\Cal{G}$, we can define
a map
$$
j:\, F_P(E_1\otimes W_1)\to F_P(E_1)\bigcirc F_P(W_1).
\eqno(2.19)
$$
To construct it, first produce for each $n$ a map
$$
P(n)\otimes (E_1\otimes W_1)^{\otimes n} \to
P(n)\otimes E_1^{\otimes n}\otimes P(n)\otimes W_1^{\otimes n}
$$
combining $\Delta$ with regrouping, and then the induced
map of colimits
$$
j(n):\, P(n)\otimes_{\bold{S}_n} (E_1\otimes W_1)^{\otimes n} \to
(P(n)\otimes_{\bold{S}_n} E_1^{\otimes n})\otimes
(P(n)\otimes_{\bold{S}_n} W_1^{\otimes n}).
$$
We put $j=\oplus_n j(n)$. One can check that this is
a morphism of $P$--algebras.

\smallskip

If $V,W$ are two $P$--algebras,
presented by their structure morphisms $\alpha :\,F_P(V)\to V$,
$\beta :\,F_P(W)\to W$, we can compose $j(n)$ with $\alpha \otimes \beta$
to define a structure of $P$--algebra on $V\otimes W$.
This gives a symmetric monoidal structure on $PALG$ still denoted
$\otimes$. We can now state

\medskip

{\bf 2.6.2. Proposition.} {\it Analog of Theorem 2.2 holds in the
category of P--algebras with presentation.}

\smallskip

We skip a proof which follows the same plan as that of
Theorem 2.2.

\medskip

{\bf 2.7. Explicit constructions of cohomomorphism objects in the
category of associative algebras.} Existence proof of
cohomomorphism objects generally is not very illuminating. In this
and the following subsections, we cite several explicit
constructions, valid under additional assumptions.

\smallskip

a) {\it Quadratic algebras.} This was the case first treated in
[Ma1] and [Ma2]. Briefly, a quadratic algebra $A$ (over a field $k$)
is defined by its presentation $\alpha :\,F_{ASS}(A_1) \to A$
where $F_{ASS}(A_1)=T(A_1)$ is the free (tensor) algebra freely generated
by a finite--dimensional vector space $A_1$, such that
$\roman{Ker}\,\alpha$ is the ideal generated by the space
of quadratic relations $R_A\subset A_1^{\otimes 2}$.
In particular, $A_1$ is embedded in $A$, each $A$ is
naturally graded ($A_1$ in degree 1), and a morphism
of quadratic presentations is the same as morphism
of algebras which preserves grading. Moreover,
each morphism is uniquely defined by its restriction
to the space of generators.

\smallskip

This category $Qalg$ has a contravariant duality involution:
$A\mapsto A^!$ where $A^!= T(A_1^*)/(R_A^{\perp})$,
$A_1^*$ denoting the dual space to $A_1$.

\smallskip

It has also two different symmetric monoidal structures,
which are interchanged by $!$:
``white product'' $\circ$ and ``black product'' $\bullet$:
$$
T(A_1)/(R_A)\circ T(B_1)/(R_B) =
T(A_1\otimes B_1)/(S_{(23)}(R_A\otimes B_1^{\otimes 2}+
A_1^{\otimes 2}\otimes R_B),
\eqno(2.20)
$$
$$
T(A_1)/(R_A)\bullet T(B_1)/(R_B) =
T(A_1\otimes B_1)/(S_{(23)}(R_A\otimes R_B),
\eqno(2.21)
$$
Here $S_{(23)}:\, A_1^{\otimes 2}\otimes
B_1^{\otimes 2}\to (A_1\otimes B_1)^{\otimes 2}$ interchanges two middle factors.

\smallskip

Both categories have unit objects: polynomials of one variable for $\circ$,
and dual numbers $k[\varepsilon ]/(\varepsilon^2)$
for $\bullet$ respectively. They are $!$--dual to each
other. The generator $\varepsilon$ combined with general
categorical constructions produces differential
in various versions of Koszul complex which is a base
of Koszul duality.

\smallskip

The monoidal category $(Qalg,\bullet )$ has   internal homomorphism objects.
Explicitly,
$$
\underline{hom}_{\bullet} (A,B)=A^!\circ B .
\eqno(2.22)
$$

The monoidal category $(Qalg,\circ )$ has  inner cohomomorphism objects.
Explicitly,
$$
\underline{cohom}_{\circ} (A,B)=A\bullet B^! .
\eqno(2.23)
$$

{\bf Claim.} {\it (2.23) is a valid description of the internal
cohomomorphism object of two quadratic algebras
in the total category of (algebras with) presentations.}

\smallskip

These formulas follow from a general functorial
isomorphism (adjunction formula)
$$
\roman{Hom}_{Qalg} (A\bullet B^{!},C)=
\roman{Hom}_{Qalg} (A, B\circ C)
\eqno(2.24)
$$
\smallskip

Abstract properties of $Qalg$ expressed by (2.22)--(2.24)
can be axiomatized to produce an interesting
version of the notion of rigid tensor category
of [DeMi]. They justify the use of the name
``quantum linear spaces'' for objects of $Qalg$
although the category itself is not linear or even additive.

\smallskip

It might be worthwhile to consider  $(Qalg,\circ )$
and $(Qalg,\bullet )$ as ground categories
for $\Gamma$--operads (and cooperads).

\smallskip

Notice that the cohomology spaces of $\overline{M}_{0,n+1}$,
components of the Quantum Cohomology cooperad,
are quadratic algebras (Keel's theorem).

\smallskip

b) {\it $N$--homogeneous algebras.} It was shown in
[BerDW] that similar results hold
for the category $H_Nalg$ of homogeneous algebras generated
in degree 1 with relations generated in degree $N$, for any
fixed $N\ge 2$. If one continues to denote
by $R_A\subset A_1^{\otimes N}$
generating relations, $A^!$ is given by
relations $R_A^{\perp}\subset (A_1^*)^{\otimes N}$.
Thus $!=!_N$ explicitly depends on $N$ and gives, for example,
different dual objects of a free algebra, depending on where we put
its trivial relations.
In the definitions (2.20) one should make
straightforward modifications,
replacing $2$ by $N$. Formulas (2.22)--(2.24)
still hold in the new setup. Unit object for
$\bullet$ is now $k[\varepsilon ]/(\varepsilon^N).$

\smallskip

As a consequence, Koszul complexes become
$N$--complexes leading to an interesting
new homological effects: see [BerM] and references therein.

\smallskip

c) {\it Homogeneous algebras generated in degree one.}
This case is treated in Chapter 3 of [PP] and in Section 1.3
of [GrM]; the approaches in these two papers nicely
complement each other.

\smallskip

White product
(2.21) and its $H_Nalg$--version extend to this larger
category as $(A\circ B)_n:=A_n\otimes B_n$
(Segre product). The duality  morphism $!$ and the black
product do not survive in [PP],
and  internal homomorphism  objects (2.22) perish.
The right hand side
of the formula (2.23) is replaced by a rather
long combinatorial construction which we do not reproduce here,
and (2.24) becomes simply the characteristic
property of cohom's:
$$
\roman{Hom} (\underline{cohom}\,(A,B),C)=
\roman{Hom} (A, B\circ C)
$$
However, if $A$ is quadratic, then (2.23) can be resurrected in
a slightly modified form: $B$ is the quotient of  a quadratic
algebra $\roman{q}B$ (leave the same generators and only quadratic relations),
and we have
$$
\underline{cohom}_{\circ} (A,B)=A\bullet (\roman{q}B)^!
\eqno(2.25)
$$
(see [PP], Ch.3 , Proposition 4.3).

\smallskip

An important novelty of [PP] is a treatment of white products
and cohom objects in the  categories of graded modules
over graded rings.

\smallskip

In [GrM], a version of $!$ and several versions of
the black product are introduced. They are not as neatly packed together
however as in the cases $Qalg$ and $H_Nalg$. In particular,
$!$ is not an involution. As a compensation,
Theorem 1 of [GrM] establishes a nice extension
of (2.25) for general homogeneous $A$ and $B$
involving a certain ``triangle product''.

\smallskip

A logical next step would be the introduction of these
notions into ``non--commuta\-tive projective
algebraic geometry'' where  coherent sheaves
appear \`a la Serre as quotient categories of graded
modules, and cohom's of graded rings can serve as
interesting non--commutative correspondences.

\medskip

{\bf 2.8. Explicit constructions of cohomomorphism
objects in the categories of operads.} One objective of
the brief review above was to collect a list of patterns
that could be subsequently recognized in various categories of
operads. Most of the existing results which we are
aware of concern quadratic operads and Koszul duality patterns.
Cohomomorphism objects appear as a byproduct,
although, as we have seen above, their existence is
the most persistent phenomenon, even when the neat package
$(\circ ,\bullet ,!)$ cannot be preserved.

\smallskip

a) {\it Binary quadratic operads.} The pioneering paper [GiKa]
defined $(\circ ,\bullet ,!)$ for binary quadratic (ordinary)
operads. Their construction uses a description of
operads as monoids in $(COLL, \boxtimes )$,
a category of collections endowed with a monoidal
structure $\boxtimes$ (in 1.6 above, this is description (III)).
This description was rather neglected here (see Appendix), but it makes
clear the analogy between the tensor algebra of a linear space
and  the free operad generated by a collection $V$.

\smallskip

B.~Vallette in [Va3] describes a construction
of a free monoid which allows him to treat the cases when
the relevant monoidal structure is not biadditive,
which is the case of operads.
The resulting weight grading of the free monoid can be used to define
weight graded quotients, analogs of graded associative algebras: see
[Fr], [Va1], [Va2].

\smallskip

The subcategory of ordinary operads with presentation considered
in [GiKa] consists of weight graded operads generated
by their binary parts (values on corollas with two inputs),
with relations in weight 2.  After introducing $(\circ ,\bullet ,!)$,
Ginzburg and Kapranov prove the adjunction formula
(2.24) and thus the formula (2.23) as well.

\smallskip

b) Developing this technique, B.~Vallette in [Va1] and [Va2]
defines $\circ$ for properads (see Appendix) with presentation
and studies the case of (weight) quadratic relations
for generators of arbitrary arity. The full adjunction formula (2.24) is established for quadratic operads with generators
of one and the same arity $k\ge 2$ ([Va2], Sec. 4.6, Theorem 26),
thus generalizing the $k=2$ case of [GiKa].

\smallskip

Yet another version of this result is proved for
$k$--ary quadratic {\it regular} operads, with different
definitions of black and white products
([Va2], Sec. 5.1, Theorem 40). (A regular operad is an ordinary
symmetric operad
which is induced by some non--symmetric operad
in a sense that will not be made precise here).

\bigskip

\centerline{\bf \S 3. Non--abelian constructions}

\medskip

{\bf 3.0. Introduction.} In Sec.~2 it was proved that given an
abelian category $\Cal{G}$ with a symmetric monoidal structure
$\tensor$ and its left adjoint $\cohom$, the categories of operads
in $\Cal{G}$ and algebras over Hopf operads also possess $\cohom$ for
the natural extension of $\tensor$. In this Section we prove
a version of this
result in a non--abelian setting.

\smallskip

More precisely, we start with the description of operads
as algebras over the triple $(\Cal{F},\mu ,\eta )$,
cf.~Proposition 1.7.1. Here $\Cal{F}$ is an explicit endofunctor
on the category of collections $\Gamma\Cal{G}COLL$
which inherits from $\Cal{G}$ a monoidal structure, an
abelian structure, and  inner cohomomorphisms.

\smallskip

We replace $\Gamma\Cal{G}COLL$ with an abstract monoidal category
$(\C,\otimes )$ endowed with an endofunctor $T$ which has the
structure of a triple that commutes (up to a natural
transformation) with the monoidal structure. In this case, with
some additional assumptions, extension of $\cohom$ is
straightforward and is given by the adjoint lifting theorem. In
this formulation it is unnecessary to suppose that the monoidal
product is symmetric, however, when applied to operads, the triple
itself is produced using symmetric properties of the monoidal
structure.

\smallskip

Finally we formulate the natural notion of the derived $\cohom$
and consider some cases when such a functor exists.

\smallskip

{\bf 3.1. Tensor product and $\cohom$ for algebras over a triple.}
It is well known that if a category is equipped with a monoidal
structure distributive with respect to direct sums, then
categories of algebras over ordinary Hopf operads in this category
will possess an extension of the monoidal structure. We can treat
more general operads considering them in the context of triples.
In this Section we formulate conditions on the triples, needed to
extend monoidal structure to algebras over them.

\smallskip

Let $(\C,\tensor)$ be a category with a coherently associative
product (bifunctor). We do not assume $\tensor$ to be symmetric or
to possess a unit object.

\smallskip

{\bf 3.1.1. Definition.} {\it A Hopf--like triple on $\C$ is a triple
$\triple:\C\rightarrow\C$, $\mu :\,T\circ T\to T$, together with
a natural transformation between bifunctors on $\C$:
$$
\t:\triple\circ\tensor\twomor\tensor\circ(\triple\times\triple),
$$
satisfying the following conditions:

\smallskip

a) We have a commutative diagram of natural transformations:
$$
\xymatrix{T\circ T\circ \otimes \ar[d]_{T(\tau )}
\ar[r]^{\tau (\mu\circ Id)} &
\otimes\circ (T\times T)\\
T\circ \otimes \circ (T\times T) \ar[r]^{\tau} &
\otimes\circ (T\times T)\circ (T\times T)
\ar[u]_{\otimes (\mu\times \mu )}}
$$
where $\monpro:\triple^2\twomor\triple$ is the structure natural
transformation of $\triple$.

\medskip

b) Let
$\ass:\tensor\circ(\tensor\times\Id)\twomor\tensor\circ(\Id\times\tensor)$
be the associativity constraint for $\tensor$. Then
the following diagram of natural transformations is commutative:
$$
\xymatrix{T\circ\otimes\circ(\otimes\times Id) \ar[d]_{T(\alpha )}
 \ar[r]^{\tau} &
\otimes\circ (T\times T)\circ (\otimes \times Id)
\ar[r]^{\otimes (\tau \times Id)} &
\otimes \circ ((\otimes\circ(T\times T))\times T)
\ar[d]^{\alpha} \\
T\circ\otimes\circ(Id\times\otimes )
 \ar[r]^{\tau} &
\otimes\circ (T\times T)\circ (Id \times \otimes )
\ar[r]^{\otimes (Id \times \tau)} &
\otimes \circ (T\times (\otimes\circ(T\times T)))}
$$
}

\smallskip

The conditions in the definition above allow us to define an
associative product on the category of $\triple$-algebras, by
extending it from $\C$. The functor morphism $\t$ provides a definition,
condition  a) ensures that different ways of composing $\t$
produce the same result, and condition b) implies coherent
associativity of the resulting product on the category of
$\triple$-algebras by utilizing associativity isomorphisms of
$\tensor$. Here is the formal reformulation of all this.

\medskip

{\bf 3.1.2. Definition.} {\it Let $(A,\alpha:\triple(A)\rightarrow A)$,
$(B,\beta:\triple(B)\rightarrow B)$ be two $\triple$-algebras in
$\C$. We define $(A,\alpha)\white(B,\beta)$ to be $(A\tensor
B,\alpha\white\beta)$, where $\alpha\white\beta:\triple(A\tensor
B)\rightarrow A\tensor B$ is the composition}
$$
\triple(A\tensor B)\overset{\t}\to
{\rightarrow}\triple(A)\tensor\triple(B)
\overset{\alpha\tensor\beta}\to{\rightarrow}A\tensor B.
$$

\smallskip

{\bf 3.1.3. Lemma.} {\it Let $\A$ be the category of
$\triple$-algebras in $\C$. Defined as above, $\white$ is a
bifunctor on $\A$. It satisfies associativity conditions together
with coherence. The associativity isomorphisms are preserved by
the forgetful functor $\U:\A\rightarrow\C$.}

\smallskip

{\bf Proof.} To show that $\white$ is a bifunctor we have to show
first that its value on a pair of objects in $\A$ is again in
$\A$, i.e. that the following diagram is commutative
$$
\xymatrix{T^2(A\otimes B) \ar[d]_{T(\alpha\bigcirc\beta )}
\ar[r]^{\mu} &
T(A\otimes B) \ar[d]^{\alpha\bigcirc\beta }\\
T(A\otimes B) \ar[r]^{\alpha\bigcirc\beta} &
A\otimes B}
$$

This follows from the condition a) in the Definition 3.1.1 and the
commutativity of the following diagram:
$$
\xymatrix{T(T(A)\otimes T(B)) \ar[d]_{\tau}
\ar[r]^{T(\alpha\otimes\beta )} & T(A\otimes B) \ar[d]^{\tau}\\
T^2(A)\otimes T^2(B) \ar[d]_{\mu\otimes\mu}
\ar[r]^{T(\alpha )\otimes T(\beta )} &
T(A)\otimes T(B) \ar[d]^{\alpha\otimes\beta}\\
T(A)\otimes T(B) \ar[r]^{\alpha\otimes\beta} &A\otimes B}
$$

The upper square of this diagram is commutative because $\triple$
and $\tensor$ are functors and $\t$ is a natural transformation.
The lower square is commutative because it is the result of an application of
$\tensor$ to two commutative squares (representing the fact that
$\alpha$ and $\beta$ are structure morphisms for
$\triple$--algebras).

\smallskip

Thus on objects $\white$ behaves like a bifunctor on $\A$. We have
to show that $f\white g$ is a morphism in $\A$ for any two
morphisms $f:(A_1,\alpha_1)\rightarrow(A_2,\alpha_2)$ and
$g:(B_1,\beta_1)\rightarrow(B_2,\beta_2)$ in $\A$, i.e. the
following diagram is commutative
$$
\xymatrix{T(A_1\otimes B_1) \ar[d]_{\alpha_1\bigcirc\beta_1}
\ar[r]^{T(f\otimes g)} &
T(A_2\otimes B_2) \ar[d]^{\alpha_2\bigcirc\beta_2 }\\
A_1\otimes B_1 \ar[r]^{f\otimes g} &
A_2\otimes B_2}
$$
By the definition of $\white$, this diagram can be decomposed as
follows
$$
\xymatrix{T(A_1)\otimes T(B_1) \ar[d]_{\tau}
\ar[r]^{T(f\otimes g )} & T(A_2\otimes B_2) \ar[d]^{\tau}\\
T(A_1)\otimes T(B_1) \ar[d]_{\alpha_1\otimes\beta_1}
\ar[r]^{T(f)\otimes T(g)} &
T(A_2)\otimes T(B_2) \ar[d]^{\alpha_2\otimes\beta_2}\\
A_1\otimes B_1 \ar[r]^{f\otimes g} &A_2\otimes B_2}
$$

The upper square of this diagram is commutative because $\tensor$
and $\triple$ are functors and $\t$ is a natural transformation.
The lower square is commutative because it is the result of an application of
$\tensor$ to two commutative squares (expressing the fact that $f$
and $g$ are morphisms of $\triple$-algebras).

\smallskip

Clearly pairs of identities are mapped to identities by $\white$
and the composition is preserved, i.e. $\white$ is a bifunctor on
$\A$.

\smallskip

It remains to show that $\white$ is coherently associative and
that the forgetful functor to $\C$ preserves the associativity
isomorphisms. The former claim actually follows from the latter,
since $\tensor$ is coherently associative. So all we have to do is
to show that the $\tensor$-associativity isomorphisms between the
images of the forgetful functor belong to $\A$, i.e. the outer
rectangle of the following diagram is commutative for all
$A_1,A_2,A_3\in\A$
$$
\CD
T((A_1\otimes A_2)\otimes A_3)
@>{(\tau\otimes Id)\circ\tau}>>
(T(A_1)\otimes T(A_2))\otimes T(A_3)
@>{(\alpha_1\otimes\alpha_2)\otimes\alpha_3}>>
(A_1\otimes A_2)\otimes A_3 \\
@V{T(\alpha )}VV   @V{\alpha}VV  @V{\alpha}VV \\
T(A_1\otimes (A_2\otimes A_3))
 @>{(Id\otimes\tau)\circ\tau}>>
(T(A_1)\otimes (T(A_2)\otimes T(A_3))
@>{\alpha_1\otimes (\alpha_2\otimes\alpha_3)}>>
A_1\otimes (A_2\otimes A_3)
\endCD
$$
The right hand square of this picture is commutative because
$\tensor$ and $\triple$ are functors and $\ass$ is a natural
transformation. Commutativity of the left hand square is the
contents of condition b) of Definition 3.1.1.

\smallskip

Having extended the product $\tensor$ from $\C$ to $\A$, we would
like to know if this extension possesses a left adjoint, given
$\tensor$ does so on $\C$. We can infer it easily from the adjoint
lifting theorem, if we assume that $\A$ has all coequalizers. For
the question of when a category of algebras over a triple has all
coequalizers see e.g. [BarW], Section 9.3.

\medskip

{\bf 3.1.4. Proposition.} {\it Suppose that $\A$ has all
coequalizers and for every $A\in\A$ the functor
$\U(A)\tensor-:\C\rightarrow\C$ has a left adjoint. Then the
functor $A\white-:\A\rightarrow\A$ has a left adjoint as well.}

\smallskip

{\bf Proof.} It is clear that $\U$ is a monadic functor
([Bo1], Definition 4.4.1) and by our construction
$\U(A\white-)=\U(A)\tensor\U(-)$, therefore (e.g. [Bo2],
Theorem 4.5.6) since $\U(A)\tensor-$ has a left adjoint, so does
$A\white-$. This shows our assertion.

\smallskip

Here we define $\cohom(A,B)$ is an object that represents
$\hom(A,B\white-)$. From the last proposition we know that on $\A$
there is a $\cohom(-,-)$ that is functorial in the first argument.
Since Yoneda embedding is full and faithful, we conclude that our
construction is actually a bifunctor.

\medskip

{\bf 3.2. Derived $\cohom$.} By now we have constructed
$\cohom(-,A)$ for an algebra $A\in\A$, given existence of $\cohom$
on $\C$. Now assume that $\C$ carries in addition a closed model
structure, such that $C\tensor-$ is a right Quillen functor for any
$C\in\C$, i.e. it has a left adjoint and it maps fibrations and
trivial fibrations to the like. Then, by a general result,
$\cohom(-,C)$ is a left Quillen functor and we can define its left
derived version.

\smallskip

Suppose that we can transport closed model structure from $\C$ to
$\A$ through the adjunction of the free algebra and the forgetful
functor, i.e. we can introduce a closed model structure on $\A$ so
that a map in $\A$ is a weak equivalence or fibration if and only
its image under the forgetful functor is such. Then for any
$A\in\A$ we have that $A\white-$ is a right Quillen functor on
$\A$, and hence $\cohom(-,A)$ is a left Quillen functor and we can
define derived $\cohom$ as follows:
$$
\L\cohom(-,A):=\cohom(L(-),A),
$$
where $L$ is a cofibrant replacement functor on $\A$.

\smallskip

It remains only to analyze when such a transport of model
structure is possible. The general situation of a transport was
considered in several papers e.g. \cite{Bek} \cite{Bl},
\cite{CaGa}, \cite{Cra}, \cite{Q}, \cite{R}, \cite{S}. For our
purposes it is enough to consider locally presentable categories
with cofibrantly generated model structures. In this case the
conditions on the forgetful functor are quite mild.

Moreover, we are mostly interested in the particular case of
operads in the category of algebras over a Hopf operad in the
category of dg complexes of vector spaces over a field. In the
next Section we will show that categories of such operads can be
constructed as categories of algebras over certain triples on the
category of dg complexes of vector spaces. Then the machinery of
transport of model structure can be applied to these triples
directly.

\bigskip

\centerline{\bf \S 4. Iteration of algebraic constructions}

\medskip

{\bf 4.0. Introduction.} As it is described in Section 3, if a
triple is Hopf--like with respect to a monoidal structure, it is
easy to extend $\cohom$ from the ground category to the category
of algebras over this triple. Hopf operads provide examples of
such triples, given that monoidal structure on the ground category
is distributive with respect to the direct sum. However, this is
not always the case.

\smallskip

Consider the category $\As$ of associative algebras in the
category of vector spaces over a field. Since the operad of
associative algebras is Hopf, $\As$ has a symmetric monoidal
structure, given by the tensor product of algebras. On $\As$ this
monoidal structure is not distributive with respect to direct
sums, that is, free products of algebras. In fact, 
given $\as_1$, $\as_2$, $\as_2\in\As$, in general we
have no isomorphism between $\as_1\tensor(\as_2\coprod\as_3)$ and
$(\as_1\tensor\as_2)\coprod(\as_1\tensor\as_3)$. Thus we cannot
represent operads in $\As$ as algebras over a triple, and
constructions of Section 3 concerning existence of $\cohom$ do not
apply to the category of operads in $\As$.

\smallskip

However, we can overcome these difficulties by working with the
operads instead of their categories of algebras. In the example
above we can represent {\it operads} in the category of associative
algebras as {\it algebras over a colored operad} in the category of
vector spaces.

\smallskip

For example, a classical operad in $\As$ is a sequence
$\as:=\{\as_n\}_{n\in\Z_{>0}}$ of objects of $\As$ and morphisms
$$
\gamma_{m_1,\dots ,m_n}:\as_n\tensor \as_{m_1}\tensor \dots\tensor
    \as_{m_n}\rightarrow \as_{m_1+\dots +m_n},
$$
satisfying the usual axioms. The additional structures of
associative algebras on $\as_n$'s are described by a sequence of
morphisms of vector spaces:
$\alpha_n:\,\as_n\tensor\as_n\rightarrow\as_n$, that satisfy the
usual associativity axioms.

\smallskip

Compatibility of the operadic structure morphisms with the
structures of associative algebras on individual $\as_n$'s is
expressed by commutativity of the following diagrams:
$$
\CD (\as_n\otimes\as_n)\otimes (\otimes_{i=1}^n \as_{m_i}\otimes
\as_{m_i}) @>{\alpha_n\otimes (\otimes_{i=1}^n\alpha_{m_i})}>>
\as_n\otimes \as_{m_1}\otimes \dots \as_{m_n}\\
@V{(\gamma_{m_1,\dots ,m_n}\otimes \gamma_{m_1,\dots
,m_n})\circ\sigma}VV
@V{\gamma_{m_1,\dots ,m_n}}VV \\
\as_{m_1+\dots +m_n}\otimes \as_{m_1+\dots +m_n}
@>{\alpha_{m_1+\dots +m_n}}>>  \as_{m_1+\dots +m_n}
\endCD
$$
where $\sigma$ is the appropriate rearrangement of factors in the
tensor products of vector spaces.

\smallskip

It is easy to see that we can express all of these conditions on
the sequence $\{\as_n\}$ as an action of a colored symmetric
operad. Thus the category of classical operads in $\As$ is
equivalent to a category of algebras over a colored operad in
vector spaces. Therefore we have the free algebra construction,
results of Section 3 apply, and $\cohom$ extends to classical
operads in $\As$.

\smallskip

We would like to do the same in the case of
general operads, as described in Section 1. To achieve the
necessary degree of generality we will work with triples, which
will include all operadic cases, given that the ground category
has monoidal structure, which is distributive with respect to the
direct sums.

\smallskip

The usual definition of a triple on a category $\C$ consists of
three parts: a functor $\triple:\C\rightarrow\C$, and natural
transformations $\opmort:\triple\circ\triple\rightarrow\triple$,
$\opuntr:\Id\rightarrow\triple$, satisfying certain associativity
and unit axioms. The functor $\triple$ is supposed to represent
the ``free algebra'' construction, while $\opmort$ and $\opuntr$
represent composition and the identity operation respectively.

\smallskip

Often we have more information about the triple than contained in
its definition as above. We have a grading on the ``free algebra''
construction, given by the arity of operations involved, i.e. we
can decompose $\triple(C)$ ($C\in\C$) into a direct sum of
$\triple_n(C)$ ($n\in\Z_{>0}$), where $\triple_n(C)$ stands for
applying to ``generators'' $C$ all of the ``$n$-ary operations''
$\triple_n$.

\smallskip

Such decomposition of $\triple$ is very helpful, since usually
$\triple_n$'s behave better with respect to monoidal structure and
other triples than the whole $\triple$: see [Va3] for a closely
related discussion. Using this we can mimic construction of a
colored operad, as in the case of classical operads in associative
algebras above, in the more general situation of decomposed
triples.

\smallskip

In order to do so we have to formalize the notion of triples
admitting such a decomposition. We should axiomatize the relationship
between $\triple_n$'s for different $n$'s, so that the combined
object would satisfy the associativity and the unit axioms, stated
in the usual definition of a triple. The best way to do so is via
representations of operads in categories.

\smallskip

It is often the case (in particular it is so for operads as
described in Section 1) that the ``$n$-ary operations''
$\triple_n$ is not just a functor on $\C$, but is given as a
composition
$$\C\overset{\Delta}\to\rightarrow\C^{\times^n}\overset{\opelem_n}\to\rightarrow\C,$$
where $\Delta$ is the diagonal functor and $\opelem_n$ is some
functor $\C^{\times^n}\rightarrow\C$. In such cases it is
convenient to work with representations of operads on $\C$, i.e.
with morphisms of operads with codomain
$\{\Fun(\C^{\times^n},\C)\}_{n\in\Z_{n>0}}$.

\smallskip

So we will consider strict operads in the category of categories
and work with lax morphisms of such operads (i.e. with morphisms
that commute with the structure functors only up to a natural
transformation). A lax morphism into an endomorphism operad of a
category will give us a generalized triple that under some
conditions can be transformed into a usual triple. In such cases
we will say that the triple is operad--like. The technique of 
lax morphisms was  invented long ago (cf. [KS]) and applied recently
to the case of operads in [Bat].

\smallskip

We proceed as follows: we start with the well known notion of
strict pseudo--operads in categories, we organize them in a
category and then extend it to include lax morphism of operads,
which satisfy natural conditions of coherence (later these
conditions will be shown to correspond to the associativity axioms
of triples). Then we consider the notion of a strict operad in
categories (i.e. with a unit) and define lax morphisms of such
objects. Again we will need some coherence conditions, which later
will turn out to be the unit axioms of a triple. Finally we define
operad--like triples as lax morphisms into the endomorphism operad
of a category and we finish by proving that iteration of algebraic
constructions preserves existence of $\cohom$.

\medskip

{\bf 4.1. Notation.} We will denote functors usually by capital
letters (both Greek and Latin), whereas natural transformations
will be denoted by small letters.

\smallskip

Commutative diagrams of functors will be rarely commutative on the
nose, instead we will have to endow them with natural transformations,
making them
commutative. When we draw a diagram as follows
$$
\xymatrix{A \ar[d]_F \ar[rr]^G && B \ar[d]^{F'}\\
A' \ar[ru]^\alpha \ar[rr]_{G'} && B'}
$$
and say that $\alpha$ makes the diagram commutative, we mean that
$\alpha$ is a natural transformation $G'\circ F\rightarrow F'\circ
G$. Notice that $\alpha$ is not supposed to be an isomorphism
of functors, any functor morphism is acceptable. Similarly for the diagram
$$
\xymatrix{A \ar[d]_F \ar[rr]^G && B \ar[d]^{F'} \ar[ld]_{\beta}\\
A' \ar[rr]_{G'} && B'}
$$

We will often encounter one of the following two situations:
$$
\xymatrix{A \ar[d]_F \ar[rr]^G & & B \ar[d]_{F'} \ar[rr]^H & & C
\ar[d]^{F''} & A \ar[d]_F \ar[rr]^G & & B \ar[d]_{F'}
\ar[ld]_\beta \ar[rr]^H & & C \ar[ld]_{\beta'} \ar[d]^{F''}\\
A' \ar[rr]_{G'} \ar[ru]^\alpha & & B' \ar[rr]_{H'}
\ar[ru]^{\alpha'} & & C' & A' \ar[rr]_{G'} & & B' \ar[rr]_{H'} & &
C'}
$$

There is a well known procedure of pasting natural
transformations in such cases. In order to fix notation, we will
write the relevant formulas explicitly.

\smallskip

{\bf 4.1.1. Definition.} {\it These compositions are defined as follows
    $$(\alpha'\ntcomp\alpha)_a:=\alpha'_{G(a)}\circ H'(\alpha_a),\quad
    a\in A,$$
    $$(\beta\ntcomp\beta')_a:=H'(\beta_a)\circ\beta'_{G(a)},\quad a\in
    A.$$
}

\smallskip

Direct computation easily shows that both compositions are indeed
natural transformations $H'\circ G'\circ F\rightarrow F''\circ
H\circ G$ and $F''\circ H\circ G\rightarrow H'\circ G'\circ F$
respectively. Moreover, the composition $\ntcomp$
 is associative.

\medskip

{\bf 4.2. Strict pseudo-operads in categories.} Let $\Cat$ be a
small category whose objects are some small categories and
morphisms are functors. We will assume that it is sufficiently
rich so that all the following constructions make sense. In
particular, we have a symmetric monoidal structure $\times$ given
by the direct product and a choice of a category $\E$ with one
object and one morphism (identity). Thus we can consider
pseudo-operads in $\Cat$.

\smallskip

Although our goal is to prove existence of $\cohom$ for operads,
these being defined in the general way as in Section 1, the
categorical operads that we will use will be solely the classical
ones. The reason for this is that the categories of graphs
themselves, that were used in Section 1 in definition of operads,
are in fact examples of classical operads in categories: see
Section 4.9.3 below. Thus for our purposes there is no need to
consider more general categorical operads, than the classical
ones.

\smallskip

{\bf 4.2.1. Definition.} {\it A strict pseudo-operad in categories is a
classical non-symmetric pseudo-operad in $(\Cat,\times)$ (e.g.
\cite{MarShSt} Definition 1.18). We will denote the category of
strict pseudo-operads in $\Cat$ by $\PsOpSt(\Cat)$. The structure
functors of an object of $\PsOpSt(\Cat)$ will be denoted by
$\{\opcomp_{m_1,..,m_n}\}_{m_i\in\Z_{>0}}$.}

\smallskip

We could have used of course the notion of a $2$-pseudo-operad in
categories and work in the more general setting of higher operads,
but for our needs in this paper strict operads in categories will
suffice. One of the reasons for this restriction is the following
example.

\smallskip

{\bf 4.2.2. Example.} One of the most important examples of a
strict pseudo-operad in categories is the endomorphism
pseudo-operad $\Endo(\C)$ of a category $\C$. It is defined as
follows:
    $$\Endo(\C)_n:=\Fun(\C^{\times^n},\C),\quad n\in\Z_{>0},$$
where $\Fun$ stands for the category of functors. The structure
morphisms are given by compositions of functors.

\medskip

{\bf 4.3. Lax morphisms between strict pseudo-operads.} As we have
noted above we could have worked with $2$-operads in categories
instead of the usual ones, i.e. we could have incorporated natural
transformations in operadic structure. However, as the example of
the endomorphism operad of a category shows, it is enough for many
purposes, and in particular for ours, to consider only classical
operads. Yet when we start organizing these operads in categories
we have to take into account the natural transformations, that we
have omitted before.

\smallskip

We have already the category $\PsOpSt(\Cat)$, whose morphisms are
sequences of functors, that commute on the nose with the structure
functors of pseudo-operads. This rarely happens. In most cases we
have a natural transformation making these diagrams commutative.
Later we will see that it is these natural transformations that
define the multiplication for the triples that we will consider.
So we need to enlarge $\PsOpSt(\Cat)$ to include not only strict
but also lax morphisms of operads. Lax morphisms were invented a
long time ago (see e.g. \cite{KS} Section 3). For completeness we
reproduce explicit definitions here. Our treatment will deviate
from the classical one only when we will consider relative version
of the construction (i.e. categories of lax morphisms with
constant codomain) and more importantly when we introduce the
notion of an equivariant operad in categories. Both of these
constructions are specifically tailored for treatment of triples.

\smallskip

{\bf 4.3.1. Definition.} {\it Let $\operad$, $\operad'$ be two strict
pseudo-operads in categories. A lax morphism from $\operad$ to
$\operad'$ is a sequence of functors
$\{\opmorf_n:\operad_n\rightarrow\operad'_n\}_{n\in\Z_{>0}}$ and a
sequence of natural transformations $\{\opmort_{m_1,..,m_n}\}$,
making the following diagram commutative:
$$
\xymatrix {\operad_n\times
\operad_{m_1}\times..\times\operad_{m_n}
\ar[d]_{\opmorf^{\times^{n+1}}} \ar[rr]^{\opcomp_{m_1,..,m_n}} & & \operad_{m_1+..+m_n} \ar[d]^\opmorf\\
\operad'_n\times \operad'_{m_1}\times..\times\operad'_{m_n}
\ar[ru]^{\opmort_{m_1,..,m_n}} \ar[rr]^{\opcomp'_{m_1,..,m_n}} & &
\operad'_{m_1+..+m_n}}
$$

\smallskip

Given two lax morphisms
$(\{\opmorf_n\},\{\opmort_{m_1,..,m_n}\}):\operad\rightarrow\operad'$
and
$(\{\opmorf'_n\},\{\opmort'_{m_1,..,m_n}\}):\operad'\rightarrow\operad''$
we define their composition to be $(\{\opmorf'_n\circ
\opmorf_n\},\{\opmort_{m_1,..,m_n}\ntcomp\opmort'_{m_1,..,m_n}\})$,
where $\ntcomp$ is composition of natural transformations, as
defined in 4.1.1.}

\smallskip

It is easy to see that pseudo-operads in categories and lax
morphisms form a category. From the the pasting theorem of [Pow],
 we know that
composition of $\opmort$'s is associative and therefore
composition of the whole morphisms is associative. There is an
identity lax morphism for every category, given by the identity
functor and the trivial natural automorphism of it.

\smallskip

However, we are interested in a subcategory of this category,
consisting of lax morphisms, that have the property of coherence.
As usual coherence means that different ways of composing natural
transformations are equal. Later we will see that these conditions
will translate into associativity properties of the multiplication
natural transformations of triples that we will construct.

\smallskip

{\bf 4.3.2. Definition.}  {\it Let $\operad$ and $\operad'$ be two
strict pseudo-operads in categories. Let
$(\opmorf,\opmort):\operad\rightarrow\operad'$ be a lax morphism
between them. We say that $(\opmorf,\opmort)$ is coherent if the
$\ntcomp$-compositions of natural transformations in the following
two diagrams are equal:
$$
\xymatrix {\operad_m\times\operad_{\underline{m}}\times
\operad_{\underline{n}} \ar[rr]^{\opcomp_{\underline{m}}\times\Id}
\ar[d]_\opmorf && \operad_n\times\operad_{\underline{n}}
\ar[d]_\opmorf
\ar[rr]^{\opcomp_{\underline{n}}} && \operad_{\beta_1+..+\beta_n} \ar[d]^\opmorf\\
\operad'_m\times\operad'_{\underline{m}}\times\operad'_{\underline{n}}
\ar[rr]_{\opcomp'_{\underline{m}}\times\Id}
\ar[ru]^{\opmort_{\underline{m}}\times\Id} &&
\operad'_n\times\operad'_{\underline{n}}
\ar[rr]_{\opcomp'_{\underline{n}}}
\ar[ru]^{\opmort_{\underline{n}}} &&
\operad'_{\beta_1+..+\beta_n}}
$$

\smallskip

$$
\xymatrix
{\operad_m\times\operad_{\underline{m}}\times\operad_{\underline{n}}
\ar[rr]^{\Id\times\opcomp_{\underline{n}}} \ar[d]_\opmorf &&
\operad_m\times\operad_{\underline{m}'} \ar[d]_\opmorf
\ar[rr]^{\opcomp_{\underline{m}'}} && \operad_{\beta_1+..+\beta_n} \ar[d]^\opmorf\\
\operad'_m\times\operad'_{\underline{m}}\times\operad'_{\underline{n}}
\ar[rr]_{\Id\times\opcomp'_{\underline{n}}}
\ar[ru]^{\Id\times\opmort_{\underline{n}}} &&
\operad'_m\times\operad'_{\underline{m}'}
\ar[rr]_{\opcomp'_{\underline{m}'}}
\ar[ru]^{\opmort_{\underline{m}'}} &&
\operad'_{\beta_1+..+\beta_n}} \eqno(4.1)
$$
where $\underline{m}:=\{\alpha_1,..,\alpha_m\}$,
$\underline{n}:=\{\beta_1,..,\beta_n\}$,
$\underline{m}':=\{\beta_1+..+\beta_{\alpha_1},..,\beta_{\alpha_1+..+\alpha_{m-1}+1}+..+\beta_n\}$,
$\alpha_1+..+\alpha_m=n$, $\alpha_i$, $\beta_i\in\Z_{>0}$, and
$\operad_{\underline{m}}:=\operad_{\alpha_1}\times..\times\operad_{\alpha_m}$,
similarly for $\operad_{\underline{n}}$,
$\operad_{\underline{m}'}$ and $\operad'$.}

\smallskip

{\bf 4.3.3. Proposition.} {\it Strict pseudo-operads in categories
and coherent lax morphisms form a subcategory of the category of
strict pseudo-operads and lax morphisms}

\smallskip

{\bf Proof.} It is clear that the identity lax morphism for any
strict pseudo-operad is coherent. We have to prove that
composition of coherent lax morphisms is again coherent.

Let $\operad\overset{(\opmorf,\opmort)}\to{\rightarrow}\operad'
\overset{(\opmorf',\opmort')}\to{\rightarrow}\operad''$ be a
sequence of coherent lax morphisms between strict pseudo-operads
in categories. We have to show that the $\ntcomp$-product of
$\opmort$ and $\opmort'$ provides a unique way of making diagrams
commutative, i.e. compositions of the natural transformations as
in diagrams (4.1) are equal. But these diagrams are the outer
rectangles of the following diagrams
$$
\xymatrix{\operad_m\times\operad_{\underline{m}}\times\operad_{\underline{n}}
\ar[d]_\opmorf \ar[rr]^{\opcomp_{\underline{m}}\times\Id} &&
\operad_n\times\operad_{\underline{n}} \ar[d]_\opmorf
\ar[rr]^{\opcomp_{\underline{n}}} && \operad_{\beta_1+..+\beta_n} \ar[d]^\opmorf\\
\operad'_m\times\operad'_{\underline{m}}\times\operad'_{\underline{n}}
\ar[d]_{\opmorf'} \ar[ru]^{\opmort_{\underline{m}}\times\Id}
\ar[rr]^{\opcomp'_{\underline{m}}\times\Id} &&
\operad'_n\times\operad'_{\underline{n}} \ar[d]_{\opmorf'}
\ar[rr]^{\opcomp'_{\underline{n}}}
\ar[ru]^{\opmort_{\underline{n}}} &&
\operad'_{\beta_1+..+\beta_n} \ar[d]^{\opmorf'}\\
\operad''_m\times\operad''_{\underline{m}}\times\operad''_{\underline{n}}
\ar[ru]^{\opmort'_{\underline{m}}\times\Id}
\ar[rr]_{\opcomp''_{\underline{m}}\times\Id} &&
\operad''_n\times\operad''_{\underline{n}}
\ar[rr]_{\opcomp''_{\underline{n}}}
\ar[ru]^{\opmort'_{\underline{n}}} &&
\operad''_{\beta_1+..+\beta_n}}
$$

\smallskip

$$
\xymatrix{\operad_m\times\operad_{\underline{m}}\times\operad_{\underline{n}}
\ar[d]_{\opmorf} \ar[rr]^{\Id\times\opcomp_{\underline{n}}} &&
\operad_m\times\operad_{\underline{m}'} \ar[d]_{\opmorf}
\ar[rr]^{\opcomp_{\underline{m}'}} && \operad_{\beta_1+..+\beta_n} \ar[d]^\opmorf\\
\operad'_m\times\operad'_{\underline{m}}\times\operad'_{\underline{n}}
\ar[d]_{\opmorf'} \ar[ru]^{\Id\times\opmort_{\underline{n}}}
\ar[rr]^{\Id\times\opcomp'_{\underline{n}}} &&
\operad'_m\times\operad'_{\underline{m}'} \ar[d]_{\opmorf'}
\ar[rr]^{\opcomp'_{\underline{m}'}}
\ar[ru]^{\opmort_{\underline{m}'}} &&
\operad'_{\beta_1+..+\beta_n} \ar[d]^{\opmorf'}\\
\operad''_m\times\operad''_{\underline{m}}\times\operad''_{\underline{n}}
\ar[ru]^{\Id\times\opmort'_{\underline{n}}}
\ar[rr]_{\Id\times\opcomp''_{\underline{n}}} &&
\operad''_m\times\operad''_{\underline{m}'}
\ar[rr]_{\opcomp''_{\underline{m}'}}
\ar[ru]^{\opmort'_{\underline{m}'}} &&
\operad''_{\beta_1+..+\beta_n}}
$$

Since $(\opmorf,\opmort)$ is coherent, compositions of natural
transformations in the first rows of these diagrams are equal,
similarly for the second rows. Therefore compositions of first the
rows and then the columns are equal. We would like to show that
compositions of first the columns and then the rows are equal as
well.

\smallskip

This follows from the pasting theorem in [Pow].

\medskip

{\bf 4.3.4. Notation.} We will denote the category of strict
pseudo-operads in categories and coherent lax morphisms by
$\PsOpLC(\Cat)$.

\medskip

{\bf 4.4. Categories of pseudo-operads over a pseudo-operad.} As
with any category, it is necessary sometimes to consider an object
$\operad$ in $\PsOpLC(\Cat)$ and all morphisms in $\PsOpLC(\Cat)$
with codomain $\operad$. It will be very important for us when we
will work with representations of categorical operads on a
category, i.e. when $\operad$ is the endomorphism operad for some
category $\C$.

\smallskip

We would like of course to organize all these morphisms into a
category -- the category of pseudo-operads over $\operad$. But
first we have to decide what shall we call a morphism between two
such morphisms. In the standard way (i.e. without presence of
$2$-morphisms) we would define a morphism to a be a commutative
(on the nose) triangle as follows
$$
\xymatrix{\operad' \ar[rd]_{(\opmorf',\opmort')}
\ar[rr]^{(\opmorf,\opmort)} &&
\operad'' \ar[ld]^{(\opmorf'',\opmort'')}\\
& \operad &}
$$

However, we have natural transformations and we should take them
into account, i.e. we should define a morphism to be a diagram
with a natural transformation making it commutative
$$
\xymatrix{\operad' \ar[d]_{\opmorf'} \ar[rr]^{\opmorf} &&
\operad'' \ar[d]^{\opmorf''}\\
\operad \ar[ru]^{\opmoro} \ar[rr]_\Id && \operad} \eqno(4.2)
$$

\smallskip

Now we have to decide how $\opmort$, $\opmort'$, $\opmort''$ and
$\opmoro$ should relate to each other. If we did not have
$\opmoro$, then the relation would have been
$\opmort'=\opmort\ntcomp\opmort''$. Having $\opmoro$ we can put
all of these natural transformation into one big diagram
$$
\xymatrix{\operad'_m\times\operad'_{\underline{m}}
\ar[dd]_{{\opmorf'}^{\times^{m+1}}}
\ar[rrd]_{\opmorf^{\times^{m+1}}}
\ar[rrrrrr]^{\opcomp'_{\underline{m}}} &&&&&& \operad'_{|m|}
\ar[dd]^{\opmorf'} \ar[lld]_{\opmorf}\\
&& \operad''_m\times\operad''_{\underline{m}}
\ar[d]_{{\opmorf''}^{\times^{m+1}}}
\ar[rr]^{\opcomp''_{\underline{m}}} \ar[ru]^{\opmort} &&
\operad''_{|m|}
\ar[d]^{\opmorf''} &&\\
\operad_m\times\operad_{\underline{m}} \ar[ru]^{\opmoro}
\ar[rr]_\Id && \operad_m\times\operad_{\underline{m}}
\ar[ru]^{\opmort''} \ar[rr]_{\opcomp_{\underline{m}}} &&
\operad_{|m|} && \operad_{|m|} \ar[ll]_\Id \ar[lu]_{\opmoro}}
\eqno(4.3)
$$
where for typographical reasons we have omitted $\opmort'$, that
makes the outer rectangle commutative.

\smallskip

We see that there are two ways to construct natural
transformations
$$\opcomp_{\underline{m}}\circ{\opmorf'}^{\times^{m+1}}\rightarrow\opmorf''\circ\opmorf\circ\opcomp'_{\underline{m}}.$$
One is $\opmort\ntcomp\opmort''\ntcomp\opmoro$ and the other is
$\opmoro\ntcomp\opmort'$. Their equality is the natural condition
of compatibility.

\smallskip

{\bf 4.4.1. Definition.} {\it The category $\PsOpP$ has coherent lax
morphisms with codomain $\operad$ as objects, and for any two such
morphisms $(\opmorf',\opmort'):\operad'\rightarrow\operad$,
$(\opmorf'',\opmort''):\operad''\rightarrow\operad$, a morphism
from the first one to the second is a pair
$((\opmorf,\opmort),\opmoro)$, where
$(\opmorf,\opmort):\operad'\rightarrow\operad''$ is a coherent lax
morphism of strict pseudo-operads, and $\opmoro$ is a natural
transformation, making diagram (4.2) commutative and satisfying
    $$\opmort\ntcomp\opmort''\ntcomp\opmoro=\opmoro\ntcomp\opmort'.$$}

\smallskip

{\bf 4.4.2. Proposition.} {\it Constructed as above $\PsOpP$ is
indeed a category.}

\smallskip

{\bf Proof.} A morphism in $\PsOpP$ consists of a coherent lax
morphism and a natural transformation. From Definition 4.1.1 and
Proposition 4.3.3 we know how to compose both types, so the
composition in $\PsOpP$ is clear.

Since by the pasting theorem composition of natural transformations
is associative (Lemma 4.1.2) and we know that this is true for
coherent lax morphisms (Proposition 4.3.3), all we have to do now
is to show that the condition, formulated in Definition 4.4.1 is
satisfied by the composition.

So let $(\opmorf',\opmort'):\operad'\rightarrow\operad$,
$(\opmorf'',\opmort''):\operad''\rightarrow\operad$ and
$(\opmorf''',\opmort'''):\operad'''\rightarrow\operad$ be three
coherent lax morphisms, with codomain $\operad$. Suppose we have
two morphisms
$((\opmorff,\opmortt),\opmoro):(\opmorf',\opmort')\rightarrow(\opmorf'',\opmort'')$,
$((\opmorff',\opmortt'),\opmoro'):(\opmorf'',\opmort'')\rightarrow(\opmorf''',\opmort''')$,
that satisfy conditions of Definition 4.4.1. We can organize
everything into one diagram as follows
$$
\xymatrix{\operad'_m\times\operad'_{\underline{m}}
\ar[ddd]_{\opmorf'} \ar[rd]^{\opmorff}
\ar[rrrrr]^{\opcomp'_{\underline{m}}}
&&&&& \operad'_{|m|} \ar[ddd]^{\opmorf'} \ar[ld]_{\opmorff}\\
& \operad''_m\times\operad''_{\underline{m}}
\ar[ru]^{\opmortt_{\underline{m}}}
\ar[rrr]^{\opcomp''_{\underline{m}}} \ar[rd]^{\opmorff'}
\ar[dd]_{\opmorf''} &&&
\operad''_{|m|} \ar[dd]_{\opmorf''} \ar[ld]_{\opmorff'} &\\
&& \operad'''_m\times\operad'''_{\underline{m}}
\ar[ru]^{\opmortt'_{\underline{m}}}
\ar[r]^{\opcomp'''_{\underline{m}}} \ar[d]_{\opmorf'''} &
\operad'''_{|m|} \ar[d]^{\opmorf'''} &&\\
\operad_m\times\operad_{\underline{m}} \ar[ru]^{\opmoro}
\ar[r]_\Id & \operad_m\times\operad_{\underline{m}}
\ar[ru]^{\opmoro'} \ar[r]_\Id &
\operad_m\times\operad_{\underline{m}}
\ar[ru]^{\opmort'''_{\underline{m}}}
\ar[r]_{\opcomp_{\underline{m}}} & \operad_{|m|} & \operad_{|m|}
\ar[l]^\Id \ar[lu]_{\opmoro'} & \operad_{|m|} \ar[l]^\Id
\ar[lu]_{\opmoro}}
$$
Since $((\opmorff,\opmortt),\opmoro)$ and
$((\opmorff',\opmortt'),\opmoro')$ satisfy conditions of
definition 4.4.1 we have the following equalities
    $$\opmortt'\ntcomp\opmort'''\ntcomp\opmoro'=\opmoro'\ntcomp\opmort'',
    \qquad\opmortt\ntcomp\opmort''\ntcomp\opmoro=\opmoro\ntcomp\opmort'.$$
Using these equalities we get the following
    $$(\opmoro'\ntcomp\opmoro)\ntcomp\opmort'=\opmoro'\ntcomp\opmortt\ntcomp\opmort''\ntcomp\opmoro=
    \opmortt\ntcomp\opmoro'\ntcomp\opmort''\ntcomp\opmoro=
    (\opmortt\ntcomp\opmortt')\ntcomp\opmort'''\ntcomp(\opmoro'\ntcomp\opmoro),$$
where the second equality is justified by the pasting
theorem ([Pow]). The combined equality is exactly
the condition as in the Definition 4.4.1 for the composition.

\smallskip

Later we will often work with representations of categorical
operads on a category, i.e. we will study lax morphisms into the
endomorphism operad. We will want to construct a functor from the
category of certain representations to the category of triples on
that category. For that we will need a notion of the category of
representations. One candidate is obviously the category of
operads over the endomorphism operad, as constructed above.
However, it will prove to be too relaxed. We will need a somewhat
more restricted notion. Namely we will consider the subcategory
with the same objects but only strict morphisms.

\smallskip

{\bf 4.4.3. Notation.} We will denote by $\PsOpPS$ the
subcategory of $\PsOpP$, consisting of the same objects as
$\PsOpP$, but for any pair of them $(\operad',\opmorf',\opmort')$
and $(\operad'',\opmorf'',\opmort'')$, a morphism
$((\opmorf,\opmort),\opmoro)$ from the first to the second is in
$\PsOpPS$ if $(\opmorf,\opmort)$ is strict, i.e. $\opmort$ is the
identity.

Since strict morphisms form a subcategory of the category of lax
morphisms, we see that $\PsOpPS$ is indeed a subcategory of
$\PsOpP$.

\medskip

{\bf 4.5. Strict operads in categories.} Until now we have
considered pseudo-operads. Now we would like to discuss also the
unital version of our constructions. Since the category $\Cat$ is
monoidal with a unit, we have the natural notion of an operad in
$\Cat$, as before we restrict our attention only to the classical
operads. Recall that $\E$ is a choice of a category with one
object and one morphism (identity).

\smallskip

{\bf 4.5.1. Definition.} {\it A strict operad in categories is a
classical non-symmetric operad (e.g. \cite{MarShSt} Section 1.2)
in the monoidal unital category $(\Cat,\times,\E)$. We will denote
the category of strict operads in categories by $\OpSt(\Cat)$.}

\smallskip

As with pseudo-operads, we would like to extend the category
$\OpSt(\Cat)$ to include lax morphism of operads. As usual that
should mean making all diagrams, that before were commutative on
the nose, commutative only up to natural transformations. First we
list the relevant diagrams from the classical definition of
operads.

\smallskip

As defined above, a strict operad in categories is a pseudo-operad
$\operad$ with a strict morphism of strict pseudo-operads
$\opunit:\E\rightarrow\operad$, where we consider $\E$ as a
pseudo-operad with $\E_n:=\varnothing$ for $n>1$, and the obvious
structure morphism. The strict morphism $\opunit$ should make the
following diagrams commutative
$$
\xymatrix{\operad_n \ar[r] \ar[rrd]_\Id & \E\times\operad_n
\ar[r]^{\opunit\times\Id} & \operad_1\times\operad_n
\ar[d]^{\opcomp_n} & \operad_n \ar[r] \ar[rrd]_\Id &
\operad_n\times\E^{\times^n}
\ar[r]^{\Id\times\opunit^{\times^n}} & \operad_n\times\operad_1^{\times^n} \ar[d]^{\opcomp_{1,..,1}}\\
&& \operad_n &&& \operad_n}
$$

\smallskip

A strict morphism between strict operads commutes with these
$\opunit$'s on the nose, i.e. we have the following commutative
diagram for a strict morphism of operads
$\opmorf:\operad\rightarrow\operad'$
$$
\xymatrix{\E \ar[d]_\Id \ar[rr]^{\opunit} && \operad
\ar[d]^{\opmorf}\\
\E \ar[rr]_{\opunit'} && \operad'}
$$
We would like our lax morphisms to do that only up to a natural
transformation, that satisfies some coherence conditions. Later we
will see that these conditions are translated into the usual unit
axioms for triples.

\smallskip

{\bf 4.5.2. Definition.} {\it Let $\opunit:\E\rightarrow\operad$ and
$\opunit':\E\rightarrow\operad'$ be two strict operads in
categories. A coherent lax morphism from $\operad$ to $\operad'$ is
a coherent lax morphism of pseudo--operads
$(\opmorf,\opmort):\operad\rightarrow\operad'$ and a natural
transformation $\opuntr$, making the following diagram commutative
$$
\xymatrix{\E \ar[rr]^\opunit \ar[d]_{\Id} && \operad \ar[d]^\opmorf\\
\E \ar[rr]_{\opunit'} \ar[ru]^\opuntr && \operad'}
$$
such that compositions of natural transformations in the following
diagrams are identities
$$
\xymatrix{\operad_n \ar[d]^{\opmorf_n} \ar[rr] &&
\E\times\operad_n \ar[d]_{\Id\times\opmorf_n}
\ar[rr]^{\opunit\times\Id} && \operad_1\times\operad_n
\ar[d]_{\opmorf}
\ar[rr]^{\opcomp_n} && \operad_n \ar[d]^{\opmorf_n}\\
\operad'_n \ar[rr] && \E\times\operad'_n
\ar[rr]_{\opunit'\times\Id} \ar[ru]^{\opuntr\times\Id} &&
\operad'_1\times\operad'_n \ar[rr]_{\opcomp_n} \ar[ru]^{\opmort}
&& \operad'_n}
$$

\smallskip

$$
\xymatrix{\operad_n \ar[d]_{\opmorf_n} \ar[rr] &&
\operad_n\times\E^{\times^n} \ar[d]_{\opmorf_n\times\Id}
\ar[rr]^{\Id\times \opunit^{\times^n}} &&
\operad_n\times\operad_1^{\times^n}
\ar[d]_\opmorf \ar[rr]^{\opcomp_{1,\dots,1}} && \operad_n \ar[d]^{\opmorf_n}\\
\operad'_n \ar[rr] && \operad'_n\times\E^{\times^n}
\ar[rr]_{\Id\times{\opunit'}^{\times^n}}
\ar[ru]_{\Id\times\opuntr^{\times^n}} &&
\operad'_n\times{\operad'}_1^{\times^n}
\ar[rr]_{\opcomp'_{1,\dots,1}} \ar[ru]^{\opmort} && \operad'_n}
$$
}
\smallskip

{\bf 4.5.3. Proposition.} {\it Strict operads in categories and
lax morphisms between them constitute a category}

\smallskip

{\bf Proof.} Let
$(\opmorf,\opmort,\opuntr):\operad\rightarrow\operad'$ and
$(\opmorf',\opmort',\opuntr'):\operad'\rightarrow\operad''$ be two
lax morphisms between strict operads in categories. We define
their composition as
$(\opmorf'\circ\opmorf,\opmort\ntcomp\opmort',\opuntr\ntcomp\opuntr')$,
where $\ntcomp$ denotes composition of natural transformations as
defined in 4.1.1.

It is clear that for any strict operad, sequence of identity
functors and identity natural transformations in place of
$\opmort$ and $\opuntr$ constitute a lax morphism, and this
morphism satisfies the conditions of identity with respect to the
composition above. It remains to show that composition of lax
morphisms is again a lax morphism and that this composition is
associative.

To prove that composition is well defined we have to show that
composition of natural transformations in the following diagrams
are identities:
$$
\xymatrix{\operad_n \ar[d]_{\opmorf_n} \ar[rr] &&
\E\times\operad_n \ar[d]_{\Id\times\opmorf_n}
\ar[rr]^{\opunit\times\Id} && \operad_1\times\operad_n
\ar[d]_{\opmorf} \ar[rr]^{\opcomp_n} && \operad_n \ar[d]^{\opmorf_n}\\
\operad'_n \ar[rr] \ar[d]_{\opmorf'_n} && \E\times\operad'_n
\ar[d]_{\Id\times\opmorf'_n} \ar[rr]^{\opunit'\times\Id}
\ar[ru]^{\opuntr\times\Id} && \operad'_1\times\operad'_n
\ar[d]_{\opmorf'} \ar[rr]^{\opcomp'_n} \ar[ru]^{\opmort} &&
\operad'_n \ar[d]^{\opmorf'_n}\\
\operad''_n \ar[rr] && \E\times\operad''_n
\ar[rr]_{\opunit''\times\Id} \ar[ru]^{\opuntr'\times\Id} &&
\operad''_1\times\operad''_n \ar[rr]_{\opcomp''_n}
\ar[ru]^{\opmort'} && \operad''_n}
$$

\smallskip

$$
\xymatrix{\operad_n \ar[d]_{\opmorf_n} \ar[rr] &&
\operad_n\times\E^{\times^n} \ar[d]_{\opmorf_n\times\Id}
\ar[rr]^{\Id\times\opunit^{\times^n}} &&
\operad_n\times\operad_1^{\times^n}
\ar[d]_\opmorf \ar[rr]^{\opcomp_{1,..,1}} && \operad_n \ar[d]^{\opmorf_n}\\
\operad'_n \ar[d]_{\opmorf'_n} \ar[rr] &&
\operad'_n\times\E^{\times^n} \ar[d]_{\opmorf'_n\times\Id}
\ar[rr]^{\Id\times{\opunit'}^{\times^n}}
\ar[ru]^{\Id\times\opuntr} &&
\operad'_n\times{\operad'}_1^{\times^n} \ar[d]_{\opmorf'}
\ar[rr]^{\opcomp'_{1,\dots,1}} \ar[ru]^{\opmort} &&
\operad'_n \ar[d]^{\opmorf'_n}\\
\operad''_n \ar[rr] && \operad''_n\times\E^{\times^n}
\ar[rr]_{\Id\times{\opunit''}^{\times^n}}
\ar[ru]^{\Id\times\opuntr'} &&
\operad''_n\times{\operad''}_1^{\times^n}
\ar[rr]_{\opcomp''_{1,\dots,1}} \ar[ru]^{\opmort'} && \operad''_n}
$$

Arguing as in the proof of Proposition 4.3.3 we see that indeed
compositions of these natural transformations are identities and
hence composition of lax morphisms between strict operads in
categories is well defined.

Since the composition of natural transformations is
associative, in particular this is true for $\opuntr$'s, hence
strict operads in categories and lax morphisms between them indeed
form a category.

\smallskip

{\bf 4.5.4. Notation.} We will denote the category of strict
operads in categories and lax morphisms between them by
$\Operad(\Cat)$.

\medskip

{\bf 4.6. Categories of operads over an operad.} Since the
endomorphism pseudo-operad of a category is obviously an operad,
we would like to have a notion of a category of operads over an
operad, similarly to the case of pseudo-operads, that we have
considered in 4.4.

\smallskip

Let $\opunit:\E\rightarrow\operad$ be a strict operad in
categories. We want to organize morphisms in $\Operad(\Cat)$ with
codomain $(\operad,\opunit)$ into a category. Given two of them
$((\opmorf',\opmort'),\opuntr'):(\operad',\opunit')\rightarrow(\operad,\opunit)$
and
$((\opmorf'',\opmort''),\opuntr''):(\operad'',\opunit'')\rightarrow(\operad,\opunit)$
we would like to have the notion of a lax morphism from the first
to the second. If we wanted only the ones coming from
$\Operad(\Cat)$, we would have defined such a morphism as a lax
morphism of operads (Definition 4.5.2)
$((\opmorf,\opmort),\opuntr):(\opunit',\operad')\rightarrow(\operad'',\opunit'')$,
such that the following diagram is commutative
$$
\xymatrix{\E \ar[d]_\Id \ar[rrr]^{\opunit'} &&& \operad'
\ar[ld]^{\opmorf} \ar[dd]^{\opmorf'}\\
\E \ar[d]_\Id \ar[rr]^{\opunit''} \ar[ru]^{\opuntr} && \operad''
\ar[rd]^{\opmorf''} &\\
\E \ar[rrr]_{\opunit} \ar[ru]^{\opuntr''} &&& \operad}
$$
and we have the equality $\opuntr\ntcomp\opuntr''=\opuntr'$.

\smallskip

However, as in the case of pseudo-operads, we do not want an
equality $\opmorf''\circ\opmorf=\opmorf'$, but we usually have a
natural transformation $\opmoro$, making the following diagram
commutative
$$
\xymatrix{\operad' \ar[rr]^{\opmorf} \ar[d]_{\opmorf'} &&
\operad''
\ar[d]^{\opmorf''}\\
\operad \ar[ru]^{\opmoro} \ar[rr]_\Id && \operad}
$$

\smallskip

Obviously this natural transformation should satisfy the
conditions of a morphism in $\PsOpP$, spelled out in Definition
4.4.1. In addition it should respect, in a sense, the unital
structures $\opuntr$, $\opuntr'$ and $\opuntr''$.

\smallskip

In the following diagram
$$
\xymatrix{\E \ar[d]_\Id \ar[rrr]^{\opunit'} &&& \operad'
\ar[ld]^{\opmorf} \ar[rd]^{\opmorf'} &\\
\E \ar[d]_\Id \ar[rr]^{\opunit''} \ar[ru]^{\opuntr} && \operad''
\ar[rd]^{\opmorf''} && \operad \ar[l]_{\opmoro} \ar[ld]^\Id\\
\E \ar[ru]^{\opuntr''} \ar[rrr]_{\opunit} &&& \operad &}
$$
we see that there are two natural transformations
$\opunit\rightarrow\opmorf''\circ\opmorf\circ\opunit'$. One is
$\opuntr\ntcomp\opuntr''$ and the other is
$\opmoro\ntcomp\opuntr'$. Their equality is the natural
compatibility condition between $\opmoro$ and the unital
structures.

\smallskip

{\bf 4.6.1. Definition.} {\it We will denote by $\OpP$ the category
whose objects are coherent lax morphisms of strict operads with
codomain $(\operad,\opunit)$. Given two such morphisms
$((\opmorf',\opmort'),\opuntr'):(\operad',\opunit')\rightarrow(\operad,\opunit)$
and
$((\opmorf'',\opmort''),\opuntr''):(\operad'',\opunit'')\rightarrow(\operad,\opunit)$
a morphism from the first one to the second is a morphism
$((\opmorf,\opmort),\opuntr):(\operad',\opunit')\rightarrow(\operad'',\opunit'')$
in $\Operad(\Cat)$ and a natural transformation $\opmoro$, such that
$((\opmorf,\opmort),\opmoro)$ is morphism
$\opmorf'\rightarrow\opmorf''$ in $\PsOpP$, and in addition we
have
    $$\opmoro\ntcomp\opuntr'=\opuntr\ntcomp\opuntr''.$$}

{\bf 4.6.2. Proposition.} {\it Defined as above $\OpP$ is indeed a
category}

\smallskip

{\bf Proof.} Composition of two morphisms is inherited from
the category of pseudo-operads over a pseudo-operad. Identities
are obviously present. All we have to do is to show that the
composition of two morphisms, that satisfy conditions of the above
definition, also satisfies these conditions.

Let
$(\opmorf',\opmort',\opuntr'):(\operad',\opunit')\rightarrow(\operad,\opunit)$,
$(\opmorf'',\opmort'',\opuntr''):(\operad'',\opunit'')\rightarrow(\operad,\opunit)$,
$(\opmorf''',\opmort''',\opuntr'''):(\operad''',\opunit''')\rightarrow(\operad,\opunit)$
be three objects of $\OpP$. Let
$((\opmorff,\opmortt),\opuntrr,\opmoro)$ be a morphism from the
first to the second, and let
$((\opmorff',\opmortt'),\opuntrr',\opmoro')$ be a morphism from
the second to the third. We can organize everything into one
diagram
$$
\xymatrix{\E \ar[d]_\Id \ar[rrr]^{\opunit'} &&& \operad'
\ar[ld]_{\opmorff} \ar[rd]^{\opmorf'} &\\
\E \ar[d]_\Id \ar[ru]^{\opuntrr} \ar[rr]^{\opunit''} && \operad''
\ar[d]^{\opmorff'} \ar[rrd]^{\opmorf''} && \operad
\ar[l]_{\opmoro} \ar[d]^{\Id}\\
\E \ar[ru]^{\opuntrr'} \ar[rr]^{\opunit'''} \ar[d]_\Id &&
\operad'''
\ar[rd]^{\opmorf'''} && \operad \ar[l]^{\opmoro'} \ar[ld]^\Id\\
\E \ar[ru]^{\opuntr'''} \ar[rrr]_\opunit &&& \operad}
$$
By assumption we have
    $$\opmoro\ntcomp\opuntr'=\opuntrr\ntcomp\opuntr''\qquad\opmoro'\ntcomp\opuntr''=\opuntrr'\ntcomp\opuntr'''.$$
Using these equalities we get a sequence
    $$(\opmoro'\ntcomp\opmoro)\ntcomp\opuntr'=\opmoro'\ntcomp\opuntrr\ntcomp\opuntr''=
    \opuntrr\ntcomp\opmoro'\ntcomp\opuntr''=(\opuntrr\ntcomp\opuntrr')\ntcomp\opuntr''',$$
where in the second equality we have used the pasting theorem
of [Pow]. The composite equality is exactly
the condition, which $\opmoro'\ntcomp\opmoro$ should satisfy
according to Definition 4.6.1.

\smallskip

As it was noted in 4.4, when we will consider representations of
categorical operads on a category, we would like to consider lax
morphisms into the endomorphism operad of that category, and we
will want to organize these representations into a category, where
as morphisms we take a strict subset of morphisms in $\OpP$.

\smallskip

{\bf 4.6.3. Definition.} {\it Let $(\operad,\opunit)$ be a strict
operad in categories. We define the category $\OpPS$ as a
subcategory of $\OpP$, consisting of the same objects, but for any
morphism $((\opmorf,\opmort),\opuntr,\opmoro)$ in $\OpP$, it is
also a morphism in $\OpPS$ if $\opmort$ and $\opuntr$ are
identities.}

\medskip

{\bf 4.7. Operad-like triples as lax representations.} So far we
have considered categorical operads abstractly. In this subsection
we will work with specific operads, namely the endomorphism
operads of categories. Example 4.2.2 shows that for any category
$\C$, $\Endo(\C)$ is a pseudo-operad. Mapping $\E$ to the identity
functor on $\C$ obviously defines a structure of an operad on
$\Endo(\C)$.

\smallskip

In this Section we are interested in representations of
categorical operads, i.e. with lax morphisms
$\operad\rightarrow\Endo(\C)$. As Proposition 4.6.2 shows, such
representations form a category. We will work with them a lot, so
we introduce a special term.

\smallskip

{\bf 4.7.1. Definition.} {\it A generalized triple on a category $\C$ is
a coherent lax representation on it of a strict operad in categories
$\operad$. The category of generalized triples on $\C$ will be
denoted by $\GenTri(\C)$.}

\smallskip

To justify the term ``generalized triple'' we give the following
example, which is illustrative but inessential in our
considerations. It was considered in \cite{Ben} Section 5.4.

\smallskip

{\bf 4.7.2. Example.} Let $\operad$ be the strict operad $\E$. A
lax representation of $\E$ on a category $\C$ is simply
triple on $\C$ in the usual meaning of the term.

Indeed such representation consists first of a functor
$\opmorf:\E\rightarrow\Fun(\C,\C)$, which amounts to choosing a
functor $\triple:\C\rightarrow\C$, secondly of a natural
transformation $\opmort:\triple\circ\triple\rightarrow\triple$,
thirdly of a natural transformation
$\opuntr:\Id_\C\rightarrow\triple$, such that the conditions stated in
Definitions 4.5.2 and 4.3.2 are satisfied.

\smallskip

The condition spelled out in Definition 4.3.2 translates into
associativity of $\opmort$, i.e. into commutativity of the
following diagram
$$
\xymatrix{\triple^3 \ar[d]_{\triple(\opmort)} \ar[rr]^{\opmort} &&
\triple^2 \ar[d]^{\opmort}\\
\triple^2 \ar[rr]_{\opmort} && \triple}
$$

\smallskip

The condition stated in Definition 4.5.2 means that $\opuntr$ is
a unit for the operation $\opmort$, i.e. the following diagram is
commutative
$$
\xymatrix{\triple \ar[rrd]_{\Id} \ar[rr]^{\opuntr\times\Id} &&
\triple^2 \ar[d]^{\opmort} && \triple \ar[ll]_{\Id\times\opuntr} \ar[lld]^{\Id}\\
&& \triple}
$$

\smallskip

The main reason for the development of the theory of categorical
operads, that we have done, is the notion of an ``operad--like''
triple on a category. As we have explained in introductory Section 4.0,
we want to construct triples as colimits of a sequence of functors
where each element of the sequence represents ``operations of some
arity''.

\smallskip

However, in order for the combined object to satisfy the usual
axioms of a triple, the individual elements should behave in a
certain prescribed way with respect to colimits. So first we
describe the conditions, which these individual functors should
satisfy. Obviously we can consider colimits of any diagrams, but
we will restrict our attention only to colimits of groupoids. Most
of our results can be generalized to the case of arbitrary
diagrams.

\smallskip

Consider a functor $F:\C\rightarrow\C$, and a diagram
$D:\D\rightarrow\C$ in $C$. Since $F$ is a functor we have a
natural transformation
    $$\colfun_F:\colim(F\circ D)\rightarrow F(\colim(D)),$$
where we consider both sides as functors from $\Fun(\D,\C)$ to
$\C$. We will say that a functor $F$ commutes with colimits of
groupoids if $\colfun_F$ is an isomorphism, whenever $\D$ is a
groupoid.

\smallskip

In dealing with generalized triples we have a more general case of
functors of the type $F:\C^{\times^n}\rightarrow\C$. We would like
to extend the notion of commutativity with colimits to this case
too. There is an obvious way to do that, namely we will say that
$F$ commutes with colimits of groupoids if for every $1\leq i\leq
n$ and every $(n-1)$-tuple $\{C_1,..,C_{i-1},C_{i+1},..,C_n\}$ of
objects of $\C$ the functor $F(C_1,..,C_{i-1},-,C_{i+1},..,C_n)$
commutes with colimits of groupoids.

\smallskip

Having a functor $F:\C^{\times^n}\rightarrow\C$ that commutes with
colimits of groupoids, we will encounter situations when we have
$n$-diagrams  $\{\D_i\}_{1\leq i\leq n}$,
$\{D_i:\D_i\rightarrow\C\}$, and we will consider
$F(D_1,\dots,D_n)$. Since $F$ is a functor we have a morphism
    $$\colim(F(D_1,\dots,D_n))\rightarrow
    F(\colim(D_1),\dots,\colim(D_n)).$$
It is easy to see that this morphism factors through colimits in
each variable of $F$, which are by assumption isomorphisms.
Therefore this morphism is an isomorphism as well.

\medskip

{\it Remark.} Our construction of generalized triples was specifically tailored
for description of triples, i.e. monoids in the monoidal category
of functors. However, since we work with operads in categories we
obviously can use generalized triples to represent other objects,
for example monoidal structures.

\smallskip

Indeed, consider a non--symmetric operad in categories, generated
by one binary operation, and having isomorphisms in the category
of ternary operations, connecting the two different ways of
composing the binary operation with itself. If we demand that
these isomorphisms satisfy the usual pentagon conditions of
coherence, then a strict representation of this operad on a
category is nothing else but a coherently associative product on
this category. If we start with two generating operations and
demand coherent associativity of both and in addition certain
compatibility morphisms between their mixed compositions (these
morphisms do not have to be isomorphisms), then a representation
of such operad would be a $2$--monoidal category as described in
\cite{Va2}.

\medskip

Before we proceed with the definition of operad--like triples and
provide a way of constructing ordinary triples from them, we need a
technical preparation. We need to prove a lemma, that allows us to
combine functors, that commute with colimits of groupoids, and get
a functor commuting with such colimits as well. This will be
needed in the proof of associativity of the structure natural
transformation of the triple, that we construct from an
operad-like one.

\smallskip

{\bf 4.7.3. Lemma.} {\it Let $F$, $G$ be two functors
$\C\rightarrow\C$, that commute with colimits of groupoids. Then
$F\circ G$ commutes with such colimits as well. Moreover, let
$D':\D'\rightarrow\Fun(\C,\C)$ be a diagram of functors, commuting
with colimits of groupoids, and $\D'$ being a groupoid itself.
Then for any groupoid $\D$ and any diagram $D:\D\rightarrow\C$ we
have the following commutative diagram of natural transformations}
$$
\xymatrix{\underset{G_m\in D'(\D')}\to\colim(colim(F\circ G_m\circ
D)) \ar[d]_{\colfun_{F\circ G}} \ar[rr]^{\colfun_F} &&
F(\underset{G_m\in D'(\D')}\to\colim(\colim(G_m\circ D)))
\ar[d]^{F(\colfun_G)}\\
\underset{G_m\in D'(\D)}\to\colim(F\circ G_m(\colim(D)))
\ar[rr]_{\colfun_F} && F(\underset{G_m\in
D'(\D')}\to\colim(G_m(\colim(D))))}
$$

\smallskip

{\bf Proof.} Let $C$ be in the image of $F\circ G\circ D$. There
are two morphisms going out of $C$: one to $F(\colim(G\circ D))$
and another one to $F\circ G(\colim(D))$. Since $F$ is a functor,
$F(\colfun_G)$ completes these morphisms to a commutative
triangle. Therefore we have a factorization of the natural
transformation $\colim(F\circ G\circ D)\rightarrow F\circ
G(\colim(D))$ as $\colfun_F$ applied to $G\circ D$, followed by
$F(\colfun_G)$. Both of these are isomorphisms, therefore so is
their composition. This proves the first claim of the lemma.

The second claim is proved in a similar manner. One traces
different ways to get from a object in the image of $D$ to
$F(\underset{G_m\in D'(\D')}\to\colim(G_m(\colim(D))))$ and finds
that they are equal, due to functoriality of $F$ and $G_m$'s and
the assumption that these functors commute with colimits of
groupoids.

\smallskip

{\bf 4.7.4. Definition.} {\it An operad-like triple on a category $\C$
is a generalized triple $\operad\rightarrow\Endo(\C)$, such that each
component of $\operad$ is a groupoid and for each object of
$\operad_n$ its image in $\Fun(\C^{\times^n},\C)$ commutes with
colimits of groupoids.}

\smallskip

We have defined the category $\GenTri(\C)$ of all generalized
triples on $\C$ by utilizing all possible morphisms of operads
over an operad. With operad-like triples we want to restrict our
attention to only strict subcategories as in definition 4.6.3.

\smallskip

We will denote by $\OpTrip(\C)$ the
category whose objects are operad-like triples on $\C$ and whose
morphisms are strict morphisms over $\Endo(\C)$ as defined in
Definition 4.6.3.

\smallskip

{\bf 4.7.5. Proposition.} {\it Let
$((\opmorf,\opmort),\opuntr):\operad\rightarrow\Endo(\C)$ be an
operad-like triple on $\C$. Then if we define a functor
$\Tot(\opmorf):\C\rightarrow\C$ as follows
    $$\Tot(\opmorf)(C):=\underset{n\in\Z_{>0}}\to\coprod\underset{\opelem_n\in
    \opmorf_n(\operad_n)}\to\colim(\opelem_n(C^{\times^n})),$$
we get a triple on $\C$, with the multiplication and the unit
given by $\opmort$ and $\opuntr$ respectively. In this way we get
a functor $\Tot:\OpTrip(\C)\rightarrow\Triple(\C)$ from the
category of operad-like triples on $\C$ to the category of triples
on it.}

\smallskip

{\bf Proof.} First we give the definition of the multiplication
and the unit natural transformations for $\Tot(\opmorf)$. The
multiplication $\Tot(\opmort)$ is defined as composition of the
following sequence of natural transformations:
    $$\underset{n\in\Z_{>0}}\to\coprod\underset{\opelem_n\in\opmorf_n(\operad_n)}\to
    \colim(\opelem_n((\underset{m\in\Z_{>0}}\to\coprod\underset{\opelem'_m\in
    \opmorf_m(\operad_m)}\to\colim(\opelem'_m(C^{\times^m})))^{\times^n}))\rightarrow$$
    $$\rightarrow\underset{n,m_1,..,m_n}\to\coprod
    \underset{\opelem_n,\opelem'_{m_i}\in \opmorf(\operad)}\to\colim(\opelem_n
    (\opelem'_{m_1}(C^{\times^{m_1}}),..,\opelem'_{m_n}(C^{\times^{m_n}})))\rightarrow$$
    $$\rightarrow\underset{n,m_1,..,m_n}\to\coprod\underset{\opelem_n,\opelem'_{m_i}\in
    \opmorf(\operad)}\to\colim(\opmorf(\opelem_n\circ(\opelem'_{m_1},..,\opelem'_{m_n}))
    (C^{\times^{m_1+..+m_n}})),$$
where the first arrow is given by commutativity of $\opelem_n$
with colimits of groupoids and the second arrow is the sum of
natural transformations $\opmort_{m_1,..,m_n}$, given by the lax
representation.

The unit $\Tot(\opuntr)$ is given as composition of the following
sequence of natural transformations
    $$\Id_\C\rightarrow
    \opmorf_1(\opunel)\rightarrow\underset{n\in\Z_{>0}}\to\coprod\underset{\opelem_n\in
    \opmorf_n(\operad_n)}\to\colim(\opelem_n),$$
where $\opunel$ is the image under $\opmorf$ of the identity in
$\operad$, the first arrow is given by the unit $\opuntr$ in the
lax representation $((\opmorf,\opmort),\opuntr)$, and the second
arrow is the natural inclusion of an object of a diagram into the
colimit of the diagram.

We have to show that associativity and unit axioms hold. We know
from definition of coherent lax morphisms of operads (Definition
4.3.2) that the natural transformation $\opmort$ satisfies
associativity conditions. The multiplication for $\Tot(\opmorf)$
is given as a colimit of $\opmort$'s, using commutativity with
colimits of groupoids of individual functors in
$\opmorf(\operad)$. From lemma 4.7.3 we know that different ways
of composing $\colfun$'s for colimits of functors produce the same
result, therefore from associativity of $\opmort$ follows
associativity of $\Tot(\opmort)$. Similarly unit properties of
$\opuntr$ with respect to $\opmort$ imply the same for
$\Tot(\opuntr)$ with respect to $\Tot(\opmort)$.

It remains to show that $\Tot$ is a functor from $\OpTrip(\C)$ to
$\Triple(\C)$. Given two operad-like triples
$((\opmorf,\opmort),\opuntr):\operad\rightarrow\Endo(\C)$,
$((\opmorf',\opmort'),\opuntr):\operad'\rightarrow\Endo(\C)$ and a
morphism $(\opmorf'',\opmoro)$ from the first to the second we
have a natural transformation
$\Tot(\opmoro):\Tot(\opmorf)\rightarrow\Tot(\opmorf')$, given by
$\opmoro$. Indeed, each $\opelem\in\opmorf(\operad)$ is mapped by
$\opmoro$ to $\opmorf'(\opmorf''(\opelem))$ in
$\opmorf'(\operad')$. The latter is canonically included into
$\Tot(\opmorf')$.

Now we see that the compatibility conditions for $\opmoro$ with
$\opmort$, $\opmort''$ and with $\opuntr$, $\opuntr'$ translate
exactly to the fact that $\Tot(\opmoro)$ is a map between monoids
in the monoidal category of endofunctors on $\C$.

\medskip

{\bf 4.8. Symmetric operad--like triples.} So far we have dealt
with non-symmetric operads, and therefore with non-symmetric
operad-like triples. Now we would like to introduce action of
symmetric groups in our construction. Since the monoidal category
$(\Cat,\times,\E)$ is symmetric there is a standard notion of a
symmetric strict operad in categories.

\smallskip

{\bf 4.8.1. Definition.} A symmetric strict operad in categories
is a classical symmetric operad (e.g. \cite{MarShSt}, Section 1.2)
in the category $(\Cat,\times,\E)$, i.e. it is a strict operad
$\operad$ and for each $n\in\Z_{>0}$ an action of the symmetric
group $\Sym_n$ on the category $\operad_n$ is given. We will denote the
functor on $\operad_n$, that corresponds to an element
$\sym_n\in\Sym_n$, by the same symbol $\sym_n$. A coherent lax
morphism between two symmetric strict operads is a coherent lax
morphism between the operads (as defined in Definition 4.5.2),
such that the functors $\opmorf_n$ commute (on the nose) with the
action of symmetric groups.

\smallskip

Action of symmetric groups by means of functors provides
definition of a symmetric operad in categories, but it is not
useful for defining symmetric operad--like triples, since we need
natural transformations for that. Therefore we introduce the
notion of an equivariant symmetric operad in categories.

\smallskip

In order to do that we need  more notation. We will denote by
$\SYM_n$ the category whose objects are elements of the symmetric
group $\Sym_n$ and the set of morphisms $\hom(\sym_n,\sym'_n)$
between any two of them consists of one element:
$\sym_n^{-1}\sym'_n\in\Sym_n$. Composition is obvious. One could
call this ``a regular groupoid'' version of the symmetric group
$\Sym_n$.

\smallskip

{\bf 4.8.2. Definition.} {\it An equivariant symmetric operad in
categories is a symmetric strict operad $\operad$, such that for every
object $\oper_n\in\operad_n$ there is a functor
$\symact_n:\SYM_n\rightarrow\operad_n$, such that any object
$\sym_n\in\SYM_n$ is mapped to $\sym_n(\oper_n)$, and the
following compatibility conditions are satisfied. Let
$\oper_{m_i}\in\operad_{m_i}$ $1\leq i\leq n$. Then for every
$(n+1)$-tuple of morphisms $\sym_n\in\SYM_n$,
$\sym_{m_i}\in\SYM_{m_i}$, the composition functor
$\opcomp_{m_1,\dots,m_n}:\operad_n\times\operad_{m_1}\times\dots\times\operad_{m_n}\rightarrow\operad_{m_1+\dots+m_n}$
maps $\symact_n(\sym_n)\times\Id^{\times^n}$ to
$\symact_{m_1+\dots+m_n}(\sym_{\underline{n}})$, where the
morphism $\sym_{\underline{n}}\in\SYM_{m_1+\dots+m_n}$ is
$\sym_n$-permutation of the $n$ blocks. Also the composition
functor maps
$\Id\times\symact_{m_1}(\sym_{m_1})\times\dots\times\symact_{m_n}(\sym_{m_n})$
to
$\symact_{m_1+\dots+m_n}(\sym_{m_1}\times\dots\times\sym_{m_n})$,
where the morphism
$\sym_{m_1}\times\dots\times\sym_{m_n}\in\SYM_{m_1+\dots+m_n}$
corresponds to the product of permutations.}

\smallskip

It is clear that the compatibility conditions in the last definition
are meant to reflect the standard equivariance properties of
operads in symmetric categories. Indeed, when we will define
symmetric operad-like triples we will see that these compatibility
conditions translate into the usual equivariance.

\smallskip

{\bf 4.8.3. Example.} Let $\C$ be a category and consider the
endomorphism operad $\Endo(\C)$ of $\C$ (Example 4.2.2.). There is
one obvious symmetric structure on $\Endo(\C)$, namely for any
$\sym_n\in\Sym_n$ we define a functor
$\sym_n:\Fun(\C^{\times^n},\C)\rightarrow\Fun(\C^{\times^n},\C)$
as follows
    $$\Fun(\C^{\times^n},\C)\ni F\mapsto F\circ\sym_n,$$
where we consider $\sym_n$ as a functor
$\C^{\times^n}\rightarrow\C^{\times^n}$, permuting the variables.
It is clear that in this way we get a symmetric structure on
$\Endo(\C)$ and we will always consider endomorphism functors with
symmetric structures chosen in this way.

\smallskip

Note that in general an endomorphism operad is not an equivariant
symmetric operad. However, we do not require equivariance in the
definition of a lax morphism between two symmetric operads in
categories, and consequently we can consider lax morphisms from an
equivariant operad to one which is not. When the codomain is an
endomorphism operad as in the last example we will have a special
name for it.

\smallskip

{\bf 4.8.4. Definition.} {\it We will call a lax morphism from an
equivariant symmetric operad $\operad$ to an endomorphism operad
$\Endo(\C)$ a symmetric generalized triple. If every component of
$\operad$ is a groupoid and every functor in the image of the
generalized triple commutes with colimits of groupoids, we will
call such a generalized triple a symmetric operad-like triple. We
organize symmetric operad-like triples into a category, where
morphisms are lax morphisms as in Definition 4.7.4, and in
addition commuting (on the nose) with the symmetric structure (as
in Definition 4.8.1.) We will denote this category by
$\OpTriS(\C)$.}

\smallskip

We have proved a proposition (Proposition 4.7.7) stating that we
can get a usual triple from an operad-like one, and that this
correspondence is a functor. Since equivariant symmetric operads
differ from the non-symmetric ones by presence of an action of
symmetric groups (functors) and representations of the ``regular
symmetric groupoids'' (invertible morphisms) we see that the same
proof applies to equivariant operad-like triples as well. So we
get a functor $\Tot:\OpTriS(\C)\rightarrow\Triple(\C)$.

\smallskip

We would like to illustrate the role of equivariance in
representation of an equivariant symmetric operad on a category.
So let $\operad$ be an equivariant symmetric operad. Let $\C$ be a
category and let $((\opmorf,\opmort),\opuntr)$ be a representation
of $\operad$ on $\C$. Let $\oper_n,\{\oper_{m_i}\}_{1\leq i\leq
n}$ be elements of $\operad_n$ and $\{\operad_{m_i}\}_{1\leq i\leq
n}$ respectively. Let $\sym_n,\{\sym_{m_i}\}$ be morphisms in
$\SYM_n$ and $\{\SYM_{m_i}\}$. Writing explicitly the conditions
for $\opmort$ to be a natural transformation we get the following
commutative diagrams of natural transformations between functors
on $\C$
$$
\xymatrix{\opmorf(\oper_n)\circ(\opmorf(\oper_{m_1})\times\dots\times\opmorf(\oper_{m_n}))
\ar[d]_{\opmort_{m_1,\dots,m_n}} \ar[r] &
\opmorf(\oper_n)\circ\sigma_n\circ(\opmorf(\oper_{m_1})\times\dots\times\opmorf(\oper_{m_n}))
\ar[d]^{\opmort_{m_1,\dots,m_n}}\\
\opmorf(\opcomp(\oper_n\times\oper_{m_1}\times\dots\times\oper_{m_n}))
\ar[r] &
\opmorf(\opcomp(\sym_n(\oper_n)\times\oper_{m_1}\times\dots\times\oper_{m_n}))}
$$
where the upper horizontal arrow is
$\opmorf(\symact(\sym_n))\circ(\Id^{\times^n})$, and the lower
horizontal arrow is $\opmorf(\symact(\sym_{\underline{n}}))$.

\smallskip

$$
\xymatrix{\opmorf(\oper_n)\circ(\opmorf(\oper_{m_1})\times\dots\times\opmorf(\oper_{m_n}))
\ar[d]_{\opmort_{m_1,\dots,m_n}} \ar[r] &
\opmorf(\oper_n)\circ(\opmorf(\oper_{m_1})\circ\sym_{m_1}\times\dots\times\opmorf(\oper_{m_n})\circ\sym_{m_n})
\ar[d]^{\opmort_{m_1,\dots,m_n}}\\
\opmorf(\opcomp(\oper_n\times\oper_{m_1}\times\dots\times\oper_{m_n}))
\ar[r] &
\opmorf(\opcomp(\oper_n\times\sym_{m_1}(\oper_{m_1})\times\dots\times\sym_{m_n}(\oper_{m_n})))}
$$
where the upper horizontal arrow is
$\Id\circ(\opmorf(\symact(\sym_{m_1})\times\dots\times\symact(\sym_{m_n}))$,
and the lower horizontal arrow is
$\opmorf(\symact(\sym_{m_1}\times\dots\times\sym_{m_n}))$.

When we apply $\Tot$ to a symmetric operad-like triple we see that
these diagrams translate to the usual equivariance diagrams for
operads.

\medskip

{\bf 4.9. Example: operads as algebras over symmetric operad-like
triples.} In this subsection we would like to show that operads,
as they were defined in Section 1, can be described as algebras
over certain triples, that lie in the image of $\Tot$, as
described above. Until now we have considered operad--like triples
as lax representations of classical operads. However, operadic
constructions require working with colored operads, rather than
the classical ones. The passage to the colored context is
straightforward and we indicate the main steps below.

\smallskip

One could define colored operads in categories as colored operads
in the monoidal category $(\Cat,\times,\E)$ in the usual meaning
of the term. However, because we have $2$--morphisms in the
background, there are some minor adjustments to be made. As in the
usual case, a colored operad is different from a classical one by
a restriction on possible compositions. Before we define it we
need a technical preparation.

\smallskip

Suppose we have three categories $\C$, $\C'$, $\C''$, and each
object in all of them is given two colors from a set of colors
$\Colors$. One of the colors will be called incoming and the other
outgoing. We will say that we have a colored functor
$\C\times\C'\rightarrow\C''$ if for every two objects $\c$, $\c'$
of $\C$ and $\C'$ respectively, such that the incoming color of
$\c$ is equal to the outgoing one of $\c'$ (we will call such
pairs composable), we are given an object $\c''$ of $\C''$ whose
incoming color is that of $\c'$ and outgoing -- that of $\c$. On
morphisms such a functor should act as a usual functor, where we
allow morphisms only between composable pairs in $\C\times\C'$.

\smallskip

{\bf 4.9.1. Definition.} {\it Let $\Colors$ be a set. A strict
$\Colors$-colored categorical operad is a sequence of categories
$\{\operad_n\}_{n\in\Z_{>0}}$, for each object in $\operad_n$ a
set of incoming colors $(\color_1,\dots,\color_n)$ and an outgoing
color $\color$, and a set of colored functors
$\{\opcomp_{m_1,\dots,m_n}:\operad_n\times\operad_{m_1}\times\dots\times\operad_{m_n}\rightarrow
\operad_{m_1+\dots+\operad_{m_n}}\}$, satisfying the obvious
associativity and unit axioms.}

\smallskip

Just as in case of monochrome operads we can introduce the notion
of a symmetric colored operad in categories. For that we need an
action of symmetric groups on the components. However, since we
need to be able to permute different sets of incoming colors
differently we have to use colored symmetric groups, i.e. for
every set of incoming colors we have a copy of the symmetric
group, that acts (by functors) on all objects with the same set of
incoming colors. We omit writing explicitly the colored extension
of the usual equivariance axioms (it is straightforward but long).

\smallskip

Similarly the equivariance structure (Definition 4.8.2) can be
generalized to the colored context in a very straightforward
manner. Indeed it requires connecting objects by isomorphisms with
their images under permutation functors, such that the
compatibility conditions from Definition 4.8.2. are satisfied. We
leave writing the details explicitly to the reader.

\smallskip

Let $\C$ and $\C'$ be two categories. Suppose that objects in $\C$
have incoming and outgoing colors from a set of colors $\Colors$
and objects in $\C'$ - from $\Colors'$. Suppose we have a map
$\colorm:\Colors\rightarrow\Colors'$. Then we will say that a
functor $\opmorf:\C\rightarrow\C'$ preserves colors if an object
$\c\in\C$ with colors $(\color_{in},\color_{out})$ is mapped by
$\opmorf$ to an object $\c'\in\C'$ with colors
$(\colorm(\color_{in}),\colorm(\color_{out}))$.

\smallskip

{\bf 4.9.2. Definition.} {\it A coherent lax morphism between symmetric
strict colored operads in categories (colored by $\Colors$ and
$\Colors'$) is a set of color preserving functors (for a choice of
a map $\colorm:\Colors\rightarrow\Colors'$) and natural
transformations $((\opmorf,\opmort), \opuntr)$, satisfying the
colored versions of coherence conditions as in definitions 4.3.2
and 4.5.2.}

\smallskip

{\bf 4.9.3. Example.} The main example for us of a colored
categorical operad is the operad, produced by abstract categories
of labelled graphs, described in Section 1. Let $\Gamma$ be a
category as in Definition 1.3. Let $\Colors$ be the set of
$\Gamma$-corollas. This is our set of colors.

We define a strict $\Colors$-colored operad $\Graphs$ as follows.
We set $\Graphs_n$ to be the category, whose objects are pairs of
morphisms in $\Gamma$:
    $$\underset{\vertex\in\Vertices_\graph}\to\coprod\corolla_\vertex\rightarrow\graph\rightarrow\corolla,$$
where $\graph$ is an object of $\Gamma$, $\corolla$ is a
$\Gamma$-corolla, and the first arrow is one of the possible
atomizations of $\graph$, provided by property (iv) of Definition
1.3. We will denote such object simply by $\graph$. Note that
corollas in the direct product are ordered.

A morphism from one such object to another is a pair of morphisms
between $\graph$'s and coproducts of corollas, such that together with
the identity on $\corolla$ they make up a commutative ladder.

The coloring on each object of $\Graphs_n$ is obvious: as it is
written above the outgoing color is $\corolla$ and the incoming
colors are $\{\corolla_\vertex\}_{\vertex\in\Vertices_\graph}$.
Actions by symmetric groups are obvious as well: we just rearrange
summands in the direct sum of corollas for the atomization.

To define composition functors we use property (vi) of Definition
1.3. Suppose we have $n+1$ objects of $\Graphs$:
$\graph_1,\dots,\graph_n,\graph$, such that the incoming colors of
the latter are exactly the outgoing ones of the former $n$-tuple.
According to property (vi) we have an object $\graph'$ of
$\Gamma$, and a morphism $\graph'\rightarrow\graph$, fitting into
diagram of the type (1.3). Here we assume that a choice of a
particular $\graph'$ is made in each case. We will call this {\it
a choice of grafting}. It is clear that we can choose atomization
of $\graph'$ to be the direct sum of atomizations of $\graph_i$'s.
Doing that and taking
$\graph'\rightarrow\graph\rightarrow\corolla$ as the outgoing
color we get a composition on $\Graphs$. Colored units are chosen
in the obvious way: they are identity maps of the corollas.

Now we have to check the associativity and unit axioms. The unit
ones are obvious. Associativity conditions are obvious if we take
$Gr$ itself as $\Gamma$, indeed all we do is substituting graphs
in place of corollas, and this operation is associative. In case
of a general $\Gamma$ we make this associativity condition part of
the choice of grafting.

Actions of symmetric groups obviously satisfy the conditions of
Definition 4.8.1. We also have a natural equivariance structure
(Definition 4.8.2) on $\Graphs$. Indeed if we take a direct sum
and rearrange the summands the result is connected to the original
sum by a unique isomorphism, that gives us the representation of
$\SYM$. Here we should check that compatibility conditions from
4.8.2 are satisfied. In case $\Gamma$ is $Gr$ itself they are
obvious. In general we make them part of the requirements for the
choice of grafting. So we have a structure of an equivariant
symmetric operad on $\Graphs$.

\smallskip

Now let $(\C,\tensor)$ be a symmetric monoidal category. We would
like to define a representation of
$\grarep:\Graphs\rightarrow\Endo(\Fun(\Colors,\C))$ as follows:
given an object of $\Graphs$
    $$\underset{1\leq i\leq n}\to\coprod\corolla_i\rightarrow\graph\rightarrow\corolla,$$
and an object $\coldia:\Colors\rightarrow\C$, $\grarep(\graph)$
acts on $\coldia$ by mapping it to the functor
$\Colors\rightarrow\C$, whose value on $\corolla$ is
$\underset{1\leq i\leq n}\to\bigotimes\coldia(\corolla_i)$ and on
the rest of colors the value is the initial object of $\C$.

If we did not have non--identity morphisms on $\Colors$, then this
definition would have been obviously correct. Indeed, then a
functor from the category of colors would have been equivalent to
just a choice of objects in $\C$, and the above choice is
obviously functorial in $\coldia$. However, we have to take into
account the non--identity morphisms in $\Colors$.

\smallskip

Let $\corolla$ and $\corolla'$ be two corollas. And let
$\cormor:\corolla\rightarrow\corolla'$ be an isomorphism. Then for
any object $\underset{1\leq i\leq
n}\to\coprod\corolla_i\rightarrow\graph\rightarrow\corolla$ of
$\Graphs$ we have a new object $\underset{1\leq i\leq
n}\to\coprod\corolla_i\rightarrow\graph\rightarrow\corolla'$,
where the last arrow is the last arrow in the original object,
followed by $\cormor$. We will denote this new object of $\Graphs$
by $\cormor_*(\graph)$. In this way, given an object of $\Graphs$
we get a sort of $\Colors$--diagram of such objects (it is not
exactly a diagram because we have excluded morphisms from
$\Graphs$, that are non--identities on the corollas).

Note that for any isomorphism $\cormor$ as above the values of
$\grarep(\graph)(\coldia)$ on $\corolla$ and of
$\grarep(\cormor_*(\graph))(\coldia)$ on $\corolla'$ are the same.
Therefore, if given a $\coldia$, a $\graph\in\Graphs$ we define
for each corolla $\corolla'\in\Colors$, that is isomorphic to
$\corolla$
    $$\corolla'\mapsto\coprod\underset{1\leq i\leq n}\to\bigotimes
    \coldia(\corolla'_i),$$
where the coproduct is taken over all objects from $\Graphs$, that
are in the $\Colors$-diagram corresponding to $\graph$ as
described above, then we would get a new functor
$\Colors\rightarrow\C$. Indeed, every morphism in $\Colors$ (as
for example $\cormor$) is mapped to the identity automorphism of
$\underset{1\leq i\leq n}\to\bigotimes\coldia(\corolla_i)$.

\smallskip

So we get a representation of $\Graphs$ on $\Fun(\Colors,\C)$, and
we will denote it by $\grarep$. As it was noted above $\Graphs$ is
an equivariant symmetric operad and hence the operad-like triple
$\grarep$ is symmetric.

From Proposition 4.7.5 we conclude that there is a triple
$\Tot(\grarep)$ on $\Fun(\Colors,\C)$. This triple is exactly the
triple $\Cal{F}$ from Section 1.5.4, and the algebras over it are
the $\Gamma\C$-operads.

\medskip

{\bf 4.10. Existence of $\cohom$ for operads in algebras.} Let
$\triple$ be a Hopf-like triple (Definition 3.1.1) on a category
$\C$, that commutes with colimits of groupoids. From Lemma 3.1.3
we know that the category $\A$ of algebras over $\triple$ has
monoidal structure, and hence for any abstract category of
labelled graphs $\Gamma$ we can consider the category $\Gamma\A
OPER$ of $\Gamma\A$-operads. Now we would like to show that if
$\C$ possesses $\cohom$, so does $\Gamma\A OPER$.

\smallskip

The key to the proof is the observation that the forgetful functor
$\U:\A\rightarrow\C$ maps $\Gamma\A$-operads to
$\Gamma\C$-operads. Therefore an object of $\Gamma\A OPER$ is a
sequence (parameterized by $\Gamma$-corollas) of objects in $\C$,
such that each one of them is a $\triple$-algebra and altogether
they make up a $\Gamma$-operad in $\C$. Of course certain
compatibility conditions between these two structures should be
satisfied. This situation is just a triple version of the usual
instance of an action of a colored operad.

\smallskip

First we are going to consider sequences of objects in $\C$, that
have both of the above structures, but with the compatibility
conditions omitted. We need a very simple lemma for this, whose
proof is straightforward and we leave it to the reader.

\smallskip

{\bf 4.10.1. Lemma.} {\it Let $\triple$ and $\triple'$ be two
triples on $\C$, such that both commute with colimits of groupoids.
Then the category $\A''$ of objects in $\C$, that are
simultaneously algebras over $\triple$ and $\triple'$ is
equivalent to the category of algebras over the following triple:
    $$\triple\coprod\triple':\ \c\mapsto\coprod\triple(\dots\triple'(\dots(\c))),$$
where the coproduct is taken over all possible words of positive
length, composed of $\triple$ and $\triple'$.}

\smallskip

In our situation we have one $\triple$ (one for each
$\Gamma$-corolla) but instead of $\triple'$ we have a sequence
(parameterized by objects of $\Graphs$) of functors
$\C^{\times^n}\rightarrow\C$ (for all $n\in\Z_{>0}$). Each of
these functors commutes with colimits of groupoids so we can form
a triple out of them and $\triple$ by forming all possible
compositions and summing them up. This is an obvious
generalization of Lemma 4.10.1. We will denote the resulting
triple by $\triple\coprod\Tot(\grarep)$.

\smallskip

By construction $\triple\coprod\Tot(\grarep)$ is the direct
product, i.e. its algebras are equivalent to operads in $\C$ and
algebras over $\triple$, and these two structures being unrelated.
Now we make this direct product into an amalgamated sum. The
needed relations are provided by the Hopf-like properties of
$\triple$. Recall that for $\triple$ to be Hopf-like means that
there is a natural transformation
    $$\hopf:\triple\circ\tensor\rightarrow\tensor\circ(\triple\times\triple),$$
satisfying certain conditions, spelled out in Definition 3.1.1.

According to property b) in Definition 3.1.1. there is a definite
natural transformation
    $$\triple\circ(\tensor\circ(\Id\times\tensor))\rightarrow\tensor\circ(\Id\times\tensor)\circ(\triple^{\times^3}),$$
and similarly for all other possible iterations of the monoidal
structure. Note that functors on both sides of the last arrow are
summands in $\triple\coprod\Tot(\grarep)$ (for each
$\Gamma$-corolla), therefore there are two ways to include the
left side in the sum, i.e. we have a pair of parallel natural
transformations (for each $\Gamma$-corolla separately):
    $$\underset{n\geq
    1}\to\coprod\triple\circ\tensor^{\circ^n}\rightrightarrows\triple\coprod\Tot(\grarep).$$
These are our relations. The left side of the two arrows is just a
functor, but we have an adjunction from functors (commuting with
coproducts) to triples, therefore we have a pair of morphisms
between triples
    $$\frak{F}(\underset{n\geq
    1}\to\coprod\triple\circ\tensor^{\circ^n})\rightrightarrows\triple\coprod\Tot(\grarep),$$
where $\frak{F}$ denotes the free triple. The coequalizer (in the
category of triples) of these two morphisms is a triple whose
algebras are exactly $\Gamma$-operads in $\A$.

\smallskip

Here we should discuss existence of coequalizers in the category
of triples on $\C$. The opposite category of $\Fun(\C,\C)$ is
equivalent to $\Fun(\C^{op},\C^{op})$, and hence the question of
existence of colimits in $\Fun(\C,\C)$ is equivalent to the
question of existence of limits in $\Fun(\C^{op},\C^{op})$. The
latter can be answered by existence of colimits in $\C$ (e.g.
\cite{Bo1}, Proposition 2.15.1.) Thus if we assume that $\C$ has
coequalizers, so does $\Fun(\C,\C)$. So the question of existence
of coequalizers in the category of triples on $\C$ is the usual
question of lifting colimits from a category to a category of
algebras over a triple. General conditions for their existence are
very restrictive, so we will assume existence of the coequalizers
above as a condition imposed on $\C$ itself and on $\triple$.

\smallskip

So far we have considered algebras over a triple, that commutes
with colimits of groupoids. However, it is not always the case, as
for example the triple of associative algebras in vector spaces
does not commute with such colimits (in general it does not
commute even with coproducts). Yet often triples that do not
commute can be represented as colimits of ones, which do commute, such as
the operad--like triples.

\smallskip

We need a reformulation of the property of a triple to be
Hopf--like in the language of operad--like triples. The following
definition expresses in the operad--like triple setting the
property of an operad to be Hopf. This notion is well known, so we
will only outline the main parts of the definition and omit the
necessary coherence properties. Examples of such operad--like
triples are provided by Hopf operads.

\smallskip

{\bf 4.10.2. Definition.} {\it Let
$\opmorf:\operad\rightarrow\Endo(\C)$ be an operad-like triple on
$\C$. It is Hopf-like if for every
$\opelem_n\in\opmorf_n(\operad_n)$, $n\in\Z_{>0}$ we have a
natural transformation
    $$\hopf_n:\opelem_n\circ\tensor^{\times^n}\rightarrow\tensor\circ(\opelem_n\times\opelem_n)\circ\sym_{2n},$$
where $\sym_{2n}$ is the permutation that moves all elements in
the even places of a sequence to the end of it. This natural
transformation should satisfy coherence conditions expressing its
associativity and compatibility with the composition natural
transformation on $\opmorf(\operad)$.}

\smallskip

Now assume that we have a Hopf--like symmetric operad--like triple
$\opmorf:\operad\rightarrow\Endo(\C)$. The family $\{\hopf_n\}$
provides us with a Hopf--like structure on $\Tot(\opmorf)$. Thus the
category $\A$ of algebras over $\Tot(\opmorf)$ has a symmetric
monoidal structure. Let $\Gamma$ be an abstract category of
labelled graphs. We would like to have objects of $\Gamma\A OPER$
as algebras over a triple on $\C$.

Just as we did in case of a single $\triple$ we first consider the
coproduct $\Tot(\opmorf)\coprod\Tot(\grarep)$. Here, as before, we
take all possible compositions, but now we have to compose
functors in several variables. Also we want to take coproduct of
symmetric operad-like triples, i.e. we add up not only all
compositions, but also applications to them of permutations of
variables.

Finally, as before, we have a pair of parallel morphisms of
triples, with codomain $\Tot(\opmorf)\coprod\Tot(\grarep)$. Their
coequalizer (if it exists) is the required triple. In total we
have the following proposition.

\smallskip

{\bf 4.10.3. Proposition.} {\it Let $(\C,\tensor)$ be a symmetric
monoidal category. Let $\opmorf:\operad\rightarrow\Endo(\C)$ be a
symmetric operad-like triple on $\C$. Suppose that $\opmorf$ is
Hopf-like, and let $\Gamma$ be an abstract category of labelled
graphs. Then, if the category of triples on $\C$ has necessary
coequalizers, the category of $\Gamma$-operads in the category of
algebras over $\Tot(\opmorf)$ is equivalent to a category of
algebras over a triple on $\C$.}

\smallskip

Now using results of Section 3.1. we establish existence of
$\cohom$ for operads in algebras over operad--like triples, given
that $\cohom$ exists on $\C$.

\newpage

\centerline{\bf Appendix. Labeled graphs }

\smallskip

\centerline{\bf corresponding to various operads.}

\medskip

$\phantom{00000000000000000000000000}$This will last out a night in Russia

$\phantom{00000000000000000000000000}$When nights are longest there.

\smallskip

$\phantom{0000000000000000000000000000}${\it W.Shakespeare, Measure for measure, 2.1.132-3}

\medskip

{\bf 0. Operads, cyclic operads, modular operads.}
The graph geometry behind these structures is basically
well known, and we will only briefly repeat it.

\smallskip

{\it Operads.} Objects: disjoint unions of
directed trees with one output each.
Morphisms: (generated by) contractions and graftings
of an output to an input. If one considers
only linear directed graphs (each vertex carries one input and one output),
one gets associative algebras.

\smallskip

{\it Cyclic operads.} Objects: disjoint unions
of (unlabeled) trees. Morphisms: contractions and graftings.

\smallskip

If one adds cyclic labeling, one gets the
non--symmetric version of operads, resp. cyclic operads.

\smallskip

{\it Modular operads.} Objects: graphs of
arbitrary topology
with genus labeling. Morphisms: contractions and graftings
compatible
with labelings in the following sense.

\smallskip

Contraction of an edge having two distinct vertices
of genera $g_1$, $g_2$,  produces a new vertex of
genus $g_1+g_2$. Contraction of a loop
augments the genus of its vertex by one.
The effect of a general contraction is the result
of the composition of contraction of edges.
Grafting does not change labels.

\medskip

{\bf 1. PROPs}. Consider first the category $\Gamma_c$
whose objects are disjoint unions of
oriented corollas, and morphisms are
mergers (including isomorphisms). Any tensor functor
$(\Gamma_c ,\coprod ) \to (\Cal{G}, \otimes )$ is determined up to an
isomorphism by the following data:

\smallskip

(i) Its values  on corollas with inputs
$\{1, \dots ,n\}$ and outputs $\{1, \dots ,m\}$
($m=0$ and $n=0$ are allowed). Let such a value be
denoted $P(m,n)$.

\smallskip

(ii) Its values upon automorphisms of such corollas. This means that
each  $P(m,n)$ is endowed by commuting actions
of $\bold{S}_m$ (left) and $\bold{S}_n$ (right).

\smallskip

(iii) Its values upon merger morphisms of such
corollas which are called {\it horizontal compositions}:
$$
P(m_1,n_1)\otimes \dots \otimes P(m_r,n_r) \to
P(m_1+\dots +m_r, n_1+\dots +n_r).
\eqno(A.1)
$$

Consider now a larger category $\Gamma$
of directed graphs without oriented wheels.
This puts restrictions to morphisms compatible
with orientations. In particular, if we contract
an edge, we must simultaneously contract all edges
connecting its ends. Mergers of two vertices
connected by an oriented path also are excluded.

\smallskip

A tensor functor $\Gamma \to \Cal{G}$ then produces
data (i)--(iii) and moreover,

\smallskip

(iv) {\it Vertical compositions:} values of the functor
upon full contractions of two--vertex directed graphs
such that all inputs belong to one vertex,
all outputs to another, and edges are oriented from
inputs to outputs:
$$
P(m,n)\otimes P(n,k) \to P(m,k), \ n\ne 0.
\eqno(A.2)
$$

These data must satisfy some compatibility conditions
which can be rephrased as existence of a monoidal category
with objects $\emptyset$, $\dots$, $\{1,\dots ,n\}$, $\dots$
(as in 1.2.3)
enriched over $\Cal{G}$ in such a way that its
morphisms become $P(m,n)$ and their composition
is given by (A.2).

\smallskip

 Allowing mergers in PROPs,
we get generally big categories $\Rightarrow \sigma$
which are main building blocks of the triple
$(\Cal{F}, \mu,\eta )$ and the respective operads.
In the following three operadic structures, we
again exclude them.

\medskip

{\bf 2. Properads.} Objects: all directed graphs as above.
Morphisms: contractions and graftings.

\medskip

{\bf 3. Dioperads.} Objects: all directed graphs with
whose connected components are simply connected. Morphisms:
contractions and graftings.

\medskip

{\bf 4. $\frac{1}{2}$--PROPs.} Objects: directed graphs with
simply connected components trees such that
each edge is either unique output of its source,
or unique input of its target. Morphisms:
contractions and graftings.

\medskip

{\bf 5. Monoidal structures on the collections.}
Following [Va1], we will introduce the following
definition, working well for the categories
of directed graphs without mergers.

\smallskip

A directed graph $\tau$ is called {\it two--level} one,
if there exists a partition
of its vertices $V_{\tau}=V_{\tau}^1\coprod V_{\tau}^2$
such that

\smallskip

{\it a) Tails at $V_{\tau}^1$ are all inputs of $\tau$,
tails at $V_{\tau}^2$ are all outputs of $\tau$.

\smallskip

b) Any edge starts at  $V_{\tau}^1$ and ends at  $V_{\tau}^2$.}

\smallskip

Clearly, such a partition is unique, if it exists at all.

\smallskip

Denote by $\Rightarrow^{(2)} \sigma$ the full subcategory
of $\Rightarrow \sigma$ consisting of objects whose sources
are two--level graphs.

\smallskip

For any two collections $A^1, A^2$, define the third one by
$$
(A^2\boxtimes_c A^1)(\sigma ):=
\roman{colim}\, (\otimes_{v\in V_{\tau}^1} A^1(\tau_v))\otimes
(\otimes_{v\in V_{\tau}^2} A^2(\tau_v))
$$
where $\roman{colim}$ is taken over $\Rightarrow^{(2)} \sigma$.

\smallskip

B.~Vallette proves that this is a monoidal structure
on collections,
and that the respective operads are monoids in the resulting
monoidal category.

\smallskip

B.~Vallette treats also the case of PROPs, but here
one must restrict oneself to ``saturated'' collections.

\bigskip

\centerline{\bf References}

\medskip

[BarW] M.~Barr, Ch.~Wells. {\it Toposes, triples and theories.}
Grundlehren der mathematischen Wissenschaften 278, Springer
Verlag, 1985, 358 pp.

[Bat] M.~Batanin. {\it The Eckmann--Hilton argument, higher operads and
$E_n$--spaces.} math.CT/0207281 Sep 2003, 58 pp.

[Bek] T.~Beke {\it Sheafifiable homotopy model categories. II.}
J. Pure Appl. Algebra  164  (2001),  no. 3, 307--324.

[Ben] J.~B\'{e}nabou. {\it Introduction to bicategories} in {\it
1967  Reports of the Midwest Category Seminar}, Springer Verlag,
1--77.

[BeMa] K.~Behrend, Yu. Manin. {\it Stacks of stable maps and Gromov--Witten
invariants.} Duke Math. Journ, vol. 85, No. 1 (1996), 1--60.

[BerDW]  R.~Berger, M.~Dubois--Violette, M.~Wambst.
{\it Homogeneous algebras.} J. Algebra, 261 (2003), 172--185.
Preprint math.QA/0203035

[BerM] R.~Berger, N.~Marconnet. {\it Koszul and Gorenstein properties
for homogeneous algebras.} Algebras and representation theory,
9 (2006), 67--97.

[BerMo] C.~Berger, I.~Moerdijk. {\it Resolution of colored operads
and rectification of homotopy algebras.} math.AT/0512576
(to appear in Cont. Math., vol. in honor of Ross Street.)

[Bl] D.~Blanc. {\it New model categories from old.} Journal of
Pure and Applied Algebra 109 (1996),  37--60

[Bo1] F.~Borceux. {\it Handbook of categorical algebra 1. Basic
category theory.} Encyclopedia of Mathematics and its
applications, Cambridge University Press, 1994, 360 pp.

[Bo2] F.~Borceux. {\it Handbook of categorical algebra 2.
Categories and structures.} Encyclopedia of Mathematics and its
applications, Cambridge University Press, 1994, 360 pp.

[CaGa] J.G.~Cabello, A.R.~Garz�n. {\it Closed model structures for
algebraic models of $n$-types.} Journal of Pure and Applied
Algebra 103 (3) (1995), 287--302.

[Cra] S.E.~Crans. {\it Quillen closed model structures for
sheaves.} Journal of Pure and Applied Algebra 101, 1995, pp. 35-57.

[DeMi] P.~Deligne, J.~Milne. {\it Tannakian categories.} In:
Hodge cycles, motives and Shimura varieties, Springer LN in Math.,
900 (1982), 101--228.

[Fr] B. Fresse. {\it Koszul duality for operads
and homology of partition posets.}

[Ga] W.~L.~Gan. {\it Koszul duality for dioperads.}
Math. Res. Lett. 10:1 (2003), 109 -- 124.

[GeKa1] E.~Getzler, M.~Kapranov. {\it Cyclic operads and cyclic
homology.} In: Geometry, Topology and Physics for
Raoul Bott (ed. by S.-T.~Yau), International Press
1995, 167--201.

[GeKa2] E.~Getzler, M.~Kapranov. {\it Modular operads.}
Compositio Math., 110:1 (1998),  65--126.

[GiKa] V.~Ginzburg, M.~Kapranov. {\it Koszul duality for
operads.} Duke Math. J., 76:1 (1994), 203--272.

[GoMa]  A.~Goncharov, Yu.~Manin. {\it  Multiple zeta-motives and moduli spaces $\overline{M}_{0,n}$}. Compos. Math.
140:1 (2004), 1-14.  Preprint math.AG/0204102

[GrM] S.~Grillo, H.~Montani. {\it Twisted internal  COHOM objects
in the category of quantum spaces.} Preprint math.QA/0112233

[Hin] V.~Hinich. {\it Homological algebra of homotopy algebras.}
Communications in Algebra 25(10), 1997, pp. 3291-3323.

[KaMa] M.~Kapranov, Yu.~Manin. {\it Modules and Morita theorem for
operads.} Am. J. of Math., 123:5 (2001), 811--838.
Preprint math.QA/9906063

[KoMa] M.~Kontsevich, Yu.~Manin. {\it Gromov--Witten classes, quantum cohomology, and enumerative
geometry.} Comm. Math. Phys.,
164:3 (1994), 525--562.

[KS] G.M.~Kelly, R.~Street. {\it Review of the elements of
$2$-categories} in {\it Category Seminar} Lecture Notes in
Mathematics 420, Springer 1974, 75--103.

[LoMa]  A.~Losev, Yu.~Manin. {\it Extended modular operad.}
In: Frobenius Manifolds,
ed. by C. Hertling and M. Marcolli,  Vieweg \& Sohn Verlag,
Wiesbaden, 2004, 181--211. Preprint
math.AG/0301003

[Ma1] Yu.~Manin. {\it Some remarks on Koszul algebras and
quantum groups.} Ann. Inst. Fourier, Tome XXXVII, f. 4 (1987), 191--205.

[Ma2] Yu.~Manin. {\it Quantum groups and non--commutative geometry.}
Publ. de CRM, Universit\'e de Montr\'eal, 1988, 91 pp.

[Ma3] Yu.~Manin. {\it Topics in noncommutative geometry.}
Princeton University Press, 1991, 163 pp.

[Ma4] Yu.~Manin. {\it Notes on quantum groups and quantum de Rham complexes.}
Teoreticheskaya i Matematicheskaya Fizika 92:3 (1992),
425--450. Reprinted in {\it Selected papers of Yu.~I.~Manin,
World Scientific, Singapore 1996, 529--554.}

[Mar] M.~Markl. {\it Operads and PROPs.} Preprint math.AT/0601129.

[MarShSt] M.~Markl, St.~Shnider, J.~Stasheff.
{\it Operads in Algebra, Topology and Physics.}
Math. Surveys and Monographs, vol. 96, AMS 2002.

[MarkSh] I.~Markov, Y.~Shi. {\it Simulating quantum computation
by contracting tensor network.} Preprint quant-ph/0511069

[Mer] S.~Merkulov. {\it PROP profile of deformation quantization
and graph complexes with loops and wheels.}
Preprint math.QA/0412257

[PP] A.~Polishchuk, L.~Positselski. {\it Quadratic algebras.}
University Lecture series, No. 37, AMS 2005.

[Pow] A.~J.~Power. {\it A 2--categorical pasting theorem.}
Journ. of Algebra, 129 (1990), 439--445.

[Q] D.~Quillen {\it Homotopical Algebra.} Lecture Notes in
Mathematics, Vol. 43, Springer, Berlin, 1967.

[R] C.~Rezk. {\it Spaces of algebra structures and cohomology of
operads.} Ph.D. Thesis, Massachusetts Institute of Technology,
Cambridge, MA, 1996.

[S] J.~Spali\'nski. {\it Strong homotopy theory of cyclic sets.}
Journal of Pure and Applied Algebra 99 (1) (1995) 35--52.

[Va1] B.~Vallette. {\it A Koszul duality for PROPs.}
Preprint math.AT/0411542 (to appear in the Transactions of the AMS).

[Va2] B.~Vallette. {\it Manin's products, Koszul duality,
Loday algebras and Deligne conjecture.} Preprint math.QA/0609002

[Va3] B.~Vallette. {\it Free monoid in monoidal abelian categories.}
Preprint math.CT/0411543

[Zo] P.~Zograf. {\it Tensor networks and the enumeration
of regular subgraphs.} Preprint math.CO/0605256

\enddocument